\documentclass{article}
\usepackage{mathrsfs}
\usepackage{amsmath}
\usepackage{amssymb, amsmath}
\usepackage{graphicx}
\usepackage{fancyhdr}
\catcode`\@=11 \@addtoreset{equation}{section}

\catcode`\@=12
\usepackage{color}

\usepackage[colorlinks,
            linkcolor=blue,
            anchorcolor=blue,
            citecolor=red,
            ]{hyperref}
\newtheorem{Theorem}{Theorem}[section]
\newtheorem{Proposition}{Proposition}[section]
\newtheorem{Lemma}{Lemma}[section]
\newtheorem{Corollary}{Corollary}[section]
\newtheorem{Definition}{Definition}[section]

\newtheorem{Remark}{Remark}[section]

\newcommand{\newcom}{\newcommand}
\newcommand{\bTheorem}[1]{
\begin{Theorem} \label{T#1} }
\newcommand{\eT}{\end{Theorem}}

\newcommand{\bProposition}[1]{
\begin{Proposition} \label{P#1}}
\newcommand{\eP}{\end{Proposition}}

\newcommand{\bLemma}[1]{
\begin{Lemma} \label{L#1} }
\newcommand{\eL}{\end{Lemma}}

\newcommand{\bCorollary}[1]{
\begin{Corollary} \label{C#1} }
\newcommand{\eC}{\end{Corollary}}

\newcommand{\beq}{\begin{equation}}
\newcommand{\eeq}{\end{equation}}
\newcom{\ben}{\begin{eqnarray}}
\newcom{\een}{\end{eqnarray}}
\newcom{\beno}{\begin{eqnarray*}}
\newcom{\eeno}{\end{eqnarray*}}
\newcom{\bali}{\begin{aligned}}
\newcom{\eali}{\end{aligned}}

\newcommand{\bFormula}[1]{
\begin{equation} \label{#1}}
\newcommand{\eF}{\end{equation}}

\newcommand{\De}{\Delta}
\newcommand{\f}{\frac}

\newcommand{\p}{\partial}

\newcommand{\vr}{\varrho}
\newcommand{\vae}{a^\ep}
\newcommand{\vthe}{\vartheta^\ep}
\newcommand{\vre}{\vr^\ep}
\newcommand{\vte}{\theta^\ep}
\newcommand{\vue}{\vu^\ep}
\newcommand{\Pvue}{\mathbb{P}\vu^\ep}
\newcommand{\Qvue}{\mathbb{Q}\vu^\ep}

\newcommand{\vt}{\vartheta}
\newcommand{\vu}{\vc{u}}

\newcommand{\vv}{\mathbf{v}}

\newcommand{\tep}{\tau^\ep}
\newcommand{\sep}{\sigma^\ep}

\newcommand{\vc}[1]{{\boldsymbol #1}}

\newcommand{\Div}{{\rm div}}
\newcommand{\Grad}{\nabla}

\newcommand{\bProof}{{\bf Proof: }}

\newcommand{\ep}{\varepsilon}
\newcommand{\de}{\delta}

\font\F=msbm10 scaled 1000
\newcommand{\R}{\mbox{\F R}}

\newcommand\Cbox[2]{%
    \newbox\contentbox%
    \newbox\bkgdbox%
    \setbox\contentbox\hbox to \hsize{%
        \vtop{
            \kern\columnsep
            \hbox to \hsize{%
                \kern\columnsep%
                \advance\hsize by -2\columnsep%
                \setlength{\textwidth}{\hsize}%
                \vbox{
                    \parskip=\baselineskip
                    \parindent=0bp
                    #2
                }%
                \kern\columnsep%
            }%
            \kern\columnsep%
        }%
    }%
    \setbox\bkgdbox\vbox{
        \color{#1}
        \hrule width  \wd\contentbox %
               height \ht\contentbox %
               depth  \dp\contentbox
        \color{black}
    }%
    \wd\bkgdbox=0bp%
    \vbox{\hbox to \hsize{\box\bkgdbox\box\contentbox}}%
    \vskip\baselineskip%
}

\begin{document}


\pagestyle{fancy} \lhead{\color{blue}{The low Mach number limits of the full compressible Navier-Stokes equations}} \rhead{\emph{S.Li}}

\title{\bf The low mach number limit of global solutions to the full compressible Navier-Stokes system in critical Besov spaces with large initial data}

\author{
Sai Li     \\
 Department of Applied Mathematics, \\ Nanjing Forestry University, Nanjing, 210037, People's Republic of China \\ Email: lsmath@njfu.edu.cn
}

\maketitle

{\centerline {\bf Abstract }}
\vspace{2mm}
{We are concerned with global existence of regular solutions to full compressible Navier-Stokes equations and their asymptotic behavior when the Mach number is sufficiently small. We establish global existence in critical Besov spaces for arbitrary large initial date provided that the divergence-free component of initial velocity and the difference between initial temperature and density generate a global regular solution to incompressible Boussinesq systems. Moreover, we rigorously justify the convergence to the incompressible model as the Mach number tends to zero. The proof relies on a fine-grained analysis of the high-middle-low frequencies of density, velocity and temperature.
Our result can be seen as an improvement on Danchin and He [\emph{Math. Ann.}, 366 (2016), no. 3-4, pp. 1365-1402], including the extension from small initial data to large initial data and new convergence results which hold at the level of critical regularity.
\vspace{2mm}

{\bf Keywords: }{Full compressible Navier-Stokes equations, Low Mach number limit, Critical spaces, Large solutions, Global existence}

\vspace{2mm}

{\bf Mathematics Subject Classification:}{ 35Q30, 35B25, 76N10, 76N06  }

\tableofcontents

\section{Introduction and main result}
The motion of a compressible viscous and heat-conductive fluid is governed by the following full compressible Navier-Stokes systems:
\begin{equation}\label{equ1}
  \left\{\bali
  &\p_t\vr+\Div(\vr\vu)=0,\\
  &\p_t(\vr\vu)+\Div(\vr\vu\otimes\vu)+\Grad P=\Div\mathbb{S}(\Grad \vu ),\\
  &\p_t(\vr E)+\Div\left(\vr E\vu+P\vu\right)=\Div(\kappa\Grad \theta)+\Div(\mathbb{S}(\Grad \vu )\vu),\\
  &(\vr,\vu,\vt)|_{t=0}=(\vr_0,\vu_0,\vt_0).
  \eali
  \right.
\end{equation}

Here, the unknowns $\vr(t,x),\vu(t,x),\theta(t,x)$ represent the density of fluids, the velocity field and the temperature respectively where $t\geq0$ is the time coordinate and $x$ is the spatial coordinate. $\mathbb{S}(\Grad \vc{u} )$ is the viscous stress tensor given by $\mathbb{S}(\Grad \vc{u} )=2\mu  D\vc{u}+\lambda \Div \vc{u} \mathbb{I} $ where $Du=\f{\Grad u+(\Grad u)^\top} 2$ is the deformation tensor and viscosity coefficients  satisfy $\mu>0,2\mu+\lambda>0$. The total energy $E=e+\f1 2 |\vu|^2$ and $e=c_\nu\theta=\f R {\gamma-1}\theta$ is the internal energy, where $\gamma>1$ is the adiabatic constant and $R>0$ is a generic gas constant. $\kappa>0$ is the thermal conductivity coefficient.
We restrict our study to the case of a perfect gas where the pressure $P=R\vr\theta$ and to the case where region occupied by the fluid is the whole space $\mathbb{R}^d$ $(d\geq3)$. After rescaling, it is non restrictive to take $R=c_\nu=1$.

It is clear that \eqref{equ1} is invariant by the transform
\begin{align}\label{trans}
(\vr(t,x),\ \ \vu(t,x),\ \ \theta(t,x))\rightarrow(\vr(l^2t,lx),\ \ l\vu(l^2t,lx),\ \ l^2\theta(l^2t,lx)),\ \ l>0,
\end{align}
provided the pressure term $P$ has been changed to $l^2P$. Critical spaces for systems \eqref{equ1} are norm invariant for all
$l$ (up to an irrelevant constant) with respect to the transform \eqref{trans}.

 Danchin proved local well-posdness in \cite{D7} and global well-posedness in \cite{D6} for full compressible Navier-Stokes systems in critical Besov spaces with the initial data close to a stable constant equilibrium. See \cite{Fe,HL,MN} and the references therein for more excellent results on the well-posdness theory.

We focus on the mathematical theory of low Mach number flows. Accordingly, we introduce the following scaling (see \cite{A,LM})
\[
\vr(t,x)=\vre(\ep t,x),\ \ \vu(t,x)=\ep\vue(\ep t,x ),\ \ \theta(t,x)=\vte(\ep t,x )
\]
and
\[
\mu=\ep\mu^\ep,\lambda=\ep\lambda^\ep,\kappa=\ep\kappa^\ep,
\]
where $\ep>0$ is the Mach number.
Then target systems read as
\begin{equation}\label{equ2.1}
  \left\{\bali
    &\p_t\vre+\Div(\vre\vue)=0,\\
    &\vre(\p_t\vue+\vue\cdot\Grad\vue)-\mu^\ep\Delta\vue-(\lambda^\ep+\mu^\ep)\Grad\Div\vue+\f{\Grad(\vre\vte)}{\ep^2}=0,\\
    &\vre(\p_t\vte+\vue\cdot\Grad\vte)-\kappa^\ep\Delta\vte+\vre\vte\Div\vue=\ep^2\Big(2\mu^\ep|D\vue|^2+\lambda^\ep(\Div\vue)^2\Big)\\
    &(\vre,\vue,\vte)|_{t=0}=(\vre_0,\vue_0,\vte_0).
  \eali
  \right.
\end{equation}
We assume that coefficients $\mu^\ep,\lambda^\ep,\kappa^\ep$ are constants which obey  $\mu^\ep>0, \nu^\ep:=2\mu^\ep+\lambda^\ep>0,\kappa^\ep>0$ and for simplicity we throw away the corner marker $\ep$. We set the reference density and temperature to be 1 and consider the ill-prepared data of the form $\vre_0=1+\ep a_0$ and $\vte=1+\ep\vartheta_0$. Correspondingly, we introduce the unknowns $\vae$ and $\vthe$ which meet
\[
\vre=1+\ep\vae, \ \ \vte=1+\ep\vthe.
\]
Then $(\vae,\vue,\vthe)$  satisfy following systems:
\begin{equation}\label{equ2}
  \left\{\bali
  &\p_t\vae+\f{\Div\vue}\ep=-\Div(\vae\vue),\\
  &\p_t\vue+\vue\cdot\Grad\vue-\mathcal{L}\vue+\f{\Grad(\vae+\vthe)}\ep=I(\ep\vae)\left(\f{\Grad(\vae+\vthe)}\ep+\Grad(\vae\vthe)-\mathcal{L}\vue\right)-\Grad(\vae\vthe),\\
  &\p_t\vthe-\kappa\Delta\vthe+\f{\Div\vue}\ep=\ep(1-I(\ep\vae))\big(2\mu|D\vue|^2+\lambda|\Div\vue|^2\big)-I(\ep\vae)\kappa\Delta\vthe-\Div(\vthe\vue),\\
  &(\vae,\vue,\vthe)\mid_{t=0}=(a_0,\vu_0,\vartheta_0),
  \eali
  \right.
\end{equation}
where
\[
\mathcal{L}(\vu):=\mu\Delta {\vu} +{(\mu+\lambda) }\nabla \Div{\vu},\ \ I(a):=\f a {1+a}.
\]
In order to deal with the singular terms $\f{\Div\vue}\ep$ and $\f{\Grad(\vae+\vthe)}\ep$ in systems \eqref{equ2}, we introduce the unknowns
\[
\tau^\ep=\f{\vae+\vthe}{\sqrt{2}} \ \ \text{ and } \ \ \sigma^\ep=\f{\vthe-\vae} {\sqrt{2}}.
\]
where $\sigma^\ep$ and $\tau^\ep$ comes from the sum of the equation \eqref{equ2}$_3$ and \eqref{equ2}$_1$ and the difference between  \eqref{equ2}$_3$ and \eqref{equ2}$_1$ respectively. Accordingly, we can match $\tau^\ep$ and $\mathbb{Q}\vue$ as well as $\sigma^\ep$ and $\mathbb{P}\vue$ where $\mathbb{P}$ and $\mathbb{Q}$ are the standard projection onto the divergence free and the curl free vector fields respectively given by
\[
\mathbb{Q}:=-\Grad(-\Delta)^{-1}\Div,\ \ \mathbb{P}:=I-\mathbb{Q}.
\]
Then resulting systems read as
\begin{equation}\label{equ3}
  \left\{\bali
  &\p_t\tau^\ep+\f{\sqrt{2}}{\ep}\Div\mathbb{Q}\vue=\f{\sqrt{2}\ep}2\Big(1-I(\ep\vae)\Big)\Big(2\mu|D\vue|^2+\lambda|\Div\vue|^2\Big)-\Div(\tep\vue)\\
  &\p_t\mathbb{Q}\vue+\f{\sqrt{2}}{\ep}\Grad\tep=\mathbb{Q}\left(I(\ep\vae)\left(\f{\Grad(\vae+\vthe)}\ep+\Grad(\vae\vthe)-\mathcal{L}\vue\right)+\nu\Delta\vue-\vue\cdot\Grad\vue-\Grad(\vae\vthe)\right)
  \eali
  \right.
\end{equation}
and
\begin{equation}\label{equ4}
  \left\{\bali
  &\p_t\sep+\mathbb{P}\vue\cdot\Grad\sep-\f{\kappa}2\Delta\sep=\f{\kappa}2\Delta\tep-\Div(\sep\mathbb{Q}\vue)-\f{\sqrt{2}\kappa}2 I(\ep\vae)\Delta\vthe\\
  &\qquad\qquad\qquad\qquad\qquad\qquad\qquad\qquad+\f{\sqrt{2}\ep}2 (1-I(\ep\vae))\big(2\mu|D\vue|^2+\lambda|\Div\vue|^2\big),\\
  &\p_t\mathbb{P}\vue-\mu\Delta\mathbb{P\vue}+\mathbb{P}(\mathbb{P}\vue\cdot\Grad\mathbb{P}\vue)=\mathbb{P}\left(I(\ep\vae)\left(\f{\Grad(\vae+\vthe)}\ep+\Grad(\vae\vthe)-\mathcal{L}\vue\right)\right)\\
  &\qquad\qquad\qquad\qquad\qquad\qquad\qquad\qquad-\mathbb{P}\Big(\mathbb{Q}\vue\cdot\Grad\mathbb{P}\vue+\vue\cdot\Grad\mathbb{Q}\vue\Big).
  \eali
  \right.
\end{equation}
From the dispersive estimates for wave equations, one may expect $(\tep,\Qvue)$ to tend to zero in some spaces when the Mach number tends to zero.
We also expect $(\sep,\Pvue)$ to converge to the solution $(\Theta,\vv)$ to following Bousssinesq systems:
\begin{equation}\label{equ5}
  \left\{\bali
  &\p_t\Theta+\vv\cdot\Grad\Theta-\f{\kappa}2\Delta\Theta=0,\\
  &\p_t\vv-\mu\Delta\vv+\vv\cdot\Grad\vv+\Grad \pi=0,\ \ \Div \vv=0.
  \eali
  \right.
\end{equation}
Detailed convergence results will be given in Theorem \ref{main}. Next, we review several results on the low Mach number limit in Besov spaces that are closely related to our study and see \cite{A,DG,Fe1,Fe2,HLo,LM} for results in the framework of Sobolev spaces. For compressible Navier-Stokes systems in the whole space and periodic domain, Danchin \cite{D1,D3} investigated the question of convergence to the incompressible Navier-Stokes systems when the Mach number goes to zero and the convergence is verified for small initial date in critical spaces and large initial date in supercritical spaces. Then their works is generalized by Danchin and He \cite{D2} to $L^p$ type critical spaces and to full compressible Navier-Stokes systems but for small initial data. However, as mentioned in \cite{D1,D3}, whether convergence holds for large data with critical regularity is open. 
Until recently, Fujii \cite{F} solved the case of compressible Navier-Stokes systems in the whole space\footnote{The situation in the periodic domain is completely different because of the lack of dispersive estimates.}.
For compressible Navier-Stokes-Korteweg systems, Fujii  and Li \cite{F1} also verify the convergence to the incompressible models in some $L^p$ type critical spaces. Beyond that, there are very few results involving this issue. Thus here we develop the study of the fully compressible Navier-Stokes equations which include global existence of regular solutions in critical Besov spaces for large initial data and its asymptotic behavior when the Mach number is sufficiently small.

Our main result read as follows. Note that the definition of  norms involving $ \|\cdot\|^{l,\f{\beta_0}\ep}$ and $ \|\cdot\|^{h,\f{\beta_0}\ep}$ will be introduced in Section \ref{s2}.
\begin{Theorem}\label{main} Let
\[
2<p<\min{\{d,\f{2d}{d-2}\}},\ \ 0<\f1 q\leq\min{\left\{\f1 2-\f d 2\left(\f1 2-\f1 p\right),\f{d-1} 2\left(\f1 2-\f1 p\right)\right\}}
\]
and
\begin{align}\label{cond1}
(a_0,\vu_0,\vartheta_0)\in\left(\dot{B}_{2,1}^{\f{d} 2-1}(\mathbb{R}^d)\cap\dot{B}_{2,1}^{\f{d} 2}(\mathbb{R}^d)\right)\times\left(\dot{B}_{2,1}^{\f{d} 2-1}(\mathbb{R}^d)\right)^d\times\left(\dot{B}_{2,1}^{\f{d} 2-2}(\mathbb{R}^d)\cap\dot{B}_{2,1}^{\f{d} 2-1}(\mathbb{R}^d)\right).
\end{align}
with $\inf_{x\in\mathbb{R}^d} (1+\ep a_0(x))>0.$

If systems \eqref{equ5} supplement with initial data $\Theta_0=\f{\vartheta_0-a_0}{\sqrt{2}}$ and $\vv_0=\mathbb{P}\vu_0$ admits a global regular solution $(\Theta,\vv)$ in the class
\begin{align}\label{assum1}
(\Theta,\vv)\in C\Big([0,\infty;\dot{B}_{p,1}^{\f{d} p-1}\Big)\cap L^1\Big(0,\infty;\dot{B}_{p,1}^{\f{d} p+1}\Big)\times C\Big([0,\infty;\dot{B}_{p,1}^{\f{d} p-1}\Big)\cap L^1\Big(0,\infty;\dot{B}_{p,1}^{\f{d} p+1}\Big),
\end{align}
then there exists constants $\ep_0$ and $C$ depending on $p,q,d,\mu,\lambda,\kappa,a_0,\vu_0,\vartheta_0,\Theta,\vv$ such that for any $\ep\leq\ep_0$, systems \eqref{equ2} has a unique global solution $(\vae,\vue,\vthe)$ in the class\footnote{$\vthe \in L^1_{loc}\Big(0,\infty;\dot{B}_{2,1}^{\f{d} 2}\Big)$ means that for any $0<T<\infty$, $\|\vthe\|_{L^1\Big(0,T;\dot{B}^{\f d 2}_{2,1}\Big)}<\infty$.}
\begin{align}\label{cond7}
  \vae\in C\Big([0,\infty;\dot{B}_{2,1}^{\f{d} 2-1}&\cap\dot{B}_{2,1}^{\f{d} 2}\Big),\ \ \vue\in C\Big([0,\infty;\dot{B}_{2,1}^{\f{d} 2-1}\Big)\cap L^1\Big(0,\infty;\dot{B}_{2,1}^{\f{d} 2+1}\Big),\nonumber\\
  &\vthe\in C\Big([0,\infty;\dot{B}_{2,1}^{\f{d} 2-2}\Big)\cap L^1_{loc}\Big(0,\infty;\dot{B}_{2,1}^{\f{d} 2}\Big),
\end{align}
with $\inf_{(t,x)\in[0,\infty)\times\mathbb{R}^d} (1+\ep \vae(t,x))>0$, and the solution satisfies that for some positive constant $\beta_0$ depending on $d$ and $\nu$,
\begin{align}
   \ep\|\vae\|&^{h;\f{\beta_0}\ep}_{L^\infty\Big(0,\infty;\dot{B}^{\f d 2}_{2,1}\Big)}+\f1 \ep \|\vae\|^{h;\f{\beta_0}\ep}_{L^1\Big(0,\infty;\dot{B}^{\f d 2}_{2,1}\Big)}+\f 1 {\ep}\|\vthe\|^{h;\f{\beta_0}\ep}_{L^\infty\Big(0,\infty;\dot{B}^{\f d 2-2}_{2,1}\Big)\cap L^1\Big(0,\infty;\dot{B}^{\f d 2}_{2,1}\Big)}\nonumber\\
  & +\|(\vae,\vthe)\|^{l;\f{\beta_0}\ep}_{L^\infty\Big(0,\infty;\dot{B}^{\f d 2-1}_{2,1})\cap L^1(I;\dot{B}^{\f d 2+1}_{2,1}\Big)}+\|\vue\|_{L^\infty\Big(0,\infty;\dot{B}^{\f d 2-1}_{2,1}\Big)\cap L^1\Big(I;\dot{B}^{\f d 2+1}_{2,1}\Big)}\leq C, \label{est1}
\end{align}
\begin{align}
&\|\mathbb{P}\vue\|_{L^\infty\Big(0,\infty;\dot{B}^{\f d p-1-(\f1 2-\f1 p)}_{p,1}\Big)\cap L^1\Big(0,\infty;\dot{B}^{\f d p+1-(\f1 2-\f1 p)}_{p,1}\Big)}+\|\mathbb{Q}\vue\|_{L^2\Big(0,\infty;\dot{B}^{\f d p-(\f1 2-\f1 p)}_{p,1}\Big)}\nonumber\\
&+\|\tau^\ep\|^{l;\f{\beta_0}\ep}_{L^2\Big(0,\infty;\dot{B}^{\f d p-(\f1 2-\f1 p)}_{p,1}\Big)}+\|\sigma^\ep-\Theta\|^{l;\f{\beta_0}\ep}_{L^\infty\Big(0,\infty;\dot{B}^{\f d p-1-(\f1 2-\f1 p)}_{p,1})\cap L^2\Big(0,\infty;\dot{B}^{\f d p-(\f1 2-\f1 p)}_{p,1}\Big)}\nonumber\\
&+\ep^{-\f1 2-\f1 p}\|(\sigma^\ep-\Theta,\tau^\ep)\|^{h;\f{\beta_0}\ep}_{L^\infty\Big(0,\infty;\dot{B}^{\f d p-2}_{p,1}\Big)\cap L^1\Big(0,\infty;\dot{B}^{\f d p}_{p,1}\Big)}\leq C\ep^{\f1 2-\f 1 p}.\label{est3}
\end{align}
\begin{align}\label{est66}
 \lim_{\ep\rightarrow0}\|\mathbb{Q}\vue\|_{L^q\left(0,\infty:\dot{B}_{p,1}^{\f{d} p-1+\f2 q}\right)\cap L^{q'}\left(0,\infty:\dot{B}_{p,1}^{\f{d} p-1+\f2 {q'}}\right)}&+\ep\|a\|_{L^\infty(0,\infty;\dot{B}^{\f d p}_{p,1})}\\
 &+ \|\mathbb{P}\vue-\vv\|_{L^\infty\left(0,\infty;\dot{B}^{\f d p-1}_{p,1}\right)\cap L^1\left(0,\infty;\dot{B}^{\f d p+1}_{p,1}\right)}=0.\nonumber
\end{align}
\end{Theorem}

\begin{Remark} For our study, one of the critical regularity conditions for initial data is
\begin{align}\label{cond2}
   (a_0,\vu_0,\vartheta_0)\in \dot{B}_{2,1}^{\f{d} 2}(\mathbb{R}^d)\times\dot{B}_{2,1}^{\f{d} 2-1}(\mathbb{R}^d)\times\dot{B}_{2,1}^{\f{d} 2-2}(\mathbb{R}^d).
\end{align}
The assumption \eqref{cond1} has two additional requirements $a_0\in\dot{B}_{2,1}^{\f{d} 2-1}(\mathbb{R}^d)$ and $\vartheta_0\in {B}_{2,1}^{\f{d} 2-1}(\mathbb{R}^d)$ compared to \eqref{cond2}.
The requirement $a_0\in\dot{B}_{2,1}^{\f{d} 2-1}(\mathbb{R}^d)$ is prepared for low and middle frequencies of $\vae$, which also appears in \cite{D4,D6} for the same purpose. There are two main reasons for the requirement $\vartheta_0\in\dot{B}_{2,1}^{\f{d} 2-1}(\mathbb{R}^d)$.
The first is to improve the regularity of $\Theta$. Toward the same end, Danchin and He \cite{D2} assume $\Theta_0\in\dot{B}_{2,1}^{\f{d} 2-1}(\mathbb{R}^d)$ which allows $\Theta$ to be in the class
\[
\Theta\in  C\Big([0,\infty;\dot{B}_{2,1}^{\f{d} 2-1}\Big)\cap L^1\Big(0,\infty;\dot{B}_{2,1}^{\f{d} 2+1}\Big).
\]
Here  the requirement $\vartheta_0 \in\dot{B}_{2,1}^{\f{d} 2-1}(\mathbb{R}^d)$ together with $a_0\in\dot{B}_{2,1}^{\f{d} 2-1}(\mathbb{R}^d)$ implies $\Theta_0\in\dot{B}_{2,1}^{\f{d} 2-1}(\mathbb{R}^d)$ which also allows $\Theta$ to be in a larger class \eqref{assum1}, i.e.
\[
\Theta\in  C\Big([0,\infty;\dot{B}_{p,1}^{\f{d} p-1}\Big)\cap L^1\Big(0,\infty;\dot{B}_{p,1}^{\f{d} p+1}\Big),
\]
since the embedding $\dot{B}_{2,1}^{\f{d} 2+s}(\mathbb{R}^d)\hookrightarrow\dot{B}_{p,1}^{\f{d} p+s}(\mathbb{R}^d)$ holds for any $p>2,s\in \mathbb{R}$. The second is to establish uniform estimates, with respect to $\ep$, \eqref{est1} whose left-hand side is our energy functional which can be regard as a rescaling version of energy functional in \cite{D6}. This forces us to require that the norms of initial date is uniform bounded with respect to $\ep$, i.e.
\[
\ep\|a_0\|^{h;\f{\beta_0}\ep}_{\dot{B}^{\f d 2}_{2,1}}+\f1\ep\|\vartheta_0\|^{h;\f{\beta_0}\ep}_{\dot{B}^{\f d 2-2}_{2,1}}+\|(a_0,\vartheta_0)\|^{l;\f{\beta_0}\ep}_{\dot{B}^{\f d 2-1}_{2,1}}+\|\vu_0\|_{\dot{B}^{\f d 2-1}_{2,1}}\leq C,
\]
where the constant $C$ need not be small.
 The obstacles come from $\f1\ep\|\vartheta_0\|^{h;\f{\beta_0}\ep}_{\dot{B}^{\f d 2-2}_{2,1}}$ and $\|\vartheta_0\|^{l;\f{\beta_0}\ep}_{\dot{B}^{\f d 2-1}_{2,1}}$\footnote{If $\vartheta_0 \in\dot{B}_{2,1}^{\f{d} 2-2}(\mathbb{R}^d)$, $\|\vartheta_0\|^{l;\f{\beta_0}\ep}_{\dot{B}^{\f d 2-1}_{2,1}}\leq \f{\beta_0}\ep \|\vartheta_0\|^{l;\f{\beta_0}\ep}_{\dot{B}^{\f d 2-2}_{2,1}}$  }. Benefiting from the assumption $\vartheta_0 \in\dot{B}_{2,1}^{\f{d} 2-1}(\mathbb{R}^d)$, it holds that
\[
\f1\ep\|\vartheta_0\|^{h;\f{\beta_0}\ep}_{\dot{B}^{\f d 2-2}_{2,1}}+\|\vartheta_0\|^{l;\f{\beta_0}\ep}_{\dot{B}^{\f d 2-1}_{2,1}}\leq \f1\ep \Big(\f{\beta_0}\ep\Big)^{-1}\|\vartheta_0\|^{h;\f{\beta_0}\ep}_{\dot{B}^{\f d 2-1}_{2,1}}+\|\vartheta_0\|^{l;\f{\beta_0}\ep}_{\dot{B}^{\f d 2-1}_{2,1}}=\f1{\beta_0}\|\vartheta_0\|^{h;\f{\beta_0}\ep}_{\dot{B}^{\f d 2-1}_{2,1}}+\|\vartheta_0\|^{l;\f{\beta_0}\ep}_{\dot{B}^{\f d 2-1}_{2,1}}\leq C.
\]
\end{Remark}

\begin{Remark} In \cite{D2}, for small initial data satisfying
\begin{align}\label{cond3}
\ep\|\vae_0\|^{h;\f{\beta_0}\ep}_{\dot{B}^{\f d p}_{p,1}}+\f1\ep\|\vthe_0\|^{h;\f{\beta_0}\ep}_{\dot{B}^{\f d p-2}_{p,1}}+\|\mathbb{Q}\vue_0\|^{h;\f{\beta_0}\ep}_{\dot{B}^{\f d p-1}_{p,1}}+\|(\vae_0,\vthe_0,\mathbb{Q}\vue_0)\|^{l;\f{\beta_0}\ep}_{\dot{B}^{\f d 2-1}_{2,1}}+\|\mathbb{P}\vue_0\|_{\dot{B}^{\f d 2-1}_{2,1}}\ll1
\end{align}

with $2\leq p<d$ and $p\leq \f{2d}{d-2}$, the authors showed following global estimates
\begin{align}
   \ep\|\vae\|&^{h;\f{\beta_0}\ep}_{L^\infty\Big(0,\infty;\dot{B}^{\f d p}_{p,1}\Big)}+\f1 \ep \|\vae\|^{h;\f{\beta_0}\ep}_{L^1\Big(0,\infty;\dot{B}^{\f d p}_{p,1}\Big)}+\f 1 {\ep}\|\vthe\|^{h;\f{\beta_0}\ep}_{L^\infty\Big(0,\infty;\dot{B}^{\f d p-2}_{p,1}\Big)\cap L^1\Big(0,\infty;\dot{B}^{\f d p}_{p,1}\Big)}\nonumber\\
  &+\|\mathbb{Q}\vue\|^{h;\f{\beta_0}\ep}_{L^\infty\Big(0,\infty;\dot{B}^{\f d p-1}_{p,1}\Big)\cap L^1\Big(0,\infty;\dot{B}^{\f d p+1}_{p,1}\Big)}+\|(\vae,\vthe,\mathbb{Q}\vue)\|^{l;\f{\beta_0}\ep}_{L^\infty\Big(0,\infty;\dot{B}^{\f d 2-1}_{2,1})\cap L^1(0,\infty;\dot{B}^{\f d 2+1}_{2,1}\Big)}\nonumber\\
  &+\|\mathbb{P}\vue\|_{L^\infty\Big(0,\infty;\dot{B}^{\f d p-1}_{p,1}\Big)\cap L^1\Big(0,\infty;\dot{B}^{\f d p+1}_{p,1}\Big)}<\infty \label{est2}
\end{align}
and proved \eqref{est3}. 
It is clear that the case of $p=2$ in \eqref{est2} is estimates \eqref{est1}.  We also establish \eqref{est2} for the case of $p>2$, see \eqref{est49} and \eqref{est51}. There is no restriction on the size of the initial data, so Theorem \ref{main} can be seen as an improvement of \cite{D2} from small initial data to large initial data. The new convergence results \eqref{est66} hold at the level of critical regularity. 
\end{Remark}

The key idea in proving Theorem \ref{main} is the three-part analysis method, i.e., high-frequency, middle-frequency and low-frequency analyses which is inspired by Fujii \cite{F}. This lead to the design of quantify $M^{\ep,\alpha}_{p,q}[a,\vu,\vartheta;\tau,\sigma](I)$ (see Section \ref{s4} for more details) which will be small for suitable parameters and small time intervals $I$ . The quantify $M^{\ep,\alpha}_{p,q}[a,\vu,\vartheta;\tau,\sigma](I)$ would be a suitable replacement for the energy functional $E^\ep[\vae,\vue,\vthe](I)$ which is the left-hand side of \eqref{cond2} in some product estimates. This can transform the square term of $E^\ep[\vae,\vue,\vthe](I)$ into $E^\ep[\vae,\vue,\vthe](I)\cdot M^{\ep,\alpha}_{p,q}[a,\vu,\vartheta;\tau,\sigma](I)$. Since the energy functional $E^\ep[\vae,\vue,\vthe](I)$ need not be small in the case of large initial data, this transformation is very effective for establishing global estimates \eqref{cond2}. The quantify $M^{\ep,\alpha}_{p,q}[\vae,\vue,\vthe;\tau^\ep,\sigma^\ep](I)$ can be divided into two parts such that
\begin{align*}
M^{\ep,\alpha}_{p,q}[\vae,\vue,\vthe;\tau^\ep,\sigma^\ep](I)&\leq M^{\ep,\alpha}_{p,q}[\vae,\vue-\vv,\vthe;\tau^\ep,\sigma^\ep-\Theta](I)+\|(\Theta,\vv)\|_{L^q\Big(I;\dot{B}^{\f d p-1+\f2 q}_{p,1}\Big)\cap L^1\Big(I;\dot{B}^{\f d p+1}_{p,1}\Big)}
\end{align*}
for some $q\in[2,\infty)$. For sufficiently small time intervals $I$, we have
\begin{align}\label{cond4}
 \|(\Theta,\vv)\|_{L^q\Big(I;\dot{B}^{\f d p-1+\f2 q}_{p,1}\Big)\cap L^1\Big(I;\dot{B}^{\f d p+1}_{p,1}\Big)}\ll1
\end{align}
which cause the proof to be divided into $N$ steps where
\[
 \|(\Theta,\vv)\|_{L^q\Big(T_{n-1},T_{n};\dot{B}^{\f d p-1+\f2 q}_{p,1}\Big)}\ll1
\]
for some sequence $\{T_n\}^{N}_{n=0}$ such that $0=T_0<T_1<\cdots<T_{N-1}<T_{N}=\infty$.
The biggest difficulty when dealing with $M^{\ep,\alpha}_{p,q}[\vae,\vue-\vv,\vthe;\tau^\ep,\sigma^\ep-\Theta](I)$ is the difference between high, middle frequencies which are about the unknowns $\vae$, $\mathbb{Q}\vue, \vthe$ and low frequencies involve the unknown $\tau^\ep, \mathbb{Q}\vue, \sigma^\ep-\Theta$. This can be resolved by using \eqref{cond4}, the product estimate for truncated Besov norm in Scetion \ref{s2} and by picking sufficient suitable large $\alpha$ and  sufficient small $\ep$ such that
\[
\alpha\ep\ll1 \text{ and }\|(\Theta,\vv\|^{h;\alpha}_{L^\infty\Big(0,\infty;\dot{B}^{\f d p-1}_{p,1}\Big)\cap L^1\Big(0,\infty;\dot{B}^{\f d p+1}_{p,1}\Big)}\ll1.
\]

The rest of this paper is organized as follows. In Section \ref{s2}, we introduce some concepts and analytical tools in Besov spaces. In Section \ref{s3}, linear estimates are established. In Section 4,  we introduce the energy functional $E^\ep[a,\vu,\vartheta](I)$, quantity $M^{\ep,\alpha}_{p,q}[a,\vu,\vartheta;\tau,\sigma](I)$ and present some of their properties. In Section \ref{s5}, a priori estimates are established. In Section \ref{s6}, we complete the proof of Theorem \ref{main}. Throughout this paper, $C$ always denotes a generic positive constant which may change from line to line and $C(a_1,\cdots,a_n)$ means that $C$ only relies on $a_1,\cdots,a_n$. We denote $C^{-1}a\leq b\leq Ca$ as $a\approx b$ and  denote the conjugate of $q$ as $q'$ such that $\f1 {q'}+\f1 q=1$. For two Banach spaces $A$ and $B$,  we set $\|\cdot\|_{A\cap B}=\|\cdot\|_{A}+\|\cdot\|_{B}$.

\section{Preliminaries}\label{s2}

\subsection{Function spaces}
Let $\varphi$ be a radial, non-negative smooth function supported in $\{\xi\in\mathbb{R}^d:\f3 4\leq|\xi|\leq \f8 3\}$ and satisfying $\sum_{j\in\mathbb{Z}}\varphi(2^{-j}\xi)=1\text{ for any } \xi\in\mathbb{R}^d\backslash\{0\}$. (See e.g. \cite{BCD} for the existence of such functions). 
$\mathcal{S}(\R^d)$ and $\mathcal{S}'(\R^d)$ denote the set of Schwartz functions and tempered distributions on $\mathbb{R}^d$. For $f\in\mathcal{S}'(\R^d)$, the Littlewood-Paley decomposition is defined as follows
\[
\dot{\Delta}_{j} f  :=\mathcal{F}^{-1}(\varphi(2^{-j}|\cdot|)\mathcal{F}{f}),\ \  \dot{S}_j f := \sum\limits_{j'\leq j-1}\dot{\Delta}_{j'} f
\]
where $\mathcal{F}$ and $\mathcal{F}^{-1}$ are Fourier transform and inverse Fourier transform respectively.

Let $1\leq p,r\leq \infty $ and $s\in\mathbb{R}$ and $f\in\mathcal{S}'(\R^d)$,
\begin{align*}
 \dot{B}_{p,r}^s(\mathbb{R}^d) & :=
 \left\{f\in\mathcal{S}'(\mathbb{R}^d):\|f\|_{\dot{B}_{p,r}^s}<\infty \ \ \text{and}
\mathop {\lim }\limits_{j \to  - \infty } \left\| {{S_j}f} \right\|_{L^\infty(\mathbb{R}^d) }=0 \right\},  \\
\|f\|_{\dot{B}_{p,r}^s}  &
:=
\left(\sum\limits_{j\in\mathbb{Z}}2^{jsr}\|\dot{\Delta}_{j} f\|^{r}_{L^p(\mathbb{R}^d)}\right)^{\f1 r}.
\end{align*}
For $1\leq q,p,r\leq\infty$, $s\in\mathbb{R}$, a interval $I\subset\mathbb{R}$ and a time-dependent functions $F$ whose value is in Besov spaces, we introduce following norms:
\begin{align}
\|F\|_{L^q(I;\dot{B}_{p,r}^s)} :=\left\|\|F(t,\cdot)\|_{\dot{B}_{p,r}^s}\right\|_{L^q(I)},\ \ \|F\|_{\widetilde{L}^{q}(I;\dot{B}_{p,r}^s)}:=\left\|\left(2^{js}\|\dot{\Delta}_{j} F\|_{L^{q}(I;L^p)}\right)_{j\in\mathbb{Z}}\right\|_{\ell^r}.
\nonumber
\end{align}
The second norms was first introduced by Chemin and Lerner in \cite{CL}.
Thanks to Minkowski inequality, there holds
\begin{align}\label{minski}
\|F\|_{L^q(I;\dot{B}_{p,r}^s)}\leq\|F\|_{\widetilde{L}^{q}(I;\dot{B}_{p,r}^s)}, \ \ \text{if}\ \ r\leq q;\quad \|F\|_{L^q(I;\dot{B}_{p,r}^s)}\geq\|F\|_{\widetilde{L}^{q}(I;\dot{B}_{p,r}^s)}, \ \ \text{if}\ \ r\geq q.
\end{align}
For $0\leq\alpha \leq \beta\leq+\infty$, we define
\[
f^{l;\alpha}:=\sum\limits_{2^j<\alpha}\dot{\Delta}_{j} f,\ \ f^{m;\alpha,\beta}:=\sum\limits_{\alpha\leq2^j<\beta}\dot{\Delta}_{j} f,\ \
f^{h;\beta}:=\sum\limits_{2^j\geq\beta}\dot{\Delta}_{j} f.
\]
Next, we introduce truncated Besov norms\footnote{It is clear that $\|f\|^{l;\alpha}_{\dot{B}_{p,r}^s}\approx \|f^{l;\alpha}\|_{\dot{B}_{p,r}^s}, \|F\|^{l;\alpha}_{L^q(I;\dot{B}_{p,r}^s)}\approx \|F^{l;\alpha}\|_{L^q(I;\dot{B}_{p,r}^s)},\|F\|^{l;\alpha}_{\widetilde{L}^q(I;\dot{B}_{p,r}^s)}\approx\|F^{l;\alpha}\|_{\widetilde{L}^q(I;\dot{B}_{p,r}^s)}$. The inequalities for the middle and high frequencies are exactly similar to those for the low frequencies.}
\[
\|f\|^{h;\beta}_{\dot{B}_{p,r}^s}:=\left\|(2^{js}\|\dot{\Delta}_{j} f\|_{L^p})_{2^j\geq\beta}\right\|_{\ell^{r}},
\]
\[
\|f\|^{m;\alpha,\beta}_{\dot{B}_{p,r}^s}:=\left\|(2^{js}\|\dot{\Delta}_{j} f\|_{L^p})_{\alpha\leq2^j<\beta}\right\|_{\ell^{r}},
\]
\[
\|f\|^{l;\alpha}_{\dot{B}_{p,r}^s}:=\left\|(2^{js}\|\dot{\Delta}_{j} f\|_{L^p})_{2^j<\alpha}\right\|_{\ell^{r}},
\]
\[
\|F\|^{h;\beta}_{L^q(I;\dot{B}_{p,r}^s)}:=\left\|\|F\|^{h;\beta}_{\dot{B}_{p,r}^s}\right\|_{L^q(I)}=\left\|\|(2^{js}\|\dot{\Delta}_{j} F\|_{L^p})_{2^j\geq\beta}\|_{\ell^{r}}\right\|_{L^q(I)},
\]
\[
\|F\|^{m;\alpha,\beta}_{L^q(I;\dot{B}_{p,r}^s)}:=\left\|\|F\|^{m;\alpha,\beta}_{\dot{B}_{p,r}^s}\right\|_{L^q(I)}=\left\|\|(2^{js}\|\dot{\Delta}_{j} F\|_{L^p})_{\alpha\leq2^j<\beta}\|_{\ell^{r}}\right\|_{L^q(I)},
\]
\[
\|F\|^{l;\alpha}_{L^q(I;\dot{B}_{p,r}^s)}:=\left\|\|F\|^{l;\alpha}_{\dot{B}_{p,r}^s}\right\|_{L^q(I)}=\left\|\|(2^{js}\|\dot{\Delta}_{j} F\|_{L^p})_{2^j<\alpha}\|_{\ell^{r}}\right\|_{L^q(I)},
\]
\[
\|F\|^{h;\beta}_{\widetilde{L}^q(I;\dot{B}_{p,r}^s)}:=\left\|(2^{js}\|\dot{\Delta}_{j} F\|_{L^q(I;L^p)})_{2^j\geq\beta}\right\|_{\ell^{r}}
\]
\[
\|F\|^{m;\alpha,\beta}_{\widetilde{L}^q(I;\dot{B}_{p,r}^s)}:=\left\|(2^{js}\|\dot{\Delta}_{j} F\|_{L^q(I;L^p)})_{\alpha\leq2^j<\beta}\right\|_{\ell^{r}},
\]
\[
\|F\|^{l;\alpha}_{\widetilde{L}^q(I;\dot{B}_{p,r}^s)}
:=
\left\|(2^{js}\|\dot{\Delta}_{j} F\|_{L^q(I;L^p)})_{2^j<\alpha}\right\|_{\ell^{r}}.
\]

In the end, we introduce Bony's para-product decomposition. For $f,g\in\mathcal{S}'(\mathbb{R}^d)$,
\begin{align}
fg=T_{f}g+T_{g}f+R(g,f),
\nonumber
\end{align}
where
\begin{align}
T_{f}g=\sum\limits_{j\in\mathbb{Z}}\dot{S}_{j-2} f\dot{\Delta}_{j} g,\quad T_{g}f=\sum\limits_{j\in\mathbb{Z}}\dot{S}_{j-2} g\dot{\Delta}_{j} f,\quad
R(g,f)=\sum\limits_{|j-j'|\leq2}\dot{\Delta}_{j} g\dot{\Delta}_{j'} f.
\nonumber
\end{align}
The following almost orthogonality properties are well-known:
\beq\label{ao1}
\dot{\Delta}_{j'}\dot{\Delta}_{j} f\equiv0, \text{ if }\, |j'-j|\geq2, \ \   \Delta_{j'}( \dot{S}_{j-2}f\dot{\Delta}_{j} g)\equiv0,\text{ if }\, |j'-j|\geq3.
\eeq
and
\begin{align}\label{ao2}
\dot{\Delta}_{k}(\dot{\Delta}_{j} g\dot{\Delta}_{j'} f)\equiv0, \text{ if } k-j\geq5,\,|j'-j|\leq2.
\end{align}

\subsection{Analytical tools in Besov spaces}
In this subsection, we collect some useful lemmas and auxiliary product estimates.
\begin{Lemma}\emph{(\cite{BCD})}\label{Be1}.
Let $1\leq p,p_1,p_2,r,r_1,r_2\leq\infty$ and $s,s_1,s_2\in\mathbb{R}$. If $f\in\dot{B}_{p,r}^{s}(\mathbb{R}^d)$,
then there exists a positive constant $C=C(d,s)$ such that
\[
C^{-1}\| f\|_{\dot{B}_{p,r}^s}\leq\|\nabla f\|_{\dot{B}_{p,r}^{s-1}}\leq C\| f\|_{\dot{B}_{p,r}^s};
\]
If $f\in\dot{B}_{p,r}^{s}(\mathbb{R}^d)$, there exists a positive constant $C=C(d,s,s_1)$ such that
\[
C^{-1}\| f\|_{\dot{B}_{p,r}^s}\leq \|\Lambda^{s_1} f\|_{\dot{B}_{p,r}^{s-s_1}}\leq C\| f\|_{\dot{B}_{p,r}^s};
\]
where $\Lambda^{s_1}:=(\sqrt{-\Delta})^{s_1}$;

If $p_1\leq p_2$ and $f\in\dot{B}_{p_1,r}^{s}(\mathbb{R}^d)$, then there exists a positive constant $C=C(d,p_1,p_2)$ such that
\[
\| f\|_{\dot{B}_{p_2,r}^{s-d(\f1 {p_1}-\f1 {p_2})}}\leq C \| f\|_{\dot{B}_{{p_1},r}^s};
\]

If $\f1 p=\f1 {p_1}+\f1 {p_2},\,s=s_1+s_2,\,s_1\leq0$ and $f\in\dot{B}_{{p_1},1}^{s_1}(\mathbb{R}^d),\,g\in\dot{B}_{{p_2},1}^{s_2}(\mathbb{R}^d)$, there exists a positive constant $C=C(d,s)$ such that
\[
\|T_{f}g\|_{\dot{B}_{p,r}^s}\leq C\|f\|_{\dot{B}_{{p_1},1}^{s_1}}\|g\|_{\dot{B}_{{p_2},r}^{s_2}};
\]

If $\f1 p=\f1 {p_1}+\f1 {p_2},\,s=s_1+s_2,\,s>0,\f1 r\leq\f1 {r_1}+\f1 {r_2},\,$ and $f\in\dot{B}_{{p_1},1}^{s_1}(\mathbb{R}^d),\,g\in\dot{B}_{{p_2},1}^{s_2}(\mathbb{R}^d)$, there exists a positive constant $C=C(d,\,s_1,\,s_2)$ such that
\[
\|R(g,f)\|_{\dot{B}_{p,r}^s}\leq C\|f\|_{\dot{B}_{{p_1},r_1}^{s_1}}\|g\|_{\dot{B}_{{p_2},r_2}^{s_2}}.
\]
\end{Lemma}

\begin{Lemma}\emph{(\cite{F}).}\label{Be2} Let $G$ which vanishes at $0$ be a smooth function on an open interval $I\subset\mathbb{R}$ satisfying $(-R,R)\subset I$ with some $R>0$. If $1\leq p\leq\infty,1\leq p_1<\infty,s\in\mathbb{R}$ and $s+\f d {p_1}>0,\f1 p+\f1 {p_1}\leq1$, then for any $f\in\dot{B}_{{p},1}^{s}(\mathbb{R}^d)\cap\dot{B}_{{p_1},1}^{\f d {p_1}}(\mathbb{R}^d)$ satisfying $\|f\|_{L^\infty}<R$, there exists a positive constant $C=C(d,s,R,p,p_1,G)$ such that
\[
\|G(f)\|_{\dot{B}_{p,1}^s}\leq C
\left(
1+\|f\|_{\dot{B}_{p_1,1}^{\f d {p_1}}}
\right)
\|f\|_{\dot{B}_{p,1}^s}.
\]
\end{Lemma}

\begin{Lemma}\label{Be3} If $1\leq p_2,p_3\leq p_1\leq\infty$, $1\leq q_2,q_3\leq q_1\leq\infty$, $s\in\mathbb{R}$, $\mu>0,0<\alpha\leq\infty$ and $\vu$ is a solution to Heat equations:
\[
\p_t\vu-\mu\Delta\vu=f_1+f_2 \ \ \text{ and } \ \ \vu|_{t=0}=\vu_0,
\]
then there exist a constant $C=C(d,q_1,q_2,\mu)$ such that
\[
\|\vu\|^{l;\alpha}_{\widetilde{L}^{q_1}(I;\dot{B}_{p_1,1}^{s+\f2 {q_1}})}\leq C\Big(\|\vu_0\|^{l;\alpha}_{\dot{B}_{p_2,1}^{s+\f d {p_2}-\f d {p_1}}}+\|f_1\|^{l;\alpha}_{\widetilde{L}^{q_2}\Big(I;\dot{B}_{p_2,1}^{s+\f d {p_2}-\f d {p_1}-\f2{q'_2}}\Big)}+\|f_2\|^{l;\alpha}_{\widetilde{L}^{q_3}\Big(I;\dot{B}_{p_3,1}^{s+\f d {p_3}-\f d {p_1}-\f2{q'_3}}\Big)}\Big).
\]
\end{Lemma}

\bProof Lemma \ref{Be3}  is a direct application of Corollary 2.5 from \cite{BCD}. $\Box$

\medskip

We then show several product estimates for truncated Besov norms. For $0\leq\alpha\leq\beta\leq\infty$, $j\in\mathbb{Z}$, $\mathbf{1}_{l;\alpha}(j)$, $\mathbf{1}_{m;\alpha,\beta}(j)$ and $\mathbf{1}_{h;\beta}(j)$ are defined as
\[
\mathbf{1}_{l;\alpha}(j)=1, \text{ if } 2^j<\alpha,\ \ \mathbf{1}_{l;\alpha}(j)=0,\text{ if } 2^j\geq\alpha,
\]
\[
\mathbf{1}_{m;\alpha,\beta}(j)=1, \text{ if } \alpha\leq2^j<\beta,\ \ \mathbf{1}_{m;\alpha,\beta}(j)=0,\text{ if } 2^j<\alpha \text{ or } 2^j\geq\beta.
\]
\[
\mathbf{1}_{h;\beta}(j)=1, \text{ if } 2^j\geq\beta,\ \ \mathbf{1}_{l;\alpha}(j)=0,\text{ if } 2^j<\beta.
\]
\begin{Lemma}\label{Be5}
Let $1\leq p,p_1,p_2,r\leq\infty$, $s,s_1,s_2\in\mathbb{R}$, $s=s_1+s_2$,\ \ $\sigma=d(\f1 {p_1}+\f1 {p_2}-\f1 p)$ \text{ and } $0\leq\alpha\leq\beta\leq\infty$. If
\begin{align}\label{new1}
\f1 p\leq\f1 {p_1}+\f1 {p_2}\leq1,\ \ \f1 {p_1}\leq\f1 p,\ \ s_1\leq\f d {p_1},\ \ s_1\leq\sigma,
\end{align}
there exists a positive constant $C=C(d,p,p_1,p_2,s,s_1,s_2)$ such that
\[
\|T_{f}g\|^{m;\alpha,\beta}_{\dot{B}_{p,r}^s}\leq C\|f\|^{l;\f\beta 2}_{\dot{B}_{p_1,1}^{s_1}}\|g\|^{m;\f\alpha 4,4\beta}_{\dot{B}_{p_2,r}^{s_2+\sigma}}
\]
for $f\in\dot{B}_{p_1,1}^{s_1}(\mathbb{R}^d)$ and $g\in\dot{B}_{p_2,r}^{s_2+\sigma}(\mathbb{R}^d)$.
\end{Lemma}
\bProof  We divide the proof into two cases of $s_1<0$ and $s_1\geq0$. If $s_1<0$, from the condition $\f1 {p_1}\leq\f1 p$ and \eqref{ao1}, for $\alpha\leq2^j<\beta$, we have
\begin{align*}
2^{js}\|\dot{\Delta}_j(T_{f}g)\|_{L^p}&\leq C2^{js}\sum\limits_{|j-j'|\leq2}\|\dot{S}_{j'-2}f\|_{L^{p_1}}\|\dot{\Delta}_{j'}g\|_{L^{p^\ast_1}}\\
&\leq C
\sum\limits_{|j-j'|\leq2}2^{s_2(j-j')}\mathbf{1}_{m;\f\alpha 4,4\beta}(j')2^{j's_2}\|\dot{\Delta}_{j'}g\|_{L^{p^\ast_1}}\sum\limits_{k\leq j-1}2^{s_1(j-k)}2^{ks_1}\mathbf{1}_{l;\f\beta 2}(k)\|\dot{\Delta}_{k}f\|_{L^{p_1}},
\end{align*}
where $\f1 {p^\ast_1}=\f1 p-\f1 {p_1}$.
 Condition \eqref{new1} implies $\f1 {p^\ast_1}\leq \f1 {p_2}$, then from Young's inequalities and Bernstein inequalities 
 we have
\[
\|T_{f}g\|^{m;\alpha,\beta}_{\dot{B}_{p,r}^s}\leq C\|f\|^{l;\f\beta 2}_{\dot{B}_{p_1,1}^{s_1}}\|g\|^{m;\f\alpha 4,4\beta}_{\dot{B}_{p^\ast_1,r}^{s_2}}\leq C\|f\|^{l;\f\alpha 2}_{\dot{B}_{p_1,1}^{s_1}}\|g\|^{m;\f\alpha 4,4\beta}_{\dot{B}_{p_2,r}^{s_2+\sigma}}
\]
If $s_1\geq0$, we utilize the following decomposition
\[
\f1 p=\f1 {p_1}-\f {s_1} d+\f1 {p_2}-\f {(\sigma-s_1)} d=\f1 {p_3}+\f1 {p_4}\text{ with } \f1 {p_3}=\f1 {p_1}-\f {s_1} d \text{ and } \f1 {p_4}=\f1 {p_2}-\f {(\sigma-s_1)} d.
\]
Then $s_1\geq0$ and condition \eqref{new1} imply $0\leq\f1 {p_3}\leq \f1 {p_1},0\leq\f1 {p_4}\leq \f1 {p_2}$.
Similarly, if $\alpha\leq2^j<\beta$, we obtain
\begin{align*}
2^{js}\|\dot{\Delta}_j(T_{f}g)\|_{L^p}&\leq C2^{js}\sum\limits_{|j-j'|\leq2}\|\dot{S}_{j'-2}f\|_{L^{p_3}}\|\dot{\Delta}_{j'}g\|_{L^{p_4}}\\
&\leq C
\sum\limits_{|j-j'|\leq2}2^{s(j-j')}\mathbf{1}_{m;\f\alpha 4, 4\beta}(j')2^{j's}\|\dot{\Delta}_{j'}g\|_{L^{p_4}}\sum\limits_{k\leq j-1}\mathbf{1}_{l;\f\beta 2}(k)\|\dot{\Delta}_{k}f\|_{L^{p_3}}.
\nonumber
\end{align*}
Then, from Young's inequalities and Bernstein inequalities, 
we get
\[
\|T_{f}g\|^{m;\alpha,\beta}_{\dot{B}_{p,r}^s}\leq C\|f\|^{l;\f\beta 2}_{\dot{B}_{p_3,1}^0}\|g\|^{m;\f\alpha 4, 4\beta}_{\dot{B}_{p_4,r}^s}\leq C\|f\|^{l;\f\beta 2}_{\dot{B}_{p_1,1}^{s_1}}\|g\|^{m;\f\alpha 4, 4\beta}_{\dot{B}_{p_2,r}^{s_2+\sigma}}.\ \ \Box
\]

\begin{Lemma}\label{Be6} Let $1\leq p,p_1,p_2,p_3,p_4,r,r_1,r_2,r_3,r_4\leq\infty$, $s_1,s_2,s_3,s_4\in\mathbb{R}$, $0\leq\alpha\leq\infty$ and $s>0$. If
\[
\f1 p=\f1 {p_1}+\f1 {p_2}=\f1 {p_3}+\f1 {p_4},\ \ s=s_1+s_2=s_3+s_4, \ \ \f1 r\leq\min\left\{\f1 {r_1}+\f1 {r_2},\f1 {r_3}+\f1 {r_4}\right\},
\]
there exists a positive constant $C=C(d,s,s_1,s_2,s_3,s_4)$ such that
\begin{align}\label{Be7}
\|R(f,g)\|_{\dot{B}_{p,r}^s}\leq C\|f\|^{l;\alpha}_{\dot{B}_{p_1,r_1}^{s_1}}\|g\|^{l;4\alpha}_{\dot{B}_{p_2,r_2}^{s_2}}+C\|f\|^{h;\alpha}_{\dot{B}_{p_3,r_3}^{s_3}}\|g\|^{h;\f{\alpha} 4}_{\dot{B}_{p_4,r_4}^{s_4}}
\end{align}
and
\begin{align}\label{Be8}
\|R(f,g)\|^{h;\alpha}_{\dot{B}_{p,r}^s}\leq C \|f\|^{h;\f\alpha {16}}_{\dot{B}_{p_1,r_1}^{s_1}}\|g\|^{h;\f\alpha {64}}_{\dot{B}_{p_2,r_2}^{s_2}}
\end{align}
provided that the right-hand side is finite.
\end{Lemma}

\bProof The equality \eqref{Be7} is a direct result of Lemma 2.5 in \cite{F}. Next we show \eqref{Be8}. For $\alpha\leq 2^j$, we have
\begin{align*}
2^{js}\|\dot{\Delta}_j R(f,g)\|_{L^p}&\leq C2^{js}\sum\limits_{j-j'\leq4}\sum\limits_{|k-j'|\leq2}\|\dot{\Delta}_{j'}f\|_{L^{p_1}}\|\dot{\Delta}_k g\|_{L^{p_2}}\\
&\leq C
\sum\limits_{j-j'\leq 4}2^{s(j-j')}\mathbf{1}_{h;\f\alpha {16}}(j')2^{j's_1}\|\dot{\Delta}_{j'}f\|_{L^{p_1}}\sum\limits_{k\geq j-6}\mathbf{1}_{h;\f\alpha {64}}(k)\|\dot{\Delta}_{k}g\|_{L^{p_2}}.
\end{align*}
Then Young's inequality gives \eqref{Be8}. $\Box$

\section{Linear estimates}\label{s3}
This section is devoted to the study of the linearization of systems \eqref{equ2} which is
\begin{equation}\label{equ8}
  \left\{\bali
  &\p_t\vae+\f{\Div\vue}\ep=F,\\
  &\p_t\vue-\mathcal{L}\vue+\f{\Grad(\vae+\vthe)}\ep=G,\\
  &\p_t\vthe-\kappa\Delta\vthe+\f{\Div\vue}\ep=J,\\
  &(\vae,\vue,\vthe)\mid_{t=0}=(a_0,\vu_0,\vartheta_0).
  \eali
  \right.
\end{equation}

\begin{Lemma}\label{Le3.1} Let $1< p<\infty, s\in\mathbb{R}$ and $I\subset\mathbb{R}$ be an temporal interval with the initial time $t_0$. There exist positive constants $\beta_0=\beta_0(d,\nu)$, $C=C(d,\nu)$ and we have that for any $\ep\leq1$,
\begin{align}
   \ep&\|\vae\|^{h;\f{\beta_0}\ep}_{L^\infty\Big(0,\infty;\dot{B}^{s+1}_{p,1}\Big)}+\f1 \ep \|\vae\|^{h;\f{\beta_0}\ep}_{L^1\Big(0,\infty;\dot{B}^{s+1}_{p,1}\Big)}\nonumber\\
   &\quad+\f 1 {\ep}\|\vthe\|^{h;\f{\beta_0}\ep}_{L^\infty\Big(0,\infty;\dot{B}^{s-1}_{p,1}\Big)\cap L^1\Big(0,\infty;\dot{B}^{s+1}_{p,1}\Big)}
  +\|\mathbb{Q}\vue\|^{h;\f{\beta_0}\ep}_{L^\infty\Big(0,\infty;\dot{B}^{s}_{p,1}\Big)\cap L^1\Big(0,\infty;\dot{B}^{s+2}_{p,1}\Big)}\nonumber\\
  &\leq  C\left(\ep\|\vae\|^{h;\f{\beta_0}\ep}_{\dot{B}_{p,1}^{s+1}}(t_0)+\f1 \ep\|\vthe\|^{h;\f{\beta_0}\ep}_{\dot{B}_{p,1}^{s-1}}(t_0)+
  \|\mathbb{Q}\vue\|^{h;\f{\beta_0}\ep}_{\dot{B}_{p,1}^{s}}(t_0)\right)+C\|(F,\mathbb{Q}G)\|^{h;\f{\beta_0}\ep}_{L^1\Big(I;\dot{B}_{p,1}^{s}\Big)}\nonumber\\
  &\quad+\f C \ep\|J\|^{h;\f{\beta_0}\ep}_{L^1\Big(I;\dot{B}_{p,1}^{s-1}\Big)}+C\ep\sum\limits_{2^j\geq\f{\beta_0}\ep}2^{j(s+1)}\left(\|\dot{S}_{j-2}\Div\vue\cdot\dot{ \Delta}_j\vae\|_{L^1(I;L^p)}+\|\dot{\Delta}_j F+\dot{S}_{j-2}\vue\Grad \dot{\De}_j \vae\|_{L^1(I;L^p)}\right)\label{est12}.
\end{align}
\end{Lemma}

\bProof The proof is based on an analysis of the effective velocity  $\vc{\varpi}^\ep:=\mathbb{Q}\vue+\f1 {\ep\nu}\Grad (-\Delta)^{-1}\vae$ (see \cite{D2,F,H} for more details). The target system reads as
\begin{equation}\label{equ6}
  \left\{\bali
&\ep\p_t \vae+\f \vae {\ep\nu}=\ep F-\Div \varpi^\ep,\\
&\p_t \vc{\varpi}^\ep-\nu\Delta \vc{\varpi}^\ep+\f {\Grad\vthe} \ep=\f1 {\ep\nu}\Grad (-\Delta)^{-1}F+\mathbb{Q}G+\f {\vc{\varpi}^\ep} {\nu\ep^2}-\f1 {\nu^2\ep^3}\Grad(-\Delta)^{-1}\vae,\\
&\p_t\vthe-\kappa\Delta\vthe+\f \vae {\ep\nu}=J-\Div \varpi^\ep.
\eali\right.
\end{equation}
Applying $\dot{\Delta}_j$ to system \eqref{equ6}, we have
\begin{equation}\label{equ7}
\left\{\bali
&\ep\p_t \dot{\De}_j \vae+\f {\dot{\De}_j \vae} {\ep\nu}+\ep \dot{S}_{j-2}\vue\Grad \dot{\De}_j \vae=\ep \left(\dot{\De}_j F+\dot{S}_{j-2}\vue\Grad \dot{\De}_j \vae\right)-\Div \dot{\De}_j\vc{\varpi}^\ep,\\
&\p_t \dot{\De}_j\vc{\varpi}^\ep-\nu\De \dot{\De}_j\vc{\varpi}^\ep+\f {\Grad\dot{\De}_j\vthe} \ep=\f1 {\ep\nu}\Grad (-\De)^{-1}\dot{\De}_j F+\mathbb{Q}\dot{\De}_j G+\f{\dot{\De}_j\vc{\varpi}^\ep} {\ep^2\nu}-\f{\Grad(-\De)^{-1}\dot{\De}_j \vae} {\nu^2\ep^3},\\
&\p_t\dot{\De}_j\vthe-\kappa\Delta\dot{\De}_j\vthe+\f {\dot{\De}_j\vae} {\ep\nu}=\dot{\De}_jJ-\Div \dot{\De}_j\varpi^\ep.
\eali\right.
\end{equation}
Multiplying equations $\eqref{equ7}_{1,2,3}$ by $\dot{\De}_j\vae|\dot{\De}_j\vae|^{p-2}$, $\dot{\De}_j\vc{\varpi}^\ep|\dot{\De}_j\vc{\varpi}^\ep|^{p-2}$ and $\dot{\De}_j\vthe|\dot{\De}_j\vthe|^{p-2}$ respectively, integrating over $\mathbb{R}^d$, then dividing by $|\dot{\De}_j\vae|^{p-1},|\dot{\De}_j\vc{\varpi}^\ep|^{p-1},|\dot{\De}_j\vthe|^{p-1}$ respectively and integrating over $I$, we arrive at
\begin{align}
\ep\|\dot{\De}_j\vae\|_{L^\infty(I;L^p)}+\f1 {\ep}\|\dot{\De}_j\vae&\|_{L^1(I;L^p)}\leq C\ep\|\dot{\De}_j\vae(t_0)\|_{L^p}+C\ep\|\dot{S}_{j-2}\Div\vue\cdot\dot{ \Delta}_j\vae\|_{L^1(I;L^p)}\nonumber\\
&+C\ep\|\dot{\Delta}_j F+\dot{S}_{j-2}\vue\Grad \dot{\De}_j \vae\|_{L^1(I;L^p)}+C2^j\|\dot{\De}_j\vc{\varpi}^\ep\|_{L^1(I;L^p)},\label{est4}
\end{align}
\begin{align}
\|\dot{\De}_j\vc{\varpi}^\ep\|_{L^\infty(I;L^p)}+2^{2j}&\|\dot{\De}_j\vc{\varpi}^\ep\|_{L^1(I;L^p)}\leq C\|\dot{\De}_j\vc{\varpi}^\ep(t_0)\|_{L^p}
+C\f{2^j}{\ep}\|\dot{\De}_j\vthe\|_{L^1(I;L^p)}+\f C {\ep2^j}\|\dot{\De}_j F\|_{L^1(I;L^p)}\nonumber\\
&+C\|\dot{\De}_j \mathbb{Q}G\|_{L^1(I;L^p)}+\f C{(\ep2^j)^2}2^{2j}\|\dot{\De}_j\vc{\varpi}^\ep\|_{L^1(I;L^p)}+\f C{(\ep2^j)^2}\f{2^j}\ep\|\dot{\De}_j\vae\|_{L^1(I;L^p)}\label{est5}
\end{align}
and
\begin{align}
\|\dot{\De}_j\vthe\|_{L^\infty(I;L^p)}+2^{2j}\|\dot{\De}_j\vthe\|_{L^1(I;L^p)}\leq &C\|\dot{\De}_j\vthe(t_0)\|_{L^p}+\f C {\ep}\|\dot{\De}_j\vae\|_{L^1(I;L^p)}\nonumber\\
&+C2^j\|\dot{\De}_j\vc{\varpi}^\ep\|_{L^1(I;L^p)}+C\|\dot{\De}_j J\|_{L^1(I;L^p)}\label{est6}
\end{align}

Let $\delta$ be a small positive constant which will be determined later. Multiplying  both sides of \eqref{est4} by $2^j\delta$, \eqref{est5} by $C\delta$, and \eqref{est6} by $\f1 {2^j\ep}$, and then adding the resulting formulas yields
\begin{align}
&\ep\delta  2^j\|\dot{\De}_j\vae\|_{L^\infty(I;L^p)}+\Big(\delta-\f{C^2\delta}{(\ep2^j)^2}-\f{C\ep}{(\ep2^j)^2}\Big)\f{2^j} {\ep}\|\dot{\De}_j\vae\|_{L^1(I;L^p)}+C\delta\|\dot{\De}_j\vc{\varpi}^\ep\|_{L^\infty(I;L^p)}\nonumber\\
&+\Big(C\delta-\f{C^2\delta}{(\ep2^j)^2}-\f{C\ep}{(2^{j}\ep)^2}\Big)2^{2j}\|\dot{\De}_j\vc{\varpi}^\ep\|_{L^1(I;L^p)}
+\f1 {2^j\ep}\|\dot{\De}_j\vthe\|_{L^\infty(I;L^p)}+\Big(1-C^2\delta\Big)\f{2^{2j}} {2^j\ep}\|\dot{\De}_j\vthe\|_{L^1(I;L^p)}\nonumber\\
&\leq C^\ast \Big(\ep2^j\|\dot{\De}_j\vae(t_0)\|_{L^p}+\|\dot{\De}_j\vc{\varpi}^\ep(t_0)\|_{L^p}+\f1 {2^j\ep}\|\dot{\De}_j\vthe(t_0)\|_{L^p}\Big)+C^\ast\ep2^j\|\dot{S}_{j-2}\Div\vue\cdot\dot{ \Delta}_j\vae\|_{L^1(I;L^p)}\nonumber\\
&+C^\ast\ep2^j\|\dot{\Delta}_j F+\dot{S}_{j-2}\vue\Grad \dot{\De}_j \vae\|_{L^1(I;L^p)}+C^\ast\Big(\f 1 {\ep2^j}\|\dot{\De}_j F\|_{L^1(I;L^p)}+\|\dot{\De}_j \mathbb{Q}G\|_{L^1(I;L^p)}+\f 1 {\ep2^j}\|\dot{\De}_j J\|_{L^1(I;L^p)}\Big)\label{est7}.
\end{align}
Now let $\delta=\f1{4C^2}$ and $\ep\leq1$. There exist a positive constant $\beta_0=\beta_0(d,\nu)$ and for any
$j$ satisfying $2^j\ep\geq\beta_0$, there holds that
\begin{align}\label{est8}
\delta-\f{C^2\delta}{(\ep2^j)^2}-\f{C}{(\ep2^j)^2}\geq\f{\delta} 2 \ \ \text{ and }\ \ C\delta-\f{C^2\delta}{(\ep2^j)^2}-\f{C}{(2^{j}\ep)^2}\geq \f {C\delta} 2.
\end{align}
For such $j$, we also have that for any $t\in I$,
\begin{align}
\ep2^j\|\dot{\De}_j\vae(t)\|_{L^p}+\|\dot{\De}_j\vc{\varpi}^\ep(t)\|_{L^p}\approx \ep2^j\|\dot{\De}_j\vae(t)\|_{L^p}+\|\dot{\De}_j\mathbb{Q}\vue(t)\|_{L^p} \label{est9} \\
\f{2^j}{\ep}\|\dot{\De}_j\vae(t)\|_{L^p}+2^{2j}\|\dot{\De}_j\vc{\varpi}^\ep(t)\|_{L^p}\approx \f{2^j}\ep\|\dot{\De}_j\vae(t)\|_{L^p}+2^{2j}\|\dot{\De}_j\mathbb{Q}\vue(t)\|_{L^p}.\label{est10}
\end{align}
Inserting \eqref{est8} into \eqref{est7} and utilizing \eqref{est9}, \eqref{est10}, we get
\begin{align}
\ep2^j&\|\dot{\De}_j\vae\|_{L^\infty(I;L^p)}+\f{2^j} {\ep}\|\dot{\De}_j\vae\|_{L^1(I;L^p)}+\|\dot{\De}_j\vc{\varpi}^\ep\|_{L^\infty(I;L^p)}+2^{2j}\|\dot{\De}_j\mathbb{Q}\vue\|_{L^1(I;L^p)}\nonumber\\
&\quad+\f1 {2^j\ep}\|\dot{\De}_j\vthe\|_{L^\infty(I;L^p)}+\f{2^{2j}} {2^j\ep}\|\dot{\De}_j\vthe\|_{L^1(I;L^p)}\nonumber\\
&\leq C\Big(\ep2^j\|\dot{\De}_j\vae(t_0)\|_{L^p}+\|\dot{\De}_j\mathbb{Q}\vue(t_0)\|_{L^p}+\f1 {2^j\ep}\|\dot{\De}_j\vthe(t_0)\|_{L^p}\Big)+C\ep2^j\|\dot{S}_{j-2}\Div\vue\cdot\dot{ \Delta}_j\vae\|_{L^1(I;L^p)}\nonumber\\
&\quad+C\ep2^j\|\dot{\Delta}_j F+\dot{S}_{j-2}\vue\Grad \dot{\De}_j \vae\|_{L^1(I;L^p)}+C\Big(\|(\dot{\De}_j F,\dot{\De}_j \mathbb{Q}X)\|_{L^1(I;L^p)}+\f 1 {\ep2^j}\|\dot{\De}_j J\|_{L^1(I;L^p)}\Big).\label{est11}
\end{align}
Multiplying both sides of \eqref{est11} by $2^{js}$ and summing over $j$ which satisfies $2^j\ep\geq\beta_0$, we arrive at \eqref{est12}.\ \ $\Box$

\begin{Lemma}\label{Le3.2}Let $s\in\mathbb{R},0\leq\alpha\leq\f\beta \ep<\infty$ and $I\subset\mathbb{R}$ be an temporal interval with the initial time $t_0$. There exist a constant $C=C(d,\nu,\kappa)$ such that
\begin{align}
&\|(\vae,\mathbb{Q}\vue,\vthe)\|^{m;\alpha,\f{\beta}\ep}_{\widetilde{L}^\infty(I;\dot{B}_{2,1}^{s-1})}+\f1 {(1+\beta)^2}\|(\vae,\mathbb{Q}\vue,\vthe)\|^{m;\alpha,\f{\beta}\ep}_{L^1(I;\dot{B}_{2,1}^{s+1})}\nonumber\\
&\quad\leq C\left(\|(\vae,\mathbb{Q}\vue,\vthe)(t_0)\|^{m;\alpha,\f{\beta}\ep}_{\dot{B}_{2,1}^{s-1}}
+\|(F,\mathbb{Q}G,J)\|^{m;\alpha,\f{\beta}\ep}_{L^1(I;\dot{B}_{2,1}^{s-1})}\nonumber\right).
\end{align}
\end{Lemma}
The proof of Lemma \ref{Le3.2} is similar to the argument in section 4 of \cite{D6}, so we omit the proof.

Adding equations \eqref{equ8}$_1$ and \eqref{equ8}$_3$ , dividing the resulting equation by $\sqrt{2}$, and applying the operator $\mathbb{Q}$ to equation \eqref{equ8}$_2$ , we get
\begin{equation}
  \left\{\bali
  &\p_t\tau^\ep+\f{\sqrt{2}}{\ep}\Div\mathbb{Q}\vue=\f1 {\sqrt{2}}(F+J),\\
  &\p_t\mathbb{Q}\vue-\nu\De\mathbb{Q}\vue+\f{\sqrt{2}}\ep\Grad\tau^\ep=\mathbb{Q}G.\\
  \eali
  \right.\nonumber
\end{equation}
\begin{Lemma}\label{Be9}Let $2\leq p,q\leq\infty,0<\alpha<\infty,\varsigma \in\mathbb{R}$ and $I\subset\mathbb{R}$ be an temporal interval with the initial time $t_0$. If
\[
\f1 q\leq \f{d-1} 2\left(\f1 2-\f1 p \right),\ \ (p,q,d)\neq(\infty,2,3),
\]
then there exist a positive constant $C=(d,p,q,\nu,\kappa)$ such that
\begin{align}
&\|(\tau^\ep,\mathbb{Q}\vue)\|^{l;\alpha}_{\widetilde{L}^q\left(I;\dot{B}_{p,1}^{\f d p+\f1 q+\varsigma}\right)}\leq C\ep^{\f1 q}(1+\ep\alpha)^2\left(\|(\vae,\mathbb{Q}\vue,\vthe)(t_0)\|^{l;\alpha}_{\dot{B}_{2,1}^{\f d 2+\varsigma}}+\|(F,\mathbb{Q}G,J)\|^{l;\alpha}_{L^1(I;\dot{B}_{2,1}^{\f d 2+\varsigma})}\right),\label{est16}\\
&\|(\tau^\ep,\mathbb{Q}\vue\|^{l;\alpha}_{\widetilde{L}^{2}\left(I;\dot{B}_{p,1}^{\f d p-(\f1 2-\f1 p)}\right)}\leq C\ep^{\f1 2-\f1 p}\left(\|(\vae,\mathbb{Q}\vue,\vthe)(t_0)\|^{l;\alpha}_{\dot{B}_{2,1}^{\f d 2-1}}+\|(F,\mathbb{Q}G,J)\|^{l;\alpha}_{L^1(I;\dot{B}_{2,1}^{\f d 2-1})}\right)\label{est17}.
\end{align}
\end{Lemma}
\bProof From Proposition 7.1 of \cite{D3}, we have that
\begin{align}\label{est13}
 \|(\tau^\ep,\mathbb{Q}\vue)\|^{l;\alpha}_{\widetilde{L}^q\left(I;\dot{B}_{p,1}^{\f d p+\f1 q+\varsigma}\right)}&\leq C\ep^{\f1 q}\left(\|(\tau^\ep,\mathbb{Q}\vue)(t_0)\|^{l;\alpha}_{\dot{B}_{2,1}^{\f d 2+\varsigma}}+\|(F+J,\mathbb{Q}G)\|^{l;\alpha}_{L^1(I;\dot{B}_{2,1}^{\f d 2+\varsigma})}\right.\nonumber\\
 &\quad+\left.\|\De\mathbb{Q}\vue\|^{l;\alpha}_{L^1(I;\dot{B}_{2,1}^{\f d 2+\varsigma})}.\right)
\end{align}
Utilizing Lemma \eqref{Le3.2} yields
\begin{align}
\|\De\mathbb{Q}\vue\|^{l;\alpha}_{L^1(I;\dot{B}_{2,1}^{\f d 2+\varsigma})}&\leq C\|\mathbb{Q}\vue\|^{l;\alpha}_{L^1(I;\dot{B}_{2,1}^{\f d 2+\varsigma+2})}\nonumber\\
&\leq C(1+\ep\alpha)^2\left(\|(\vae,\mathbb{Q}\vue,\vthe)(t_0)\|^{l;\alpha}_{\dot{B}_{2,1}^{\f d 2+\varsigma}}
+\|(F,\mathbb{Q}G,J)\|^{l;\alpha}_{L^1(I;\dot{B}_{2,1}^{\f d 2+\varsigma})}\label{est15}\right)
\end{align}
Combining \eqref{est13} and \eqref{est15} yields \eqref{est16}. The proof of \eqref{est17} is similar to the argument in section 5 of \cite{D2}, so we omit the proof.   \ \   $\Box$

\section{The energy functional $E^\ep$ and quantity $M^{\ep,\alpha}_{p,q}$}\label{s4}
Considering linear estimates of the previous section, we introduce the energy functional $E^\ep[a,\vu,\vartheta](I)$ and quantity $M^{\ep,\alpha}_{p,q}[a,\vu,\vartheta;\tau,\sigma](I)$. Let $\beta_0$ be the positive constant appearing in Lemma \ref{Le3.1} and $I\subset \mathbb{R}$ be a interval.
\begin{Definition} For $2\leq p<\infty,0<\alpha<\infty$ and $\alpha\leq\f{\beta_0}\ep$,
\[
\|(a_0,\vu_0,\vartheta_0)\|_{X^\ep_p}:=\|(a_0,\mathbb{Q}\vu_0,\vartheta_0)\|^{h;\alpha}_{X^\ep_p}+\|(a_0,\mathbb{Q}\vu_0,\vartheta_0)\|^{l;\alpha}_{\dot{B}^{\f d 2-1}_{2,1}}+\|\mathbb{P}\vu_0\|_{\dot{B}^{\f d 2-1}_{2,1}}
\]
where
\[
\|(a_0,\vu_0,\vartheta_0)\|^{h;\alpha}_{X^\ep_p}:=\ep\|a_0\|^{h;\f{\beta_0}\ep}_{\dot{B}^{\f d p}_{p,1}}+\f1\ep\|\vartheta_0\|^{h;\f{\beta_0}\ep}_{\dot{B}^{\f d p-2}_{p,1}}+\|\mathbb{Q}\vu_0\|^{h;\f{\beta_0}\ep}_{\dot{B}^{\f d p-1}_{p,1}}+\|(a_0,\mathbb{Q}\vu_0,\vartheta_0)\|^{m;\alpha,\f{\beta_0}\ep}_{\dot{B}^{\f d 2-1}_{2,1}}
\]
for $(a_0,\vu_0,\vartheta_0)\in{\dot{B}_{2,1}^{\f d 2}(\mathbb{R}^d)}\cap{\dot{B}_{2,1}^{\f d 2-1}(\mathbb{R}^d)}\times\Big(\dot{B}_{2,1}^{\f d 2-1}(\mathbb{R}^d)\Big)^d\times\dot{B}_{2,1}^{\f d 2-2}(\mathbb{R}^d)$.
\end{Definition}
It is clear that
\[
\|(a_0,\vu_0,\vartheta_0)\|_{X^\ep_{2}}\approx \ep\|a_0\|^{h;\f{\beta_0}\ep}_{\dot{B}^{\f d 2}_{2,1}}+\f1\ep\|\vartheta_0\|^{h;\f{\beta_0}\ep}_{\dot{B}^{\f d 2-2}_{2,1}}+\|(a_0,\vartheta_0)\|^{l;\f{\beta_0}\ep}_{\dot{B}^{\f d 2-1}_{2,1}}+\|\vu_0\|_{\dot{B}^{\f d 2-1}_{2,1}}.
\]

\begin{Definition} Let $2\leq p,q<\infty,0<\alpha\leq\f{\beta_0}\ep$.
\begin{itemize}
  \item Let $E^\ep(I)$ be the set of all space-time distributions $(a,\vu,\vartheta)$ on $I\times\mathbb{R}^d$ with the following norm finite:
 \begin{align*}
  E^\ep[a,\vu,\vartheta](I):=&\ep\|a\|^{h;\f{\beta_0}\ep}_{L^\infty(I;\dot{B}^{\f d 2}_{2,1})}+\f1 \ep \|a\|^{h;\f{\beta_0}\ep}_{L^1(I;\dot{B}^{\f d 2}_{2,1})}+\f1 \ep \|\vartheta\|^{h;\f{\beta_0}\ep}_{L^\infty(I;\dot{B}^{\f d 2-2}_{2,1})\cap L^1(I;\dot{B}^{\f d 2}_{2,1})}\\
  &\quad+\|(a,\vartheta)\|^{l;\f{\beta_0}\ep}_{L^\infty(I;\dot{B}^{\f d 2-1}_{2,1})\cap L^1(I;\dot{B}^{\f d 2+1}_{2,1})}+\|\vu\|_{L^\infty(I;\dot{B}^{\f d 2-1}_{2,1})\cap L^1(I;\dot{B}^{\f d 2+1}_{2,1})}
 \end{align*}
  \item Let $N^{\ep,\alpha}_{p,q}(I)$ be the set of all space-time distributions $(a,\vu,\vartheta,\tau,\sigma)$ on $I\times\mathbb{R}^d$ with the following norm finite:
  \begin{align*}
  \|(a,\vu,\vartheta;\tau,\sigma)\|_{N^{\ep,\alpha}_{p,q}(I)}&:= \|(a,\mathbb{Q}\vu,\vartheta)\|^{h;\alpha}_{N^{\ep,\alpha}_{p,q}(I)}+\|(\tau,\mathbb{Q}\vu,\sigma)\|^{l;\alpha}_{N^{\ep,\alpha}_{p,q}(I)}\\
  &\quad+\|\mathbb{P}\vu\|_{L^q(I;\dot{B}^{\f d p-1+\f2 q}_{p,1})\cap L^1(I;\dot{B}^{\f d p+1}_{p,1})}
 \end{align*}
  where
 \begin{align*}
   \|(a,\vu,\vartheta)\|^{h;\alpha}_{N^{\ep,\alpha}_{p,q}(I)}&:=\ep\|a\|^{h;\f{\beta_0}\ep}_{L^\infty(I;\dot{B}^{\f d p}_{p,1})}+\f1 \ep \|a\|^{h;\f{\beta_0}\ep}_{L^1(I;\dot{B}^{\f d p}_{p,1})}+\f1 \ep \|\vartheta\|^{h;\f{\beta_0}\ep}_{L^\infty(I;\dot{B}^{\f d p-2}_{p,1})\cap L^1(I;\dot{B}^{\f d p}_{p,1})}\\
   &\quad+\|\mathbb{Q}\vu\|^{h;\f{\beta_0}\ep}_{L^\infty(I;\dot{B}^{\f d p-1}_{p,1})\cap L^1(I;\dot{B}^{\f d p+1}_{p,1})}+\|(a,\mathbb{Q}\vu,\vartheta)\|^{m;\alpha;\f{\beta_0}\ep}_{L^\infty(I;\dot{B}^{\f d 2-1}_{2,1})\cap L^1(I;\dot{B}^{\f d 2+1}_{2,1})},\\
   \|(\tau,\vu,\sigma)\|^{l;\alpha}_{N^{\ep,\alpha}_{p,q}(I)}&:=\|(\tau,\mathbb{Q}\vu)\|^{l;\alpha}_{L^q(I;\dot{B}^{\f d p-1+\f2 q}_{p,1})}+\|\sigma\|^{l;\alpha}_{L^q(I;\dot{B}^{\f d p-1+\f2 q}_{p,1})\cap L^{q'}(I;\dot{B}^{\f d p-1+\f2 {q'}}_{p,1})}
  \end{align*}
  \item We set
  \begin{align*}
  M^{\ep,\alpha}_{p,q}[a,\vu,\vartheta;\tau,\sigma](I)&:=\alpha\ep E^\ep[a,\mathbb{Q}\vu,\vartheta](I)+ \|(a,\vu,\vartheta;\tau,\sigma)\|_{N^{\ep,\alpha}_{p,q}(I)}\\
  &\quad+ \left(\|(\tau,\mathbb{Q}\vu)\|^{l;\alpha}_{{L^q(I;\dot{B}^{\f d p-1+\f2 q}_{p,1})}}\right)^\f1 {q-1}\Big(  E^\ep[a,\mathbb{Q}\vu,\vartheta](I)\Big)^\f{q-2}{q-1}.
  \end{align*}
\end{itemize}
\end{Definition}
\begin{Proposition}\label{pro1}  Let $2\leq p,q<\infty,0<\alpha\leq\f{\beta_0}\ep$ and $(a,\vu,\vartheta,\tau,\sigma)$ be space-time distribution on $I\times\mathbb{R}^d$. There exist a constant $C_1=C_1(d,p,q)$ such that
\begin{align*}
&\ep\|a\|_{L^\infty(I;\dot{B}^{\f d 2}_{2,1})}+\|a\|_{L^q(I;\dot{B}^{\f d 2-1+\f2 q}_{2,1})}+\|a\|^{l;\f{4\beta_0}\ep}_{L^{q'}(I;\dot{B}^{\f d 2-1+\f2 {q'}}_{2,1})}+\|a\|_{L^2(I;\dot{B}^{\f d 2}_{2,1})}\leq C_1 E^\ep[a,\vu,\vartheta](I),\\
&\ep\|a\|_{L^\infty(I;\dot{B}^{\f d p}_{p,1})}+\|a\|^{h;\alpha}_{L^2(I;\dot{B}^{\f d p}_{p,1})}+\|a\|^{h;\alpha}_{L^q(I;\dot{B}^{\f d p-1+\f2 q}_{p,1})}+\|a\|^{m;\alpha,\f{4\beta_0}\ep}_{L^{q'}(I;\dot{B}^{\f d p-1+\f2 {q'}}_{p,1})}+\|\vartheta\|^{m;\alpha,\f{4\beta_0}\ep}_{L^\infty(I;\dot{B}^{\f d p-1}_{p,1})\cap L^{1}(I;\dot{B}^{\f d p+1}_{p,1})}\\
&+\|\vu\|_{L^q(I;\dot{B}^{\f d p-1+\f2 q}_{p,1})}+\|\vu\|_{L^{q'}(I;\dot{B}^{\f d p-1+\f2 {q'}}_{p,1})}+\|\vu\|_{L^2(I;\dot{B}^{\f d p}_{p,1})}\leq C_1  M^{\ep,\alpha}_{p,q}[a,\vu,\vartheta;\tau,\sigma](I),\\
&\|a\|_{L^q(I;\dot{B}^{\f d p-1+\f2 q}_{p,1})}+\|a\|_{L^2(I;\dot{B}^{\f d p}_{p,1})}+\|a\|^{l;\f{4\beta_0}\ep}_{L^{q'}(I;\dot{B}^{\f d p-1+\f2 {q'}}_{p,1})}+\|\vartheta\|^{l;\f{\beta_0}\ep}_{L^q(I;\dot{B}^{\f d p-1+\f2 q}_{p,1})\cap L^{q'}(I;\dot{B}^{\f d p-1+\f2 {q'}}_{p,1})}\\
&\quad\leq C_1 M^{\ep,\alpha}_{p,q}[a,\vu,\vartheta;\f{a+\vartheta}{\sqrt{2}},\f{\vartheta-a}{\sqrt{2}}](I).
\end{align*}
\end{Proposition}
\bProof The proof of Proposition \ref{pro1} is simply a matter of using the interpolation inequality and Bernstein's inequality for different frequencies, of which we will compute only a few, the rest being similarly obtainable.
\begin{align*}
\ep\|a\|_{L^\infty(I;\dot{B}^{\f d p}_{p,1})}&\leq \ep\|a\|^{h;\f{\beta_0}\ep}_{L^\infty(I;\dot{B}^{\f d p}_{p,1})}+\ep\|a\|^{m;\alpha,\f{\beta_0}\ep}_{L^\infty(I;\dot{B}^{\f d p}_{p,1})}+\ep\|a\|^{l;\alpha}_{L^\infty(I;\dot{B}^{\f d p}_{p,1})}\\
&\leq \ep\|a\|^{h;\f{\beta_0}\ep}_{L^\infty(I;\dot{B}^{\f d p}_{p,1})}+\ep\cdot\f{\beta_0}\ep\|a\|^{m;\alpha,\f{\beta_0}\ep}_{L^\infty(I;\dot{B}^{\f d p-1}_{p,1})}+\ep\alpha\|a\|^{l;\alpha}_{L^\infty(I;\dot{B}^{\f d p-1}_{p,1})}\\
&\leq \ep\|a\|^{h;\f{\beta_0}\ep}_{L^\infty(I;\dot{B}^{\f d p}_{p,1})}+C_1\|a\|^{m;\alpha,\f{\beta_0}\ep}_{L^\infty(I;\dot{B}^{\f d 2-1}_{2,1})}+C_1\ep\alpha\|a\|^{l;\alpha}_{L^\infty(I;\dot{B}^{\f d 2-1}_{2,1})}\\
&\leq C_1M^{\ep,\alpha}_{p,q}[a,\vu,\vartheta;\tau,\sigma](I)
\end{align*}
\begin{align*}
&\|a\|_{L^q(I;\dot{B}^{\f d p-1+\f2 q}_{p,1})}\\
&\leq \|a\|^{h;\f {\beta_0}\ep}_{L^q(I;\dot{B}^{\f d p-1+\f2 q}_{p,1})}+\|a\|^{m;\alpha,\f {\beta_0}\ep}_{L^q(I;\dot{B}^{\f d p-1+\f2 q}_{p,1})}+\|a\|^{l;\alpha}_{L^q(I;\dot{B}^{\f d p-1+\f2 q}_{p,1})}\\
&\leq C_1\left(\f{\beta_0}\ep\right)^{-1+\f2 q} \|a\|^{h;\f {\beta_0}\ep}_{L^q(I;\dot{B}^{\f d p}_{p,1})}+C_1\|a\|^{m;\alpha,\f {\beta_0}\ep}_{L^q(I;\dot{B}^{\f d 2-1+\f2 q}_{2,1})}+\|a\|^{l;\alpha}_{L^q(I;\dot{B}^{\f d p-1+\f2 q}_{p,1})}\\
&\leq C_1\left(\ep\|a\|^{h;\f {\beta_0}\ep}_{L^\infty(I;\dot{B}^{\f d p}_{p,1})}\right)^{1-\f1 q}\left(\f1 \ep\|a\|^{h;\f {\beta_0}\ep}_{L^1(I;\dot{B}^{\f d p}_{p,1})}\right)^{\f1 q}+C_1\left(\|a\|^{m;\alpha,\f {\beta_0}\ep}_{L^\infty(I;\dot{B}^{\f d 2-1}_{2,1})}\right)^{1-\f1 q}\left(\|a\|^{m;\alpha,\f {\beta_0}\ep}_{L^1(I;\dot{B}^{\f d 2+1}_{2,1})}\right)^{\f1 q}\\
&\quad+\|a\|^{l;\alpha}_{L^q(I;\dot{B}^{\f d p-1+\f2 q}_{p,1})}\\
&\leq C_1\left(\ep\|a\|^{h;\f {\beta_0}\ep}_{L^\infty(I;\dot{B}^{\f d p}_{p,1})}+\f1 \ep\|a\|^{h;\f {\beta_0}\ep}_{L^1(I;\dot{B}^{\f d p}_{p,1})}\right)+C_1\|a\|^{m;\alpha;\f{\beta_0}\ep}_{L^\infty(I;\dot{B}^{\f d 2-1}_{2,1})\cap L^1(I;\dot{B}^{\f d 2+1}_{2,1})}+C_1\|\f{a+\vartheta}{\sqrt{2}}-\f{\vartheta-a}{\sqrt{2}}\|^{l;\alpha}_{L^q(I;\dot{B}^{\f d p-1+\f2 q}_{p,1})}\\
&\leq  C_1 M^{\ep,\alpha}_{p,q}[a,\vu,\vartheta;\f{a+\vartheta}{\sqrt{2}},\f{\vartheta-a}{\sqrt{2}}](I),
\end{align*}

\begin{align*}
 &\|a\|^{l;\f{4\beta_0}\ep}_{L^{q'}(I;\dot{B}^{\f d p-1+\f2 {q'}}_{p,1})}\leq \|a\|^{m;\f{\beta_0}\ep,\f{4\beta_0}\ep}_{L^{q'}(I;\dot{B}^{\f d p-1+\f2 {q'}}_{p,1})}+\|a\|^{m;\alpha,\f{\beta_0}\ep}_{L^{q'}(I;\dot{B}^{\f d p-1+\f2 {q'}}_{p,1})}+\|a\|^{l;\alpha}_{L^{q'}(I;\dot{B}^{\f d p-1+\f2 {q'}}_{p,1})}\\
&\leq \left(\f{4\beta_0}\ep\right)^{-1+\f2 {q'}}\|a\|^{m;\f{\beta_0}\ep,\f{4\beta_0}\ep}_{L^{q'}(I;\dot{B}^{\f d p}_{p,1})}+C_1\|a\|^{m;\alpha;\f{\beta_0}\ep}_{L^\infty(I;\dot{B}^{\f d 2-1}_{2,1})\cap L^1(I;\dot{B}^{\f d 2+1}_{2,1})}+\|\f{a+\vartheta}{\sqrt{2}}-\f{\vartheta-a}{\sqrt{2}}\|^{l;\alpha}_{L^{q'}(I;\dot{B}^{\f d p-1+\f2 {q'}}_{p,1})}\\
&\leq C_1\left(\ep\|a\|^{h;\f {\beta_0}\ep}_{L^\infty(I;\dot{B}^{\f d p}_{p,1})}+\f1 \ep\|a\|^{h;\f {\beta_0}\ep}_{L^1(I;\dot{B}^{\f d p}_{p,1})}\right)+C_1\|a\|^{m;\alpha;\f{\beta_0}\ep}_{L^\infty(I;\dot{B}^{\f d 2-1}_{2,1})\cap L^1(I;\dot{B}^{\f d 2+1}_{2,1})}\\
&\quad+C_1\|\f{\vartheta-a}{\sqrt{2}}\|^{l;\alpha}_{L^q(I;\dot{B}^{\f d p-1+\f2 q}_{p,1})\cap L^{q'}(I;\dot{B}^{\f d p-1+\f2 {q'}}_{p,1})}
+\left({\|\f{a+\vartheta}{\sqrt{2}}\|^{l;\alpha}_{L^q(I;\dot{B}^{\f d p-1+\f2 q}_{p,1})}}\right)^{\f1 {q-1}}\left({\|\f{a+\vartheta}{\sqrt{2}}\|^{l;\alpha}_{L^1(I;\dot{B}^{\f d 2+1}_{2,1})}}\right)^{\f{q-2} {q-1}}\\
&\leq  C_1 M^{\ep,\alpha}_{p,q}[a,\vu,\vartheta;\f{a+\vartheta}{\sqrt{2}},\f{\vartheta-a}{\sqrt{2}}](I)\\
\end{align*}

\begin{align*}
&\|a\|_{L^2(I;\dot{B}^{\f d p}_{p,1})}\\
&\leq \|a\|^{h;\f {\beta_0}\ep}_{L^2(I;\dot{B}^{\f d p}_{p,1})}+\|a\|^{m;\alpha,\f {\beta_0}\ep}_{L^2(I;\dot{B}^{\f d p}_{p,1})}+\|a\|^{l;\alpha}_{L^2(I;\dot{B}^{\f d p}_{p,1})}\\
&\leq \left(\ep\|a\|^{h;\f {\beta_0}\ep}_{L^\infty(I;\dot{B}^{\f d p}_{p,1})}\right)^{\f1 2}\left(\f1 \ep\|a\|^{h;\f {\beta_0}\ep}_{L^1(I;\dot{B}^{\f d p}_{p,1})}\right)^{\f1 2}\\
&\quad+\left(\|a\|^{m;\alpha,\f {\beta_0}\ep}_{L^\infty(I;\dot{B}^{\f d 2-1}_{2,1})}\right)^{\f1 2}\left(\|a\|^{m;\alpha,\f {\beta_0}\ep}_{L^1(I;\dot{B}^{\f d 2+1}_{2,1})}\right)^{\f1 2}+\Big(\|a\|^{l;\alpha}_{L^q(I;\dot{B}^{\f d p-1+\f2 q}_{p,1})}\Big)^{\f1 2}\Big(\|a\|^{l;\alpha}_{L^{q'}(I;\dot{B}^{\f d p-1+\f2 {q'}}_{p,1})}\Big)^{\f1 2}\\
&\leq  C_1 M^{\ep,\alpha}_{p,q}[a,\vu,\vartheta;\f{a+\vartheta}{\sqrt{2}},\f{\vartheta-a}{\sqrt{2}}](I).\ \ \Box
\end{align*}

\section{A prior estimates}\label{s5}
In this section, we persent a priori estimates. Let $\beta_0$ be the positive constants appearing in Lemma \ref{Le3.1} and $I\subset\mathbb{R}$ be an temporal interval with the initial time $t_0$.  We assume that $(\vae,\vue,\vthe)$ is a regular solution to systems \eqref{equ2} on $I$ with  $\ep\|\vae\|_{L^\infty(I;L^\infty)}\leq\f1 2$ and $(\Theta,\vv)$ is a regular solution to systems \eqref{equ5} supplement with initial data $\Theta_0=\f{\vartheta_0-a_0}{\sqrt{2}}$ and $\vv_0=\mathbb{P}\vu_0$. It is clear that $(\vae,\vue,\vthe)$ satisfies systems \eqref{equ8} with
\[
F=-\Div(\vae\vue),\ \  G=I(\ep\vae)\left(\f{\Grad(\vae+\vthe)}\ep+\Grad(\vae\vthe)-\mathcal{L}\vue\right)-\Grad(\vae\vthe)-\vue\cdot\Grad\vue,
\]
\[
J=\ep(1-I(\ep\vae))\big(2\mu|D\vue|^2+\lambda|\Div\vue|^2\big)-I(\ep\vae)\kappa\Delta\vthe-\Div(\vthe\vue).
\]

\begin{Lemma}\label{Le5.1}  Let $0<\alpha<\infty, \ep\leq\min\{\f{\beta_0}{\alpha},1\},2\leq p<d,2\leq q<\infty$. Then there exist a constant $C_2=C_2(d,p,\mu,\kappa,\nu)$ such that
\begin{align*}\label{est18}
E^\ep[\vae,\vue,\vthe](I)&\leq C_2\|(\vae,\vue,\vthe)(t_0)\|_{X^\ep_2}\\
&\quad+C_2 E^\ep[\vae,\vue,\vthe](I)\cdot \left(M^{\ep,\alpha}_{p,q}[\vae,\vue,\vthe](I)+\bigg(M^{\ep,\alpha}_{p,q}[\vae,\vue,\vartheta^\ep;\tau^\ep,\sigma^\ep](I)\bigg)^2\right).
\end{align*}
\end{Lemma}
\bProof
We find that
\[
\p_t\mathbb{P}\vue-\mu\De\mathbb{P}\vue=\mathbb{P}G.
\]
Then from Lemma \ref{Be3}, Lemma \ref{Le3.1} and Lemma \ref{Le3.2}, we get that
\begin{align}\label{est19}
 &E^\ep[\vae,\vue,\vthe](I)\leq C_2\|(\vae,\vue,\vthe)(t_0)\|_{X^\ep_2}+C_2\|(F,G)\|_{L^1(I;\dot{B}_{2,1}^{\f d 2 -1})}+\f{C_2}{\ep}\|J\|^{h;\f{\beta_0}\ep}_{L^1(I;\dot{B}_{2,1}^{\f d 2 -2})}\\
 &\quad+C_2\|J\|^{l;\f{\beta_0}\ep}_{L^1(I;\dot{B}_{2,1}^{\f d 2 -1})}+C_2\ep\sum\limits_{2^j\geq\f{\beta_0}\ep}2^{j\f d 2}\left(\|\dot{S}_{j-2}\Div\vue\cdot\dot{ \Delta}_j\vae\|_{L^1(I;L^2)}+\|\dot{\Delta}_j F+\dot{S}_{j-2}\vue\Grad \dot{\De}_j \vae\|_{L^1(I;L^2)}\right).\nonumber
\end{align}
We use Proposition \ref{pro1} frequently in our later calculations. We have that
\begin{align}\label{est20}
&\|\big(\Div(\vae\vue),I(\ep\vae)\mathcal{L}\vue,\vue\cdot\Grad\vue\big)\|_{L^1(I;\dot{B}_{2,1}^{\f d 2 -1})}\leq C_2\|\vae\|_{ L^2(I;\dot{B}_{p,1}^{\f d p})}\|\vue\|_{L^2(I;\dot{B}_{2,1}^{\f d 2})}\nonumber\\
&\quad+C_2\|\vae\|_{ L^2(I;\dot{B}_{2,1}^{\f d 2})}\|\vue\|_{L^2(I;\dot{B}_{p,1}^{\f d p})}+C_2\ep\|\vae\|_{L^\infty(I;\dot{B}_{2,1}^{\f d 2})}\|\vue\|_{L^1(I;\dot{B}_{p,1}^{\f d p+1})}\nonumber\\
&\quad+C_2\ep\|\vae\|_{L^\infty(I;\dot{B}_{p,1}^{\f d p})}\|\vue\|_{L^1(I;\dot{B}_{2,1}^{\f d 2+1})}+\|\vue\|_{L^2(I;\dot{B}_{2,1}^{\f d 2})}\|\vue\|_{L^2(I;\dot{B}_{p,1}^{\f d p})}\nonumber\\
&\leq C_2 E^\ep[\vae,\vue,\vthe](I)M^{\ep,\alpha}_{p,q}[\vae,\vue,\vartheta^\ep;\tau^\ep,\sigma^\ep](I),
\end{align}
\begin{align}\label{est21}
&\|I(\ep\vae)\f{\Grad(\vae+\vthe)}\ep\|_{L^1(I;\dot{B}_{2,1}^{\f d 2 -1})}\leq C_2\|\vae\|_{L^2(I;\dot{B}_{2,1}^{\f d 2})}\|\vae\|_{L^2(I;\dot{B}_{p,1}^{\f d p})}+C_2\ep\|\vae\|_{L^\infty(I;\dot{B}_{p,1}^{\f d p})}\f1 {\ep}\|\vthe\|^{h;\f{\beta_0}\ep}_{L^1(I;\dot{B}_{2,1}^{\f d 2})}\nonumber\\
&\quad+C_2\|\vae\|_{ L^2(I;\dot{B}_{p,1}^{\f d p})}\|\vthe\|^{l;\f{\beta_0}\ep}_{L^2(I;\dot{B}_{2,1}^{\f d 2})}+C_2\ep\|\vae\|_{L^\infty(I;\dot{B}_{2,1}^{\f d 2})}\f1 {\ep}\|\vthe\|^{h;\f{\beta_0}\ep}_{L^1(I;\dot{B}_{p,1}^{\f d p})}+\|\vae\|_{ L^2(I;\dot{B}_{2,1}^{\f d 2})}\|\vthe\|^{l;\f{\beta_0}\ep}_{L^2(I;\dot{B}_{p,1}^{\f d p})}\nonumber\\
&\leq C_2 E^\ep[\vae,\vue,\vthe](I)M^{\ep,\alpha}_{p,q}[\vae,\vue,\vartheta^\ep;\tau^\ep,\sigma^\ep](I).
\end{align}
A similar calculation to the one above yields
\begin{align}\label{est22}
\|\big(I(\ep\vae)\Grad(\vae\vthe),\Grad(\vae\vthe)\big)\|_{L^1(I;\dot{B}_{2,1}^{\f d 2 -1})}&\leq C_2 E^\ep[\vae,\vue,\vthe](I)\left(M^{\ep,\alpha}_{p,q}[\vae,\vue,\vartheta^\ep;\tau^\ep,\sigma^\ep](I)\right.\nonumber\\
&\quad\left.+\big(M^{\ep,\alpha}_{p,q}[\vae,\vue,\vartheta^\ep;\tau^\ep,\sigma^\ep](I)\big)^2\right).
\end{align}
Calculation gives that
\begin{align}\label{est23}
&\|(1-I(\ep\vae))\big(2\mu|D\vue|^2+\lambda|\Div\vue|^2\big)\|^{h;\f{\beta_0}\ep}_{L^1(I;\dot{B}_{2,1}^{\f d 2 -2})}+\ep\|(1-I(\ep\vae))\big(2\mu|D\vue|^2+\lambda|\Div\vue|^2\big)\|^{l;\f{\beta_0}\ep}_{L^1(I;\dot{B}_{2,1}^{\f d 2 -1})}\nonumber\\
&\leq \|(1-I(\ep\vae))\big(2\mu|D\vue|^2+\lambda|\Div\vue|^2\big)\|_{L^1(I;\dot{B}_{2,1}^{\f d 2 -2})}\nonumber\\
&\leq C_2 E^\ep[\vae,\vue,\vthe](I)\Big(M^{\ep,\alpha}_{p,q}[\vae,\vue,\vartheta^\ep;\tau^\ep,\sigma^\ep](I)
+\big(M^{\ep,\alpha}_{p,q}[\vae,\vue,\vartheta^\ep;\tau^\ep,\sigma^\ep](I)\big)^2\Big),
\end{align}
\begin{align}\label{est24}
&\f1\ep\|I(\ep\vae)\kappa\Delta\vthe\|^{h;\f{\beta_0}\ep}_{L^1(I;\dot{B}_{2,1}^{\f d 2 -2})}+\|I(\ep\vae)\kappa\Delta\vthe\|^{l;\f{\beta_0}\ep}_{L^1(I;\dot{B}_{2,1}^{\f d 2 -1})}\nonumber\\
&\leq \f{C_2}\ep\|I(\ep\vae)\kappa\Delta\vthe\|_{L^1(I;\dot{B}_{2,1}^{\f d 2 -2})}\leq C_2 E^\ep[\vae,\vue,\vthe](I)M^{\ep,\alpha}_{p,q}[\vae,\vue,\vartheta^\ep;\tau^\ep,\sigma^\ep](I),
\end{align}
\begin{align}\label{est25}
&\f1\ep\|\Div(\vthe\vue)\|^{h;\f{\beta_0}\ep}_{L^1(I;\dot{B}_{2,1}^{\f d 2 -2})}+\|\Div(\vthe\vue)\|^{l;\f{\beta_0}\ep}_{L^1(I;\dot{B}_{2,1}^{\f d 2 -1})}\nonumber\\
&\leq \f1\ep\|\Div(\vartheta^{\ep,h;\f{\beta_0}\ep}\vue)\|^{h;\f{\beta_0}\ep}_{L^1(I;\dot{B}_{2,1}^{\f d 2 -2})}+C_2\|\Div(\vartheta^{\ep,l;\f{\beta_0}\ep}\vue)\|^{h;\f{\beta_0}\ep}_{L^1(I;\dot{B}_{2,1}^{\f d 2 -1})}+\f{C_2}\ep\|\Div(\vartheta^{\ep,h;\f{\beta_0}\ep}\vue)\|^{l;\f{\beta_0}\ep}_{L^1(I;\dot{B}_{2,1}^{\f d 2 -2})}\nonumber\\
&\quad+\|\Div(\vartheta^{\ep,l;\f{\beta_0}\ep}\vue)\|^{l;\f{\beta_0}\ep}_{L^1(I;\dot{B}_{2,1}^{\f d 2 -1})}\leq C_2 E^\ep[\vae,\vue,\vthe](I)M^{\ep,\alpha}_{p,q}[\vae,\vue,\vartheta^\ep;\tau^\ep,\sigma^\ep](I).
\end{align}
We also have that
\begin{align}\label{est26}
&\ep\sum\limits_{2^j\geq\f{\beta_0}\ep}2^{j\f d 2}\left(\|\dot{S}_{j-2}\Div\vue\cdot\dot{ \Delta}_j\vae\|_{L^1(I;L^2)}+\|\dot{\Delta}_j F+\dot{S}_{j-2}\vue\Grad \dot{\De}_j \vae\|_{L^1(I;L^2)}\right)\nonumber\\
&\leq C_2E^\ep[\vae,\vue,\vthe](I)M^{\ep,\alpha}_{p,q}[\vae,\vue,\vartheta^\ep;\tau^\ep,\sigma^\ep](I)+C_2\ep\sum\limits_{j\in\mathbb{Z}}2^{j\f d 2}\|\dot{\Delta}_j T_{\vue}\Grad\vae-\dot{S}_{j-2}\vue\Grad \dot{\De}_j \vae\|_{L^1(I;L^2)}\nonumber\\
&\leq C_2E^\ep[\vae,\vue,\vthe](I)M^{\ep,\alpha}_{p,q}[\vae,\vue,\vartheta^\ep;\tau^\ep,\sigma^\ep](I)\nonumber\\
&\quad+C_2\ep\sum\limits_{j\in\mathbb{Z}}\sum\limits_{|j-j'|\leq2}2^{j\f d 2}\left(\|(\dot{S}_{j'-2}\vue-\dot{S}_{j-2}\vue)\Grad \dot{\De}_j \dot{\De}_{j'} \vae\|_{L^1(I;L^2)}+\|[\dot{\De}_j,\dot{S}_{j'-2}\vue]\dot{\De}_{j'} \Grad\vae\|_{L^1(I;L^2)}\right)\nonumber\\
&\leq C_2E^\ep[\vae,\vue,\vthe](I)M^{\ep,\alpha}_{p,q}[\vae,\vue,\vartheta^\ep;\tau^\ep,\sigma^\ep](I),
\end{align}
where we used the commutator estimate from lemma 2.97 of \cite{BCD} in the last inequality.  Collecting the inequalities \eqref{est20}-\eqref{est26}, we obtain that
\begin{align}\label{est33}
 &\|(F,G)\|_{L^1(I;\dot{B}_{2,1}^{\f d 2 -1})}+\f{1}{\ep}\|J\|^{h;\f{\beta_0}\ep}_{L^1(I;\dot{B}_{2,1}^{\f d 2 -2})}+\|J\|^{l;\f{\beta_0}\ep}_{L^1(I;\dot{B}_{2,1}^{\f d 2 -1})}\nonumber\\
 &\quad+\ep\sum\limits_{2^j\geq\f{\beta_0}\ep}2^{j\f d 2}\left(\|\dot{S}_{j-2}\Div\vue\cdot\dot{ \Delta}_j\vae\|_{L^1(I;L^2)}+\|\dot{\Delta}_j F+\dot{S}_{j-2}\vue\Grad \dot{\De}_j \vae\|_{L^1(I;L^2)}\right)\nonumber\\
 &\leq C_2  E^\ep[\vae,\vue,\vthe](I)\cdot \left(M^{\ep,\alpha}_{p,q}[\vae,\vue,\vartheta^\ep;\tau^\ep,\sigma^\ep](I)+\bigg(M^{\ep,\alpha}_{p,q}[\vae,\vue,\vartheta^\ep;\tau^\ep,\sigma^\ep](I)\bigg)^2\right).
\end{align}

Inserting \eqref{est33} into \eqref{est19}, we complete the proof. \ \ $\Box$

\begin{Lemma}\label{Le5.2}  Let $0<\alpha<\infty, \ep\leq\min\{\f{\beta_0}{16\alpha},1\}$ and
\[
2<p<\min{\{d,\f{2d}{d-2}\}},\ \ 0<\f1 q\leq\min{\left\{\f1 2-\f d 2\left(\f1 2-\f1 p\right),\f{d-1} 2\left(\f1 2-\f1 p\right)\right\}}.
\]
Then there exist a constant $C_3=C_3(d,p,q,\mu,\kappa,\nu)$ such that
\begin{align}
&\|(\vae,\mathbb{Q}\vue,\vthe)\|^{h;\alpha}_{N^{\ep,\alpha}_{p,q}(I)}\label{est30}\\
&\leq C_3\|(\vae,\mathbb{Q}\vue,\vthe)(t_0)\|^{h;\alpha}_{X^\ep_p}+C_3(\alpha\ep)\Big(M^{\ep,\alpha}_{p,q}[\vae,\vue,\vartheta^\ep;\tau^\ep,\sigma^\ep](I)\Big)^2\nonumber\\
&\quad+C_3 M^{\ep,\alpha}_{p,q}[\vae,\vue,\vartheta^\ep;\tau^\ep,\sigma^\ep](I)\left(M^{\ep,\alpha}_{p,q}[\vae,\vue-\vv,\vartheta^\ep;\tau^\ep,\sigma^\ep-\Theta](I)
+\|(\Theta,\vv)\|^{h;\f\alpha {64}}_{L^2(I;\dot{B}_{p,1}^{\f d p})\cap L^{q'}(I;\dot{B}_{p,1}^{\f d p-1+\f2 {q'}})}\right),\nonumber
\end{align}
\begin{align*}
 \|(\tau^\ep,\mathbb{Q}\vue)\|&^{l;\alpha}_{L^q\left(I;\dot{B}_{p,1}^{\f d p-1+\f2 q}\right)}\leq C_3(\ep\alpha)^{\f1 q}\|(\vae,\mathbb{Q}\vue,\vthe)(t_0)\|_{X^\ep_2}\\
 &\quad+C_3 (\ep\alpha)^{\f1 q}E^\ep[\vae,\vue,\vthe](I)\cdot \left(M^{\ep,\alpha}_{p,q}[\vae,\vue,\vthe](I)+\bigg(M^{\ep,\alpha}_{p,q}[\vae,\vue,\vartheta^\ep;\tau^\ep,\sigma^\ep](I)\bigg)^2\right),
\end{align*}
\begin{align*}
 &\|\sep-\Theta\|^{l;\alpha}_{L^q(I;\dot{B}^{\f d p-1+\f2 q}_{p,1})
 \cap L^{q'}(I;\dot{B}^{\f d p-1+\f2 {q'}}_{p,1})}\\
  &\leq C_3\|(\sep-\Theta)(t_0)\|^{l;\alpha}_{\dot{B}_{p,1}^{\f d p-1}}+ C_3(\ep\alpha)^{\f1 {q(q-1)}}\|(\vae,\vue,\vthe)(t_0)\|_{X^\ep_2}\\
  &\quad+C_3 (\ep\alpha)^{\f1 {q(q-1)}}E^\ep[\vae,\vue,\vthe](I)\cdot \left(M^{\ep,\alpha}_{p,q}[\vae,\vue,\vthe;\tau^\ep,\sigma^\ep](I)+\bigg(M^{\ep,\alpha}_{p,q}[\vae,\vue,\vartheta^\ep;\tau^\ep,\sigma^\ep](I)\bigg)^2\right)\\
 &\quad+ C_3(\alpha\ep)\Big(M^{\ep,\alpha}_{p,q}[\vae,\vue,\vartheta^\ep;\tau^\ep,\sigma^\ep](I)\Big)^2+C_3\|\vv\|_{L^2(I;\dot{B}^{\f d p}_{p,1})}
 \|\Theta\|^{h;\alpha}_{L^2(I;\dot{B}^{\f d p}_{p,1})}\\
  &\quad+C_3
 M^{\ep,\alpha}_{p,q}[\vae,\vue-\vv,\vartheta^\ep;\tau^\ep,\sigma^\ep-\Theta](I)
 \left(M^{\ep,\alpha}_{p,q}[\vae,\vue,\vartheta^\ep;\tau^\ep,\sigma^\ep](I)+\|\vv\|_{L^2\left(I;\dot{B}_{p,1}^{\f d p}\right)}\right.\\
 &\quad+\left.\Big(M^{\ep,\alpha}_{p,q}[\vae,\vue,\vartheta^\ep;\tau^\ep,\sigma^\ep](I)\Big)^2\right),
 \end{align*}
\begin{align*}
 &\|\mathbb{P}\vue-\vv\|_{L^\infty\left(I;\dot{B}_{p,1}^{\f d p-1}\right)\cap L^1\left(I;\dot{B}_{p,1}^{\f d p+1}\right)}\\
 &\leq C_3\|(\mathbb{P}\vue-\vv)(t_0)\|_{\dot{B}_{p,1}^{\f d p-1}}
 +C_3(\alpha\ep)\Big(M^{\ep,\alpha}_{p,q}[\vae,\vue,\vartheta^\ep;\tau^\ep,\sigma^\ep](I)\Big)^2\\
 &+C_3
 M^{\ep,\alpha}_{p,q}[\vae,\vue-\vv,\vartheta^\ep;\tau^\ep,\sigma^\ep-\Theta](I)
 \left(M^{\ep,\alpha}_{p,q}[\vae,\vue,\vartheta^\ep;\tau^\ep,\sigma^\ep](I)+\|\vv\|_{L^2\left(I;\dot{B}_{p,1}^{\f d p}\right)}\right.\\
 &\quad+\left.\Big(M^{\ep,\alpha}_{p,q}[\vae,\vue,\vartheta^\ep;\tau^\ep,\sigma^\ep](I)\Big)^2\right).
\end{align*}
\end{Lemma}
\bProof From Lemma \ref{Le3.1}, we get that
\begin{align}\label{est27}
&\|(\vae,\mathbb{Q}\vue,\vthe)\|^{h;\alpha}_{N^{\ep,\alpha}_{p,q}(I)}\nonumber\\
&\leq C_3\|(\vae,\mathbb{Q}\vue,\vthe)(t_0)\|^{h;\alpha}_{X^\ep_p}+C_3\|(F,\mathbb{Q}G)\|^{h;\f{\beta_0}\ep}_{L^1\Big(I;\dot{B}_{p,1}^{\f d p-1}\Big)}+\f {C_3} \ep\|J\|^{h;\f{\beta_0}\ep}_{L^1\Big(I;\dot{B}_{p,1}^{\f d p-2}\Big)}+C_3\|(F,\mathbb{Q}G,J)\|^{m;\alpha,\f{\beta_0}\ep}_{L^1\Big(I;\dot{B}_{2,1}^{\f d 2-1}\Big)}\nonumber\\
&\quad+C_3\ep\sum\limits_{2^j\geq\f{\beta_0}\ep}2^{j(\f d p)}\left(\|\dot{S}_{j-2}\Div\vue\cdot\dot{ \Delta}_j\vae\|_{L^1(I;L^p)}+\|\dot{\Delta}_j F+\dot{S}_{j-2}\vue\Grad \dot{\De}_j \vae\|_{L^1(I;L^p)}\right).
\end{align}

We next calculate each term in the right-hand side of \eqref{est27}, which is based on Lemma \ref{Be5}, Lemma \ref{Be6} and Proposition \ref{pro1}, and also use the fact that $\vae=\f1 {\sqrt{2}}(\tau^\ep-\sigma^\ep+\Theta-\Theta), \vthe=\f1 {\sqrt{2}}(\tau^\ep+\sigma^\ep-\Theta+\Theta)$.
\begin{itemize}
  \item {Estimates for $F:$ }
  \begin{align}\label{est28}
  &\|\Div(\vae\vue)\|^{h;\f{\beta_0}\ep}_{L^1\Big(I;\dot{B}_{p,1}^{\f d p-1}\Big)}\nonumber\\
  &\leq C_3\|\vae\|_{ L^2(I;\dot{B}_{p,1}^{\f d p})}\|\vue\|^{h;\f{\beta_0}{4\ep}}_{L^2(I;\dot{B}_{p,1}^{\f d p})}+C_3\|\vue\|_{ L^2(I;\dot{B}_{p,1}^{\f d p})}\|\vae\|^{h;\f{\beta_0}{4\ep}}_{L^2(I;\dot{B}_{p,1}^{\f d p})}+C_3\|\vae\|^{h;\f{\beta_0}{16\ep}}_{ L^2(I;\dot{B}_{p,1}^{\f d p})}\|\vue\|^{h;\f{\beta_0}{64\ep}}_{L^2(I;\dot{B}_{p,1}^{\f d p})}\nonumber\\
  &\leq C_3\|\vae\|_{ L^2(I;\dot{B}_{p,1}^{\f d p})}\left(\|\mathbb{Q}\vue\|^{h;\f{\beta_0}{4\ep}}_{L^2(I;\dot{B}_{p,1}^{\f d p})}+\|\mathbb{P}\vue-\vv\|^{h;\f{\beta_0}{4\ep}}_{L^2(I;\dot{B}_{p,1}^{\f d p})}+\|\vv\|^{h;\f{\beta_0}{4\ep}}_{L^2(I;\dot{B}_{p,1}^{\f d p})}\right)\nonumber\\
  &\quad+C_3\|\vue\|_{ L^2(I;\dot{B}_{p,1}^{\f d p})}\|\vae\|^{h;\f{\beta_0}{4\ep}}_{L^2(I;\dot{B}_{p,1}^{\f d p})}
  +C_3\|\vae\|^{h;\f{\beta_0}{16\ep}}_{ L^2(I;\dot{B}_{p,1}^{\f d p})}\|\vue\|^{h;\f{\beta_0}{64\ep}}_{L^2(I;\dot{B}_{p,1}^{\f d p})}\\
  &\leq C_3 M^{\ep,\alpha}_{p,q}[\vae,\vue,\vartheta^\ep;\tau^\ep,\sigma^\ep](I)\left(M^{\ep,\alpha}_{p,q}[\vae,\vue-\vv,\vartheta^\ep;\tau^\ep,\sigma^\ep-\Theta](I)+\|\vv\|^{h;\f{\beta_0}{4\ep}}_{L^2(I;\dot{B}_{p,1}^{\f d p})}\right),\nonumber
  \end{align}
  \begin{align*}
   &\|\Div(\vae\vue)\|^{m;\alpha,\f{\beta_0}\ep}_{L^1\Big(I;\dot{B}_{2,1}^{\f d 2-1}\Big)}\\
   &\leq C_3\|\vae\|_{ L^q(I;\dot{B}_{p,1}^{\f d p-1+\f2 q})}\|\vue\|^{m;\f\alpha 4,\f{4\beta_0}{\ep}}_{L^{q'}(I;\dot{B}_{p,1}^{\f d p-1+\f2 {q'}})}+C_3\|\vue\|_{ L^q(I;\dot{B}_{p,1}^{\f d p-1+\f2 q})}\|\vae\|^{m;\f\alpha 4,\f{4\beta_0}{\ep}}_{L^{q'}(I;\dot{B}_{p,1}^{\f d p-1+\f2 {q'}})}\\
   &\quad+C_3\|\vae\|^{h;\f{\alpha}{16}}_{ L^2(I;\dot{B}_{p,1}^{\f d p})}\|\vue\|^{h;\f{\alpha}{64}}_{L^2(I;\dot{B}_{p,1}^{\f d p})}\\
   &\leq C_3\|\vae\|_{ L^q(I;\dot{B}_{p,1}^{\f d p-1+\f2 q})}\left(\|\mathbb{Q}\vue\|^{m;\f\alpha 4,\f{4\beta_0}{\ep}}_{L^{q'}(I;\dot{B}_{p,1}^{\f d p-1+\f2 {q'}})}+\|\mathbb{P}\vue-\vv\|^{m;\f\alpha 4,\f{4\beta_0}{\ep}}_{L^{q'}(I;\dot{B}_{p,1}^{\f d p-1+\f2 {q'}})}+\|\vv\|^{m;\f\alpha 4,\f{4\beta_0}{\ep}}_{L^{q'}(I;\dot{B}_{p,1}^{\f d p-1+\f2 {q'}})}\right)\\
   &\quad+C_3\|\vue\|_{ L^q(I;\dot{B}_{p,1}^{\f d p-1+\f2 q})}\left(\|\tau^\ep\|^{m;\f\alpha 4,\alpha}_{L^{q'}(I;\dot{B}_{p,1}^{\f d p-1+\f2 {q'}})}+\|\sigma^\ep-\Theta\|^{m;\f\alpha 4,\alpha}_{L^{q'}(I;\dot{B}_{p,1}^{\f d p-1+\f2 {q'}})}+\|\Theta\|^{m;\f\alpha 4,\alpha}_{L^{q'}(I;\dot{B}_{p,1}^{\f d p-1+\f2 {q'}})}\right.\\
   &\quad\left.+\|\vae\|^{m;\alpha,\f{4\beta_0}{\ep}}_{L^{q'}(I;\dot{B}_{p,1}^{\f d p-1+\f2 {q'}})}\right)+C_3\|\vue\|^{h;\f{\alpha}{64}}_{L^2(I;\dot{B}_{p,1}^{\f d p})}\left(\|\tau^\ep\|^{m;\f\alpha {16},\alpha}_{L^{q'}(I;\dot{B}_{p,1}^{\f d p-1+\f2 {q'}})}+\|\sigma^\ep-\Theta\|^{m;\f\alpha {16},\alpha}_{L^{q'}(I;\dot{B}_{p,1}^{\f d p-1+\f2 {q'}})}\right.\\
   &\quad\left.+\|\Theta\|^{m;\f\alpha {16},\alpha}_{L^{q'}(I;\dot{B}_{p,1}^{\f d p-1+\f2 {q'}})}+\|\vae\|^{h;\alpha}_{L^{q'}(I;\dot{B}_{p,1}^{\f d p-1+\f2 {q'}})}\right)\\
   &\leq C_3 M^{\ep,\alpha}_{p,q}[\vae,\vue,\vartheta^\ep;\tau^\ep,\sigma^\ep](I)\left(M^{\ep,\alpha}_{p,q}[\vae,\vue-\vv,\vartheta^\ep;\tau^\ep,\sigma^\ep-\Theta](I)
   \right.\\
   &\quad\left.+\|\vv\|^{m;\f\alpha 4,\f{\beta_0}{4\ep}}_{L^{q'}(I;\dot{B}_{p,1}^{\f d p-1+\f2 {q'}})}+\|\Theta\|^{h;\f\alpha {16}}_{L^{q'}(I;\dot{B}_{p,1}^{\f d p-1+\f2 {q'}})\cap L^2(I;\dot{B}_{p,1}^{\f d p})}\right).
  \end{align*}
  \item {Estimates for $G:$ }
   \begin{align*}
  &\|I(\ep\vae)\f{\Grad(\vae+\vthe)}\ep\|^{h;\f{\beta_0}\ep}_{L^1\Big(I;\dot{B}_{p,1}^{\f d p-1}\Big)}\\
   &\leq C_3\|\vae\|_{ L^2(I;\dot{B}_{p,1}^{\f d p})}\left(\|\vae\|^{h;\f{\beta_0}{4\ep}}_{L^2(I;\dot{B}_{p,1}^{\f d p})}+\|\vthe\|^{m;\f{\beta_0}{4\ep},\f{\beta_0}{\ep}}_{L^2(I;\dot{B}_{p,1}^{\f d p})}\right)+C_3\ep\|\vae\|_{ L^\infty(I;\dot{B}_{p,1}^{\f d p})}\f1 \ep\|\vthe\|^{h;\f{\beta_0}{\ep}}_{L^1(I;\dot{B}_{p,1}^{\f d p})}\\
   &+C_3\|\vae+\vthe\|_{ L^2(I;\dot{B}_{p,1}^{\f d p})}\|\vae\|^{h;\f{\beta_0}{4\ep}}_{L^2(I;\dot{B}_{p,1}^{\f d p})}
   +C_3\|\vae+\vthe\|^{h;\f{\beta_0}{16\ep}}_{ L^2(I;\dot{B}_{p,1}^{\f d p})}\|\vae\|^{h;\f{\beta_0}{64\ep}}_{L^2(I;\dot{B}_{p,1}^{\f d p})}\\
   &\leq C_3 M^{\ep,\alpha}_{p,q}[\vae,\vue,\vartheta^\ep;\tau^\ep,\sigma^\ep](I)
   M^{\ep,\alpha}_{p,q}[\vae,\vue-\vv,\vartheta^\ep;\tau^\ep,\sigma^\ep-\Theta](I),
  \end{align*}
   \begin{align*}
    &\|I(\ep\vae)\f{\Grad(\vae+\vthe)}\ep\|^{m;\alpha,\f{\beta_0}\ep}_{L^1\Big(I;\dot{B}_{2,1}^{\f d 2-1}\Big)}\\
    &\leq C_3\|\vae\|_{ L^q(I;\dot{B}_{p,1}^{\f d p-1+\f2 q})}\left(\|\vae+\vthe\|^{m;\f\alpha 4,\alpha}_{L^{q'}(I;\dot{B}_{p,1}^{\f d p-1+\f2 {q'}})}+\|\vae\|^{m;\alpha ,\f{4\beta_0}{\ep}}_{L^{q'}(I;\dot{B}_{p,1}^{\f d p-1+\f2 {q'}})}+\|\vthe\|^{m;\alpha ,\f{4\beta_0}{\ep}}_{L^{q'}(I;\dot{B}_{p,1}^{\f d p-1+\f2 {q'}})}\right)\\
    &\quad+C_3\left(\|\vae+\vthe\|^{l;\alpha}_{ L^2(I;\dot{B}_{p,1}^{\f d p})}+\|\vae\|^{h;\alpha}_{ L^2(I;\dot{B}_{p,1}^{\f d p})}+\|\vthe\|^{m;\alpha,\f{\beta_0}\ep}_{ L^2(I;\dot{B}_{p,1}^{\f d p})}\right)\|\vae\|_{ L^2(I;\dot{B}_{p,1}^{\f d p})}\\
    &\quad+C_3\ep\|\vae\|_{ L^\infty(I;\dot{B}_{p,1}^{\f d p})}\f1 \ep\|\vthe\|^{h;\f{\beta_0}{\ep}}_{L^1(I;\dot{B}_{p,1}^{\f d p})}+C_3\left(\|\vae+\vthe\|^{m;\f{\alpha}{16},\alpha}_{ L^2(I;\dot{B}_{p,1}^{\f d p})}+\|\vae\|^{h;\alpha}_{ L^2(I;\dot{B}_{p,1}^{\f d p})}\right.\\
    &\quad\left.+\|\vthe\|^{m;\alpha,\f{\beta_0}\ep}_{ L^2(I;\dot{B}_{p,1}^{\f d p})}\right)\|\vae\|^{h;\f{\alpha}{64}}_{L^2(I;\dot{B}_{p,1}^{\f d p})}+C_3\ep\|\vae\|^{h;\f{\alpha}{64}}_{ L^\infty(I;\dot{B}_{p,1}^{\f d p})}\f1 \ep\|\vthe\|^{h;\f{\beta_0}{\ep}}_{L^1(I;\dot{B}_{p,1}^{\f d p})}\\
     &\leq C_3 M^{\ep,\alpha}_{p,q}[\vae,\vue,\vartheta^\ep;\tau^\ep,\sigma^\ep](I)
   M^{\ep,\alpha}_{p,q}[\vae,\vue-\vv,\vartheta^\ep;\tau^\ep,\sigma^\ep-\Theta](I),
   \end{align*}
   \begin{align*}
   &\|\Grad(\vae\vthe)\|^{h;\f{\beta_0}\ep}_{L^1\Big(I;\dot{B}_{p,1}^{\f d p-1}\Big)}\\
   &\leq  C_3\|\vae\|_{ L^2(I;\dot{B}_{p,1}^{\f d p})}\|\vthe\|^{m;\f{\beta_0}{4\ep},\f{\beta_0}{\ep}}_{L^2(I;\dot{B}_{p,1}^{\f d p})}+C_3\ep\|\vae\|_{  L^\infty(I;\dot{B}_{p,1}^{\f d p})}\f{1}\ep\|\vthe\|^{h;\f{\beta_0}{\ep}}_{L^1(I;\dot{B}_{p,1}^{\f d p})}+C_3\|\vthe\|_{ L^2(I;\dot{B}_{p,1}^{\f d p})}\|\vae\|^{h;\f{\beta_0}{4\ep}}_{L^2(I;\dot{B}_{p,1}^{\f d p})}\\
   &\quad+C_3\|\vae\|^{h;\f{\beta_0}{16\ep}}_{ L^2(I;\dot{B}_{p,1}^{\f d p})}\|\vthe\|^{h;\f{\beta_0}{64\ep},\f{\beta_0}{\ep}}_{L^2(I;\dot{B}_{p,1}^{\f d p})}+C_3\ep\|\vae\|^{h;\f{\beta_0}{16\ep}}_{  L^\infty(I;\dot{B}_{p,1}^{\f d p})}\f{1}\ep\|\vthe\|^{h;\f{\beta_0}{\ep}}_{L^1(I;\dot{B}_{p,1}^{\f d p})}\\
   &\leq C_3 M^{\ep,\alpha}_{p,q}[\vae,\vue,\vartheta^\ep;\tau^\ep,\sigma^\ep](I)
   M^{\ep,\alpha}_{p,q}[\vae,\vue-\vv,\vartheta^\ep;\tau^\ep,\sigma^\ep-\Theta](I),
   \end{align*}
    \begin{align*}
     &\|\Grad(\vae\vthe)\|^{m;\alpha,\f{\beta_0}\ep}_{L^1\Big(I;\dot{B}_{2,1}^{\f d 2-1}\Big)}\\
     &=\|\Grad(T_{\vae}\vthe+T_{\vartheta^{\ep,l;\f{\beta_0}\ep}}\vae+T_{\vartheta^{\ep,h;\f{\beta_0}\ep}}\vae
     +R(\vartheta^{\ep,l;\f{\beta_0}\ep},\vae)+R(\vartheta^{\ep,h;\f{\beta_0}\ep},\vae))\|^{m;\alpha,\f{\beta_0}\ep}_{L^1\Big(I;\dot{B}_{2,1}^{\f d 2-1}\Big)}\\
     &\leq C_3\|\vae\|_{ L^q(I;\dot{B}_{p,1}^{\f d p-1+\f2 q})}\left(\|\tau^\ep\|^{m;\f\alpha 4,\alpha}_{L^{q'}(I;\dot{B}_{p,1}^{\f d p-1+\f2 {q'}})}+\|\sigma^\ep-\Theta\|^{m;\f\alpha 4,\alpha}_{L^{q'}(I;\dot{B}_{p,1}^{\f d p-1+\f2 {q'}})}\right.\\
     &\quad\left.+\|\Theta\|^{m;\f\alpha 4,\alpha}_{L^{q'}(I;\dot{B}_{p,1}^{\f d p-1+\f2 {q'}})}+\|\vthe\|^{m;\alpha,\f{4\beta_0}\ep}_{L^{q'}(I;\dot{B}_{p,1}^{\f d p-1+\f2 {q'}})}\right)+C_3\left(\|\vartheta^{\ep,l;\f{\beta_0}\ep}\|_{ L^q(I;\dot{B}_{p,1}^{\f d p-1+\f2 q})}\right.\\
     &\left.+\f1 \ep\|\vartheta^{\ep,h;\f{\beta_0}\ep}\|_{ L^q(I;\dot{B}_{p,1}^{\f d p-2+\f2 q})}\right)\cdot\left(\|\tau^\ep\|^{m;\f\alpha 4,\alpha}_{L^{q'}(I;\dot{B}_{p,1}^{\f d p-1+\f2 {q'}})}+\|\sigma^\ep-\Theta\|^{m;\f\alpha 4,\alpha}_{L^{q'}(I;\dot{B}_{p,1}^{\f d p-1+\f2 {q'}})}\right.\\
     &\quad\left.+\|\Theta\|^{m;\f\alpha 4,\alpha}_{L^{q'}(I;\dot{B}_{p,1}^{\f d p-1+\f2 {q'}})}+\|\vae\|^{m;\alpha,\f{4\beta_0}\ep}_{L^{q'}(I;\dot{B}_{p,1}^{\f d p-1+\f2 {q'}})}\right)+C_3\left(\|\tau^\ep\|^{m;\f\alpha {16},\alpha}_{L^2(I;\dot{B}_{p,1}^{\f d p})}+\|\sigma^\ep-\Theta\|^{m;\f\alpha {16},\alpha}_{L^2(I;\dot{B}_{p,1}^{\f d p})}\right.\\
     &\quad\left.+\|\Theta\|^{m;\f\alpha {16},\alpha}_{L^2(I;\dot{B}_{p,1}^{\f d p})}+\|\vae\|^{h;\alpha}_{L^2(I;\dot{B}_{p,1}^{\f d p})}\right)\left(\|\vartheta^{\ep,l;\f{\beta_0}\ep}\|^{h;\f{\alpha}{64}}_{L^2(I;\dot{B}_{p,1}^{\f d p})}+\f1 \ep\|\vartheta^{\ep,h;\f{\beta_0}\ep}\|^{h;\f{\alpha}{64}}_{L^2(I;\dot{B}_{p,1}^{\f d p-1})}\right)\\
      &\leq C_3 M^{\ep,\alpha}_{p,q}[\vae,\vue,\vartheta^\ep;\tau^\ep,\sigma^\ep](I)\left(
   M^{\ep,\alpha}_{p,q}[\vae,\vue-\vv,\vartheta^\ep;\tau^\ep,\sigma^\ep-\Theta](I)+\|\Theta\|^{m;\f\alpha {16},\alpha}_{L^2(I;\dot{B}_{p,1}^{\f d p})\cap L^{q'}(I;\dot{B}_{p,1}^{\f d p-1+\f2 {q'}})}\right),
    \end{align*}
    \begin{align*}
   &\|I(\ep\vae)\Grad(\vae\vthe)\|^{h;\f{\beta_0}\ep}_{L^1\Big(I;\dot{B}_{p,1}^{\f d p-1}\Big)}\\
   &\leq \|\ep\vae\|_{L^\infty\Big(I;\dot{B}_{p,1}^{\f d p}\Big)} \|\Grad(\vae\vthe)\|_{L^1\Big(I;\dot{B}_{p,1}^{\f d p-1}\Big)}\\
   &\leq \|\ep\vae\|_{L^\infty\Big(I;\dot{B}_{p,1}^{\f d p}\Big)} \left(\|\vae\|_{ L^2(I;\dot{B}_{p,1}^{\f d p})}\|\vthe\|^{l;\f{\beta_0}{\ep}}_{L^2(I;\dot{B}_{p,1}^{\f d p})}+C_3\ep\|\vae\|_{  L^\infty(I;\dot{B}_{p,1}^{\f d p})}\f{1}\ep\|\vthe\|^{h;\f{\beta_0}{\ep}}_{L^1(I;\dot{B}_{p,1}^{\f d p})}\right)\\
    &\leq C_3 M^{\ep,\alpha}_{p,q}[\vae,\vue,\vartheta^\ep;\tau^\ep,\sigma^\ep](I)
   M^{\ep,\alpha}_{p,q}[\vae,\vue-\vv,\vartheta^\ep;\tau^\ep,\sigma^\ep-\Theta](I),
   \end{align*}
     \begin{align*}
   &\|I(\ep\vae)\Grad(\vae\vthe)\|^{m;\alpha,\f{\beta_0}\ep}_{L^1\Big(I;\dot{B}_{2,1}^{\f d 2-1}\Big)}\\
   &\leq \|\ep\vae\|_{L^\infty\Big(I;\dot{B}_{p,1}^{\f d p}\Big)}\|\Grad(\vae\vthe)\|^{m;\f\alpha 4,\f{4\beta_0}\ep}_{L^1\Big(I;\dot{B}_{2,1}^{\f d 2-1}\Big)}+ \|\ep\vae\|_{L^\infty\Big(I;\dot{B}_{p,1}^{\f d p}\Big)} \|\Grad(\vae\vthe)\|_{L^1\Big(I;\dot{B}_{p,1}^{\f d p-1}\Big)}\\
       &\leq C_3 \left(M^{\ep,\alpha}_{p,q}[\vae,\vue,\vartheta^\ep;\tau^\ep,\sigma^\ep](I)
       +\Big(M^{\ep,\alpha}_{p,q}[\vae,\vue,\vartheta^\ep;\tau^\ep,\sigma^\ep](I)\Big)^2\right)\\
   &\quad\cdot\left(
   M^{\ep,\alpha}_{p,q}[\vae,\vue-\vv,\vartheta^\ep;\tau^\ep,\sigma^\ep-\Theta](I)+\|\Theta\|^{m;\f\alpha {64},\alpha}_{L^2(I;\dot{B}_{p,1}^{\f d p})\cap L^{q'}(I;\dot{B}_{p,1}^{\f d p-1+\f2 {q'}})}\right),
   \end{align*}
   \begin{align*}
   &\|I(\ep\vae)\mathcal{L}\vue\|^{h;\f{\beta_0}\ep}_{L^1\Big(I;\dot{B}_{p,1}^{\f d p-1}\Big)}\\
   &\leq \|I(\ep\vae)\mathcal{L}\vue\|_{L^1\Big(I;\dot{B}_{p,1}^{\f d p-1}\Big)}\\
   &\leq C_3(\alpha\ep)\|\vae\|_{ L^2(I;\dot{B}_{p,1}^{\f d p})}\|\vue\|^{l;\alpha}_{L^2(I;\dot{B}_{p,1}^{\f d p})}+\ep\|\vae\|_{ L^\infty(I;\dot{B}_{p,1}^{\f d p})}\f1\ep\|\vue\|^{h;\alpha}_{L^1(I;\dot{B}_{p,1}^{\f d p+1})}\\
   &\leq C_3(\alpha\ep)\Big(M^{\ep,\alpha}_{p,q}[\vae,\vue,\vartheta^\ep;\tau^\ep,\sigma^\ep](I)\Big)^2+M^{\ep,\alpha}_{p,q}[\vae,\vue,\vartheta^\ep;\tau^\ep,\sigma^\ep](I)
   M^{\ep,\alpha}_{p,q}[\vae,\vue-\vv,\vartheta^\ep;\tau^\ep,\sigma^\ep-\Theta](I),
   \end{align*}
   \begin{align*}
     &\|I(\ep\vae)\mathcal{L}\vue\|^{m;\alpha,\f{\beta_0}\ep}_{L^1\Big(I;\dot{B}_{2,1}^{\f d 2-1}\Big)}\\
     &\leq C_3\|\vae\|_{ L^q(I;\dot{B}_{p,1}^{\f d p-1+\f2 q})}\|\vue\|^{m;\f\alpha 4, \f{4\beta_0}\ep}_{ L^{q'}(I;\dot{B}_{p,1}^{\f d p-1+\f2 {q'}})}+C_3(\alpha\ep)\|\vae\|_{ L^2(I;\dot{B}_{p,1}^{\f d p})}\|\vue\|^{l;\alpha}_{L^2(I;\dot{B}_{p,1}^{\f d p})}\\
     &\quad+\ep\|\vae\|_{ L^\infty(I;\dot{B}_{p,1}^{\f d p})}\f1\ep\|\vue\|^{h;\alpha}_{L^1(I;\dot{B}_{p,1}^{\f d p+1})}\\
     &\leq C_3 M^{\ep,\alpha}_{p,q}[\vae,\vue,\vartheta^\ep;\tau^\ep,\sigma^\ep](I)\left(  M^{\ep,\alpha}_{p,q}[\vae,\vue-\vv,\vartheta^\ep;\tau^\ep,\sigma^\ep-\Theta](I)+\|\vv\|^{m;\f\alpha 4, \f{4\beta_0}\ep}_{ L^{q'}(I;\dot{B}_{p,1}^{\f d p-1+\f2 {q'}})}\right)\\
     &\quad+ C_3(\alpha\ep)\Big(M^{\ep,\alpha}_{p,q}[\vae,\vue,\vartheta^\ep;\tau^\ep,\sigma^\ep](I)\Big)^2,
   \end{align*}
   \begin{align*}
     &\|\vue\cdot\Grad\vue\|^{h;\f{\beta_0}\ep}_{L^1\Big(I;\dot{B}_{p,1}^{\f d p-1}\Big)}+\|\vue\cdot\Grad\vue\|^{m;\alpha,\f{\beta_0}\ep}_{L^1\Big(I;\dot{B}_{2,1}^{\f d 2-1}\Big)}\\
     &\leq C_3 \|\vue\cdot\Grad\vue\|^{h;\alpha}_{L^1\Big(I;\dot{B}_{2,1}^{\f d 2-1}\Big)}\\
     &\leq C_3\|\vue\|_{ L^q(I;\dot{B}_{p,1}^{\f d p-1+\f2 q})}\|\vue\|^{h,\f{\alpha}4}_{ L^{q'}(I;\dot{B}_{p,1}^{\f d p-1+\f2 {q'}})}+C_3\|\vue\|^{h,\f{\alpha}{64}}_{ L^q(I;\dot{B}_{p,1}^{\f d p-1+\f2 q})}\|\vue\|^{h,\f{\alpha}{16}}_{ L^{q'}(I;\dot{B}_{p,1}^{\f d p-1+\f2 {q'}})}\\
     &\leq C_3 M^{\ep,\alpha}_{p,q}[\vae,\vue,\vartheta^\ep;\tau^\ep,\sigma^\ep](I)\left(  M^{\ep,\alpha}_{p,q}[\vae,\vue-\vv,\vartheta^\ep;\tau^\ep,\sigma^\ep-\Theta](I)+\|\vv\|^{h;\f\alpha{16}}_{ L^{q'}(I;\dot{B}_{p,1}^{\f d p-1+\f2 {q'}})}\right).
   \end{align*}
  \item {Estimates for $J:$ }
  \begin{align*}
   &\|2\mu|D\vue|^2+\lambda|\Div\vue|^2\|^{h;\f{\beta_0}\ep}_{L^1\Big(I;\dot{B}_{p,1}^{\f d p-2}\Big)}+\ep\|2\mu|D\vue|^2+\lambda|\Div\vue|^2\|^{m;\alpha,\f{\beta_0}\ep}_{L^1\Big(I;\dot{B}_{2,1}^{\f d 2-1}\Big)}\\
   &\leq C_3\|2\mu|D\vue|^2+\lambda|\Div\vue|^2\|^{h;\alpha}_{L^1\Big(I;\dot{B}_{2,1}^{\f d 2-2}\Big)}\\
    &\leq C_3\|\vue\|_{ L^q(I;\dot{B}_{p,1}^{\f d p-1+\f2 q})}\|\vue\|^{h,\f{\alpha}4}_{ L^{q'}(I;\dot{B}_{p,1}^{\f d p-1+\f2 {q'}})}+C_3\|\vue\|^{h,\f{\alpha}{64}}_{ L^q(I;\dot{B}_{p,1}^{\f d p-1+\f2 q})}\|\vue\|^{h,\f{\alpha}{16}}_{ L^{q'}(I;\dot{B}_{p,1}^{\f d p-1+\f2 {q'}})}\\
     &\leq C_3 M^{\ep,\alpha}_{p,q}[\vae,\vue,\vartheta^\ep;\tau^\ep,\sigma^\ep](I)\left(  M^{\ep,\alpha}_{p,q}[\vae,\vue-\vv,\vartheta^\ep;\tau^\ep,\sigma^\ep-\Theta](I)+\|\vv\|^{h;\f\alpha{16}}_{ L^{q'}(I;\dot{B}_{p,1}^{\f d p-1+\f2 {q'}})}\right),
  \end{align*}
  \begin{align*}
  &\|I(\ep\vae)\big(2\mu|D\vue|^2+\lambda|\Div\vue|^2\big)\|^{h;\f{\beta_0}\ep}_{L^1\Big(I;\dot{B}_{p,1}^{\f d p-2}\Big)}+\ep\|I(\ep\vae)\big(2\mu|D\vue|^2+\lambda|\Div\vue|^2\big)\|^{m;\alpha,\f{\beta_0}\ep}_{L^1\Big(I;\dot{B}_{2,1}^{\f d 2-1}\Big)}\\
  &\leq C_3\ep\|\vae\|_{ L^\infty(I;\dot{B}_{p,1}^{\f d p})}\|2\mu|D\vue|^2+\lambda|\Div\vue|^2\|_{L^1\Big(I;\dot{B}_{2,1}^{\f d 2-2}\Big)}\\
  & \leq C_3 \Big(M^{\ep,\alpha}_{p,q}[\vae,\vue,\vartheta^\ep;\tau^\ep,\sigma^\ep](I)\Big)^2 M^{\ep,\alpha}_{p,q}[\vae,\vue-\vv,\vartheta^\ep;\tau^\ep,\sigma^\ep-\Theta](I),
  \end{align*}

\begin{align*}
&\f1\ep\|I(\ep\vae)\kappa\Delta\vthe\|^{h;\f{\beta_0}\ep}_{L^1\Big(I;\dot{B}_{p,1}^{\f d p-2}\Big)}\\
&\leq C_3\|I(\ep\vae)\kappa\Delta\vartheta^{\ep,l;\alpha}\|_{L^1\Big(I;\dot{B}_{p,1}^{\f d p-1}\Big)}+\f{C_3}\ep\|I(\ep\vae)\kappa\Delta\vartheta^{\ep,h;\alpha}\|_{L^1\Big(I;\dot{B}_{p,1}^{\f d p-2}\Big)}\\
&\leq C_3\left((\alpha\ep)\|\vae\|_{ L^2(I;\dot{B}_{p,1}^{\f d p})}\|\vthe\|^{l;\alpha}_{ L^2(I;\dot{B}_{p,1}^{\f d p})}+\|\vae\|_{ L^2(I;\dot{B}_{p,1}^{\f d p})}\|\vthe\|^{m;\alpha,\f{\beta_0}\ep}_{ L^2(I;\dot{B}_{p,1}^{\f d p})}+\ep\|\vae\|_{ L^\infty(I;\dot{B}_{p,1}^{\f d p})}\f1\ep\|\vthe\|^{h;\f{\beta_0}\ep}_{ L^1(I;\dot{B}_{p,1}^{\f d p})}\right)\\
&\leq C_3(\alpha\ep)\Big(M^{\ep,\alpha}_{p,q}[\vae,\vue,\vartheta^\ep;\tau^\ep,\sigma^\ep](I)\Big)^2\\
&\quad+ M^{\ep,\alpha}_{p,q}[\vae,\vue-\vv,\vartheta^\ep;\tau^\ep,\sigma^\ep-\Theta](I)
M^{\ep,\alpha}_{p,q}[\vae,\vue,\vartheta^\ep;\tau^\ep,\sigma^\ep](I),
\end{align*}
\begin{align*}
&\|I(\ep\vae)\kappa\Delta\vthe\|^{m;\alpha,\f{\beta_0}\ep}_{L^1\Big(I;\dot{B}_{2,1}^{\f d 2-1}\Big)}\\
&\leq C_3\|\vae\|_{ L^q(I;\dot{B}_{p,1}^{\f d p-1+\f2 q})}\|\vthe\|^{m;\f\alpha4,\f{4\beta_0}\ep}_{ L^{q'}(I;\dot{B}_{p,1}^{\f d p-1+\f2 {q'}})}+\|T_{\kappa\Delta\vartheta^{\ep,l;\alpha}}I(\ep\vae)\|_{L^1\Big(I;\dot{B}_{2,1}^{\f d 2-1}\Big)}\\
&\quad+\f{C_3}\ep\|T_{\kappa\Delta\vartheta^{\ep,h;\alpha}}I(\ep\vae)\|_{L^1\Big(I;\dot{B}_{2,1}^{\f d 2-2}\Big)}+\|R\big(\kappa\Delta\vartheta^{\ep,l;\alpha},I(\ep\vae)\big)\|_{L^1\Big(I;\dot{B}_{2,1}^{\f d 2-1}\Big)}\\
&\quad+\|R\big(\kappa\Delta\vartheta^{\ep,h;\alpha},I(\ep\vae)\big)\|_{L^1\Big(I;\dot{B}_{2,1}^{\f d 2-1}\Big)}\\
&\leq C_3(\alpha\ep)\Big(M^{\ep,\alpha}_{p,q}[\vae,\vue,\vartheta^\ep;\tau^\ep,\sigma^\ep](I)\Big)^2
+C_3M^{\ep,\alpha}_{p,q}[\vae,\vue,\vartheta^\ep;\tau^\ep,\sigma^\ep](I)\left(\|\Theta\|^{m;\f\alpha4,\f{4\beta_0}\ep}_{ L^{q'}(I;\dot{B}_{p,1}^{\f d p-1+\f2 {q'}})} \right.\\
&\quad \left. +M^{\ep,\alpha}_{p,q}[\vae,\vue-\vv,\vartheta^\ep;\tau^\ep,\sigma^\ep-\Theta](I)\right)
\end{align*}
\begin{align*}
&\f1\ep\|\Div(\vthe\vue)\|^{h;\f{\beta_0}\ep}_{L^1\Big(I;\dot{B}_{p,1}^{\f d p-2}\Big)}\\
&\leq \f{C_3}\ep\|\vartheta^{\ep,h;\f{\beta_0}\ep}\|_{L^2\Big(I;\dot{B}_{p,1}^{\f d p-1}\Big)}\|\vue\|^{h;\f{\beta_0}{4\ep}}_{L^2\Big(I;\dot{B}_{p,1}^{\f d p}\Big)}+C_3\|\vartheta^{\ep,l;\f{\beta_0}\ep}\|_{L^2\Big(I;\dot{B}_{p,1}^{\f d p}\Big)}\|\vue\|^{h;\f{\beta_0}{4\ep}}_{L^2\Big(I;\dot{B}_{p,1}^{\f d p}\Big)}\\
&\quad+C_3\|\vue\|_{L^2\Big(I;\dot{B}_{p,1}^{\f d p}\Big)}\left(\|\vthe\|^{m;\f{\beta_0}{4\ep},\f{\beta_0}\ep}_{L^2\Big(I;\dot{B}_{p,1}^{\f d p}\Big)}+\f1\ep\|\vthe\|^{h;\f{\beta_0}\ep}_{L^2\Big(I;\dot{B}_{p,1}^{\f d p-1}\Big)}\right)\\
&\quad+C_3\|\vue\|^{\f{\beta_0}{64\ep}}_{L^2\Big(I;\dot{B}_{p,1}^{\f d p}\Big)}\left(\|\vthe\|^{m;\f{\beta_0}{16\ep},\f{\beta_0}\ep}_{L^2\Big(I;\dot{B}_{p,1}^{\f d p}\Big)}+\f1\ep\|\vthe\|^{h;\f{\beta_0}\ep}_{L^2\Big(I;\dot{B}_{p,1}^{\f d p-1}\Big)}\right)\\
&\leq C_3M^{\ep,\alpha}_{p,q}[\vae,\vue,\vartheta^\ep;\tau^\ep,\sigma^\ep](I)
\left(M^{\ep,\alpha}_{p,q}[\vae,\vue-\vv,\vartheta^\ep;\tau^\ep,\sigma^\ep-\Theta](I)+\|\vv\|^{h;\f{\beta_0}{4\ep}}_{L^2\Big(I;\dot{B}_{p,1}^{\f d p}\Big)}\right),
\end{align*}
\begin{align*}
 &\|\Div(\vthe\vue)\|^{m;\alpha,\f{\beta_0}\ep}_{L^1\Big(I;\dot{B}_{2,1}^{\f d 2-1}\Big)}\\
 &\leq \f{C_3}\ep\|\vartheta^{\ep,h;\f{\beta_0}\ep}\|_{L^q\Big(I;\dot{B}_{p,1}^{\f d p-2+\f2 q}\Big)}\|\vue\|^{m;\f\alpha 4,\f{4\beta_0}{\ep}}_{L^{q'}\Big(I;\dot{B}_{p,1}^{\f d p-1+\f2{q'}}\Big)}+C_3\|\vartheta^{\ep,l;\f{\beta_0}\ep}\|_{L^q\Big(I;\dot{B}_{p,1}^{\f d p-1+\f2 q}\Big)}\|\vue\|^{m;\f\alpha 4,\f{4\beta_0}{\ep}}_{L^{q'}\Big(I;\dot{B}_{p,1}^{\f d p-1+\f2{q'}}\Big)}\\
 &\quad+C_3\|\vue\|_{L^q\Big(I;\dot{B}_{p,1}^{\f d p-1+\f2 q}\Big)}\left(\|\vthe\|^{m;\f\alpha 4,\alpha}_{L^{q'}\Big(I;\dot{B}_{p,1}^{\f d p-1+\f2{q'}}\Big)}+\|\vthe\|^{m;\alpha,\f{4\beta_0}\ep}_{L^{q'}\Big(I;\dot{B}_{p,1}^{\f d p-1+\f2{q'}}\Big)}\right)\\
 &\quad+C_3\|\vue\|^{\f\alpha{64}}_{L^q\Big(I;\dot{B}_{p,1}^{\f d p-1+\f2 q}\Big)}\left(\|\vthe\|^{m;\f{\alpha}{16},\alpha}_{L^{q'}\Big(I;\dot{B}_{p,1}^{\f d p-1+\f2{q'}}\Big)}+\|\vthe\|^{m;\alpha,\f{\beta_0}\ep}_{L^2\Big(I;\dot{B}_{p,1}^{\f d p-1+\f2{q'}}\Big)}+\f1\ep\|\vthe\|^{h;\f{\beta_0}\ep}_{L^{q'}\Big(I;\dot{B}_{p,1}^{\f d p-2+\f2{q'}}\Big)}\right)\\
 &\leq C_3M^{\ep,\alpha}_{p,q}[\vae,\vue,\vartheta^\ep;\tau^\ep,\sigma^\ep](I)
\left(M^{\ep,\alpha}_{p,q}[\vae,\vue-\vv,\vartheta^\ep;\tau^\ep,\sigma^\ep-\Theta](I)
\right.\\
&\quad\left.+\|\vv\|^{m;\f\alpha 4,\f{4\beta_0}{\ep}}_{L^{q'}\Big(I;\dot{B}_{p,1}^{\f d p-1+\f2{q'}}\Big)}+\|\Theta\|^{m;\f\alpha {16},\alpha}_{L^{q'}\Big(I;\dot{B}_{p,1}^{\f d p-1+\f2{q'}}\Big)}\right).
  \end{align*}
  \item {Estimates for $\ep\sum\limits_{2^j\geq\f{\beta_0}\ep}2^{j(\f d p)}\left(\|\dot{S}_{j-2}\Div\vue\cdot\dot{ \Delta}_j\vae\|_{L^1(I;L^p)}+\|\dot{\Delta}_j F+\dot{S}_{j-2}\vue\Grad \dot{\De}_j \vae\|_{L^1(I;L^p)}\right):$  }
      Calculating similarly to \eqref{est28} and using the commutator estimate from lemma 2.97 of \cite{BCD}, we have that
      \begin{align}\label{est29}
&\ep\sum\limits_{2^j\geq\f{\beta_0}\ep}2^{j(\f d p)}\left(\|\dot{S}_{j-2}\Div\vue\cdot\dot{ \Delta}_j\vae\|_{L^1(I;L^p)}+\|\dot{\Delta}_j F+\dot{S}_{j-2}\vue\Grad \dot{\De}_j \vae\|_{L^1(I;L^p)}\right)\nonumber\\
&\leq C_3 \|\vue\|^{l;\f{\beta_0}\ep}_{L^2\Big(I;\dot{B}_{p,1}^{\f d p}\Big)}\|\vae\|^{h;\f{\beta_0}\ep}_{L^2\Big(I;\dot{B}_{p,1}^{\f d p}\Big)}+C_3 \|\vue\|^{h;\f{\beta_0}\ep}_{L^1\Big(I;\dot{B}_{p,1}^{\f d p+1}\Big)}\ep\|\vae\|^{h;\f{\beta_0}\ep}_{L^\infty\Big(I;\dot{B}_{p,1}^{\f d p}\Big)}\\
&\quad+C_3 \ep\|T_{\Grad\vae}\vue+R(\Grad\vae,\vue)\|^{h;\f{\beta_0}\ep}_{L^1\Big(I;\dot{B}_{p,1}^{\f d p}\Big)}\nonumber\\
&\quad+C_3\ep\sum\limits_{2^j\geq\f{\beta_0}\ep}\sum\limits_{|j-j'|\leq2}2^{j\f d 2}\left(\|(\dot{S}_{j'-2}\vue-\dot{S}_{j-2}\vue)\Grad \dot{\De}_j \dot{\De}_{j'} \vae\|_{L^1(I;L^2)}+\|[\dot{\De}_j,\dot{S}_{j'-2}\vue]\dot{\De}_{j'} \Grad\vae\|_{L^1(I;L^2)}\right)\nonumber\\
&\leq C_3 M^{\ep,\alpha}_{p,q}[\vae,\vue,\vartheta^\ep;\tau^\ep,\sigma^\ep](I)\left(M^{\ep,\alpha}_{p,q}[\vae,\vue-\vv,\vartheta^\ep;\tau^\ep,\sigma^\ep-\Theta](I)+\|\vv\|^{h;\f{\beta_0}{4\ep}}_{L^2(I;\dot{B}_{p,1}^{\f d p})}\right).\nonumber
\end{align}
\end{itemize}
Inserting the estimates of $F,G,J$ and \eqref{est29} into \eqref{est27}, we arrive at \eqref{est30}.

\vspace{5mm}

Noting that $ \ep\leq\f{\beta_0}{16\alpha}$ then using Lemma \ref{Be9} and \eqref{est33} yields
\begin{align}
 \|(\tau^\ep,\mathbb{Q}\vue)\|^{l;\alpha}_{L^q\left(I;\dot{B}_{p,1}^{\f d p-1+\f2 q}\right)}&\leq C_3(\ep\alpha)^{\f1 q}\left(\|(\vae,\mathbb{Q}\vue,\vthe)(t_0)\|^{l;\alpha}_{\dot{B}_{2,1}^{\f d 2-1}}+\|(F,\mathbb{Q}G,J)\|^{l;\alpha}_{L^1(I;\dot{B}_{2,1}^{\f d 2-1})}\right)\nonumber\\
 &\leq C_3(\ep\alpha)^{\f1 q}\|(\vae,\mathbb{Q}\vue,\vthe)(t_0)\|_{X^\ep_2}+C_3 (\ep\alpha)^{\f1 q}E^\ep[\vae,\vue,\vthe](I)\nonumber\\
 &\quad\cdot \left(M^{\ep,\alpha}_{p,q}[\vae,\vue,\vthe](I)+\bigg(M^{\ep,\alpha}_{p,q}[\vae,\vue,\vartheta^\ep;\tau^\ep,\sigma^\ep](I)\bigg)^2\right).\label{est31}
\end{align}

\vspace{5mm}

We see that
\begin{align}\label{equ11}
 &\p_t(\sep-\Theta)-\f{\kappa}2\Delta(\sep-\Theta)=K,
\end{align}
 where
 \begin{align}\label{equ12}
   K&=\f{\kappa}2\Delta\tep-\Div(\sep\mathbb{Q}\vue)-\f{\sqrt{2}\kappa}2 I(\ep\vae)\Delta\vthe+\f{\sqrt{2}\ep}2 (1-I(\ep\vae))\big(2\mu|D\vue|^2+\lambda|\Div\vue|^2\big)\nonumber\\
   &\quad-(\mathbb{P}\vue-\vv)\Grad\sigma^\ep-\vv\Grad(\sigma^\ep-\Theta)
      \end{align}
From Lemma \ref{Be3}, Lemma \ref{Le5.1} and \eqref{est31},
\begin{align*}
 &\|\sep-\Theta\|^{l;\alpha}_{L^q(I;\dot{B}^{\f d p-1+\f2 q}_{p,1})
 \cap L^{q'}(I;\dot{B}^{\f d p-1+\f2 {q'}}_{p,1})}\\
 &\leq C_3\|(\sep-\Theta)(t_0)\|^{l;\alpha}_{\dot{B}_{p,1}^{\f d p-1}}+C_3\|\tau^\ep\|^{l;\alpha}_{L^{q'}(I;\dot{B}^{\f d p-1+\f2 {q'}}_{p,1})}+C_3(\alpha\ep)\|(1-I(\ep\vae))\big(2\mu|D\vue|^2+\lambda|\Div\vue|^2\big)\|^{l;\alpha}_{L^1(I;\dot{B}^{\f d p-2}_{p,1})}\\
 &\quad+C_3\|\Big(\Div(\sigma^{\ep,l;\f{\beta_0}\ep}\mathbb{Q}\vue),
 (\mathbb{P}\vue-\vv)\Grad\sigma^{\ep,l;\f{\beta_0}\ep},\vv\Grad(\sep-\Theta)^{l;\f{\beta_0}\ep}\Big)\|^{l;\alpha}_{L^1(I;\dot{B}^{\f d p-1}_{p,1})}\\
 &\quad+C_3\|\Big(\Div(a^{\ep,h;\f{\beta_0}\ep}\mathbb{Q}\vue),
 (\mathbb{P}\vue-\vv)\Grad a^{\ep,h;\f{\beta_0}\ep},\vv\Grad a^{\ep,h;\f{\beta_0}\ep},\vv\Grad\Theta^{h;\f{\beta_0}\ep}\Big)\|^{l;\alpha}_{L^1(I;\dot{B}^{\f d p-1}_{p,1})}\\
 &\quad +C_3\alpha\|\Big(\Div(\vartheta^{\ep,h;\f{\beta_0}\ep}\mathbb{Q}\vue),
 (\mathbb{P}\vue-\vv)\Grad \vartheta^{\ep,h;\f{\beta_0}\ep},\vv\Grad\vartheta^{\ep,h;\f{\beta_0}\ep}\Big)\|^{l;\alpha}_{L^1(I;\dot{B}^{\f d p-2}_{p,1})}+C_3\alpha\|I(\ep\vae)\Delta\vartheta^{\ep}\|^{l;\alpha}_{L^1(I;\dot{B}^{\f d p-2}_{p,1})}\\
 &\leq C_3\|(\sep-\Theta)(t_0)\|^{l;\alpha}_{\dot{B}_{p,1}^{\f d p-1}}+C_3\left(\|\tau^\ep\|^{l;\alpha}_{L^{q}(I;\dot{B}^{\f d p-1+\f2 {q}}_{p,1})}\right)^{\f1 {q-1}}\left(\|\tau^\ep\|^{l;\alpha}_{L^1(I;\dot{B}^{\f d p+1}_{p,1})}\right)^{\f{q-2} {q-1}}+C_3(\alpha\ep)\big(1+\ep\|\vae\|_{L^\infty(I;\dot{B}^{\f d p}_{p,1})}\big)\\
 &\quad \cdot\|\vue\|^2_{L^2(I;\dot{B}^{\f d p}_{p,1})}+C_3\left(\|\sigma^\ep\|^{l;\f{\beta_0}\ep}_{L^2(I;\dot{B}^{\f d p}_{p,1})}+\|\vae\|^{h;\f{\beta_0}\ep}_{L^2(I;\dot{B}^{\f d p}_{p,1})}+(\alpha\ep)\f1\ep\|\sigma^\ep\|^{h;\f{\beta_0}\ep}_{L^2(I;\dot{B}^{\f d p-1}_{p,1})}
 \right)\cdot\|(\mathbb{Q}\vue,\mathbb{P}\vue-\vv)\|_{L^2(I;\dot{B}^{\f d p}_{p,1})}\\
 &\quad+C_3\|\vv\|_{L^2(I;\dot{B}^{\f d p}_{p,1})}\left(\|\sigma^\ep-\Theta\|^{l;\alpha}_{L^2(I;\dot{B}^{\f d p}_{p,1})}+\|\sigma^\ep-\Theta\|^{m;\alpha,\f{\beta_0}\ep}_{L^2(I;\dot{B}^{\f d p}_{p,1})}+\|(\vae,\Theta)\|^{h;\f{\beta_0}\ep}_{L^2(I;\dot{B}^{\f d p}_{p,1})}+(\alpha\ep)\f1\ep\|\vthe\|^{h;\f{\beta_0}\ep}_{L^2(I;\dot{B}^{\f d p-1}_{p,1})}\right)\\
 &\quad+C_3(\alpha\ep)\left(\|\vae\|_{L^2(I;\dot{B}^{\f d p}_{p,1})}\|\vthe\|^{l;\f{\beta_0}\ep}_{L^2(I;\dot{B}^{\f d p}_{p,1})}+\ep\|\vae\|_{L^\infty(I;\dot{B}^{\f d p}_{p,1})}\f1\ep\|\vthe\|^{h;\f{\beta_0}\ep}_{L^1(I;\dot{B}^{\f d p}_{p,1})}\right)\\
 &\leq C_3\|(\sep-\Theta)(t_0)\|^{l;\alpha}_{\dot{B}_{p,1}^{\f d p-1}}+C_3(\alpha\ep)\Big(M^{\ep,\alpha}_{p,q}[\vae,\vue,\vartheta^\ep;\tau^\ep,\sigma^\ep](I)\Big)^2+C_3\|\vv\|_{L^2(I;\dot{B}^{\f d p}_{p,1})}
 \|\Theta\|^{h;\alpha}_{L^2(I;\dot{B}^{\f d p}_{p,1})}\\
  &\quad+C_3
 M^{\ep,\alpha}_{p,q}[\vae,\vue-\vv,\vartheta^\ep;\tau^\ep,\sigma^\ep-\Theta](I)
 \left(M^{\ep,\alpha}_{p,q}[\vae,\vue,\vartheta^\ep;\tau^\ep,\sigma^\ep](I)+\|\vv\|_{L^2\left(I;\dot{B}_{p,1}^{\f d p}\right)}\right.\\
 &\quad+\left.\Big(M^{\ep,\alpha}_{p,q}[\vae,\vue,\vartheta^\ep;\tau^\ep,\sigma^\ep](I)\Big)^2\right)+C_3(\ep\alpha)^{\f1 {q(q-1)}}\|(\vae,\vue,\vthe)(t_0)\|_{X^\ep_2}\\
 &\quad+C_3 (\ep\alpha)^{\f1 {q(q-1)}}E^\ep[\vae,\vue,\vthe](I)\cdot \left(M^{\ep,\alpha}_{p,q}[\vae,\vue,\vthe](I)+\bigg(M^{\ep,\alpha}_{p,q}[\vae,\vue,\vartheta^\ep;\tau^\ep,\sigma^\ep](I)\bigg)^2\right),
\end{align*}

\vspace{5mm}

We find that
\begin{align}\label{equ9}
\p_t(\mathbb{P}\vue-\vv)-\mu\De(\mathbb{P}\vue-\vv)=H,
\end{align}

where
\begin{align}\label{equ10}
  H:&=\mathbb{P}\left(I(\ep\vae)\left(\f{\Grad(\vae+\vthe)}\ep+\Grad(\vae\vthe)-\mathcal{L}\vue\right)\right)\nonumber\\
  &\quad-\mathbb{P}\Big(\mathbb{Q}\vue\cdot\Grad\vue+(\mathbb{P}\vue-\vv)\cdot\Grad\mathbb{P}\vue
+\vv\cdot\Grad(\mathbb{P}\vue-\vv)
+\mathbb{P}\vue\cdot\Grad\mathbb{Q}\vue)\Big)
\end{align}

From Lemma \ref{Be3} and Calculating similarly to estimates for $G$, we have that
\begin{align*}
 &\|\mathbb{P}\vue-\vv\|_{L^\infty\left(I;\dot{B}_{p,1}^{\f d p-1}\right)\cap L^\infty\left(I;\dot{B}_{p,1}^{\f d p+1}\right)}\\
 &\leq C_3\|(\mathbb{P}\vue-\vv)(t_0)\|_{\dot{B}_{p,1}^{\f d p-1}}+C_3\|H\|_{L^1\left(I;\dot{B}_{p,1}^{\f d p-1}\right)}\\
 &\leq C_3\|(\mathbb{P}\vue-\vv)(t_0)\|_{\dot{B}_{p,1}^{\f d p-1}}
 +C_3(\alpha\ep)\Big(M^{\ep,\alpha}_{p,q}[\vae,\vue,\vartheta^\ep;\tau^\ep,\sigma^\ep](I)\Big)^2\\
 &+C_3
 M^{\ep,\alpha}_{p,q}[\vae,\vue-\vv,\vartheta^\ep;\tau^\ep,\sigma^\ep-\Theta](I)
 \left(M^{\ep,\alpha}_{p,q}[\vae,\vue,\vartheta^\ep;\tau^\ep,\sigma^\ep](I)+\|\vv\|_{L^2\left(I;\dot{B}_{p,1}^{\f d p}\right)}\right.\\
 &\quad+\left.\Big(M^{\ep,\alpha}_{p,q}[\vae,\vue,\vartheta^\ep;\tau^\ep,\sigma^\ep](I)\Big)^2\right).\ \ \Box
\end{align*}

\begin{Lemma}\label{Le5.3}  Let $2<p<\min{\{d,\f{2d}{d-2}\}}$.
Then there exist a constant $C_4=C_4(d,p,\mu,\kappa,\nu)$ such that for $\ep\leq1$,
\begin{align}\label{est58}
 &\|\mathbb{P}\vue-\vv\|_{L^\infty\Big(I;\dot{B}^{\f d p-1-(\f1 2-\f1 p)}_{p,1}\Big)\cap L^1\Big(I;\dot{B}^{\f d p+1-(\f1 2-\f1 p)}_{p,1}\Big)}\nonumber\\
 &+\|\sigma^\ep-\Theta\|^{l;\f{\beta_0}\ep}_{L^\infty\Big(I;\dot{B}^{\f d p-1-(\f1 2-\f1 p)}_{p,1}\Big)\cap L^2\Big(I;\dot{B}^{\f d p-(\f1 2-\f1 p)}_{p,1}\Big)}+\|(\tau^\ep,\mathbb{Q}\vue)\|^{l;\f{\beta_0}\ep}_{L^2\Big(I;\dot{B}^{\f d p-(\f1 2-\f1 p)}_{p,1}\Big)}\nonumber\\
& \leq C_4\left(\|(\mathbb{P}\vue-\vv)(t_0)\|_{\dot{B}^{\f d p-1-(\f1 2-\f1 p)}_{p,1}}+\|(\sigma^\ep-\Theta)(t_0)\|^{l;\f{\beta_0}\ep}_{\dot{B}^{\f d p-1-(\f1 2-\f1 p)}_{p,1}}+\ep^{\f1 2-\f 1 p}\|(\vae,\vue,\vthe)(t_0)\|_{X^\ep_2}\right)\nonumber\\
&+ C_4 \left(\|\mathbb{P}\vue-\vv\|_{ L^2\Big(I;\dot{B}^{\f d p-(\f1 2-\f1 p)}_{p,1}\Big)}+\|\sigma^\ep-\Theta\|^{l;\f{\beta_0}\ep}_{L^2\Big(I;\dot{B}^{\f d p-(\f1 2-\f1 p)}_{p,1}\Big)}+\|(\tau^\ep,\mathbb{Q}\vue)\|^{l;\f{\beta_0}\ep}_{L^2\Big(I;\dot{B}^{\f d p-(\f1 2-\f1 p)}_{p,1}\Big)}\right)\nonumber\\
&\cdot\left(M^{\ep,\alpha}_{p,q}[\vae,\vue,\vartheta^\ep;\tau^\ep,\sigma^\ep](I)+\Big(M^{\ep,\alpha}_{p,q}[\vae,\vue,\vartheta^\ep;\tau^\ep,\sigma^\ep](I)\Big)^2
+\|\vv\|_{L^2\Big(I;\dot{B}^{\f d p}_{p,1}\Big)}\right)\nonumber\\
&+C_4\ep^{\f1 2-\f1 p}\Big(M^{\ep,\alpha}_{p,q}[\vae,\vue,\vartheta^\ep;\tau^\ep,\sigma^\ep](I)+\|(\Theta,\vv)\|_{L^2\Big(I;\dot{B}^{\f d p}_{p,1}\Big)}\Big)\nonumber\\
&\cdot\left(\Big(M^{\ep,\alpha}_{p,q}[\vae,\vue,\vartheta^\ep;\tau^\ep,\sigma^\ep](I)\Big)^2
+M^{\ep,\alpha}_{p,q}[\vae,\vue,\vartheta^\ep;\tau^\ep,\sigma^\ep](I)+\|(\Theta,\vv)\|_{L^2\Big(I;\dot{B}^{\f d p}_{p,1}\Big)}\right)\nonumber\\
&+C_4 \ep^{\f1 2-\f1 p}E^\ep[\vae,\vue,\vthe](I)\cdot \left(M^{\ep,\alpha}_{p,q}[\vae,\vue,\vthe](I)+\bigg(M^{\ep,\alpha}_{p,q}[\vae,\vue,\vartheta^\ep;\tau^\ep,\sigma^\ep](I)\bigg)^2\right).
\end{align}
\end{Lemma}
\bProof From Lemma \ref{Be3}, \eqref{equ9} and \eqref{equ10},
\begin{align}\label{est55}
 &\|\mathbb{P}\vue-\vv\|_{L^\infty\Big(I;\dot{B}^{\f d p-1-(\f1 2-\f1 p)}_{p,1}\Big)\cap L^1\Big(I;\dot{B}^{\f d p+1-(\f1 2-\f1 p)}_{p,1}\Big)}\\
 &\leq C_4\|(\mathbb{P}\vue-\vv)(t_0)\|_{\dot{B}^{\f d p-1-(\f1 2-\f1 p)}_{p,1}}\nonumber\\
&+ C_4 \|(\vue,\vv)\|_{L^2\Big(I;\dot{B}^{\f d p}_{p,1}\Big)}\left(\|\mathbb{P}\vue-\vv\|_{ L^2\Big(I;\dot{B}^{\f d p-(\f1 2-\f1 p)}_{p,1}\Big)}+\|\mathbb{Q}\vue\|^{l;\f{\beta_0}\ep}_{L^2\Big(I;\dot{B}^{\f d p-(\f1 2-\f1 p)}_{p,1}\Big)}+\|\mathbb{Q}\vue\|^{h;\f{\beta_0}\ep}_{L^2\Big(I;\dot{B}^{\f d p-(\f1 2-\f1 p)}_{p,1}\Big)}\right)\nonumber\\
&+C_4\|\vae\|_{L^2\Big(I;\dot{B}^{\f d p}_{p,1}\Big)}\left(\|\tau^\ep\|^{l;\f{\beta_0}\ep}_{L^2\Big(I;\dot{B}^{\f d p-(\f1 2-\f1 p)}_{p,1}\Big)}+\|\vae\|^{h;\f{\beta_0}\ep}_{L^2\Big(I;\dot{B}^{\f d p-(\f1 2-\f1 p)}_{p,1}\Big)}+\|\mathbb{Q}\vue\|^{l;\f{\beta_0}\ep}_{L^1\Big(I;\dot{B}^{\f d p-(\f1 2-\f1 p)}_{p,1}\Big)}\right)\nonumber\\
&+C_4\ep\|\vae\|_{L^\infty\Big(I;\dot{B}^{\f d p-(\f1 2-\f1 p)}_{p,1}\Big)}
\cdot\left(\f1 \ep\|\vthe\|^{h;\f{\beta_0}\ep}_{L^1\Big(I;\dot{B}^{\f d p}_{p,1}\Big)}+\f1 \ep\|\vthe\|^{h;\f{\beta_0}\ep}_{L^1\Big(I;\dot{B}^{\f d p}_{p,1}\Big)} \ep\|\vae\|_{L^\infty\Big(I;\dot{B}^{\f d p}_{p,1}\Big)}\right.\nonumber\\
&\left.+\|\vthe\|^{l;\f{\beta_0}\ep}_{L^2\Big(I;\dot{B}^{\f d p}_{p,1}\Big)}\|\vae\|_{L^2\Big(I;\dot{B}^{\f d p}_{p,1}\Big)}+\|\mathbb{P}\vue\|_{L^1\Big(I;\dot{B}^{\f d p+1}_{p,1}\Big)}+\|\mathbb{Q}\vue\|^{h;\f{\beta_0}\ep}_{L^1\Big(I;\dot{B}^{\f d p+1}_{p,1}\Big)}\right).\nonumber
\end{align}

From Lemma \ref{Be3}, \eqref{equ11}, \eqref{equ12}, we have
\begin{align}\label{est56}
&\|\sigma^\ep-\Theta\|^{l;\f{\beta_0}\ep}_{L^\infty\Big(I;\dot{B}^{\f d p-1-(\f1 2-\f1 p)}_{p,1}\Big)\cap L^2\Big(I;\dot{B}^{\f d p-(\f1 2-\f1 p)}_{p,1}\Big)}\nonumber\\
&\leq C_4\|(\sigma^\ep-\Theta)(t_0)\|^{l;\f{\beta_0}\ep}_{\dot{B}^{\f d p-1-(\f1 2-\f1 p)}_{p,1}}+C_4\|\tau^\ep\|^{l;\f{\beta_0}\ep}_{L^2\Big(I;\dot{B}^{\f d p-(\f1 2-\f1 p)}_{p,1}\Big)}\nonumber\\
&+C_4\left(\|\mathbb{P}\vue-\vv\|_{L^2\Big(I;\dot{B}^{\f d p-(\f1 2-\f1 p)}_{p,1}\Big)}+\|\mathbb{Q}\vue\|^{l;\f{\beta_0}\ep}_{L^2\Big(I;\dot{B}^{\f d p-(\f1 2-\f1 p)}_{p,1}\Big)}+\|\mathbb{Q}\vue\|^{h;\f{\beta_0}\ep}_{L^2\Big(I;\dot{B}^{\f d p-(\f1 2-\f1 p)}_{p,1}\Big)}\right)\nonumber\\
&\cdot\left(\|\sep\|^{l;\f{\beta_0}\ep}_{L^2\Big(I;\dot{B}^{\f d p}_{p,1}\Big)}+\|\vae\|^{h;\f{\beta_0}\ep}_{L^2\Big(I;\dot{B}^{\f d p}_{p,1}\Big)}\right)\nonumber\\
&+C_4\ep^{\f1 2-\f1 p}\|(\mathbb{P}\vue-\vv,\mathbb{Q}\vue)\|_{L^2\Big(I;\dot{B}^{\f d p}_{p,1}\Big)}\f1 \ep\|\vthe\|^{h;\f{\beta_0}\ep}_{L^2\Big(I;\dot{B}^{\f d p-1}_{p,1}\Big)}+C_4\|\sigma^\ep-\Theta\|^{l;\f{\beta_0}\ep}_{L^2\Big(I;\dot{B}^{\f d p-(\f1 2-\f1 p)}_{p,1}\Big)}\|\vv\|_{L^2\Big(I;\dot{B}^{\f d p}_{p,1}\Big)}\nonumber\\
&+C_4\ep^{\f1 2-\f1 p}\left(\|(\vae,\Theta)\|^{h;\f{\beta_0}\ep}_{L^2\Big(I;\dot{B}^{\f d p}_{p,1}\Big)}+\f1\ep\|\vthe\|^{h;\f{\beta_0}\ep}_{L^2\Big(I;\dot{B}^{\f d p-1}_{p,1}\Big)}\right)\|\vv\|_{L^2\Big(I;\dot{B}^{\f d p}_{p,1}\Big)}\nonumber\\
&+C_4\ep\|\vae\|_{L^\infty\Big(I;\dot{B}^{\f d p-(\f1 2-\f1 p)}_{p,1}\Big)}\|\vthe\|^{l;\f{\beta_0}\ep}_{L^1\Big(I;\dot{B}^{\f d p+1}_{p,1}\Big)}+C_4\ep^{\f1 2-\f1 p}\ep\|\vae\|_{L^\infty\Big(I;\dot{B}^{\f d p}_{p,1}\Big)}\f1 \ep\|\vthe\|^{h;\f{\beta_0}\ep}_{L^1\Big(I;\dot{B}^{\f d p}_{p,1}\Big)}\nonumber\\
&+C_4\ep^{\f1 2-\f1 p}\left(1+\ep\|\vae\|_{L^\infty\Big(I;\dot{B}^{\f d p}_{p,1}\Big)}\right)\|\vue\|^2_{L^2\Big(I;\dot{B}^{\f d p}_{p,1}\Big)}
\end{align}
Note that
\begin{align}\label{est60}
&\ep\|\vae\|_{L^\infty\Big(I;\dot{B}^{\f d p-(\f1 2-\f1 p)}_{p,1}\Big)}\leq C_4\ep^{\f1 2-\f1 p}\|\vae\|^{l;\f{\beta_0}\ep}_{L^\infty\Big(I;\dot{B}^{\f d p-1}_{p,1}\Big)}+C_4\ep^{\f1 2-\f1 p}\|\vae\|^{h;\f{\beta_0}\ep}_{L^\infty\Big(I;\dot{B}^{\f d p}_{p,1}\Big)}\nonumber\\
&\leq C_4 \ep^{\f1 2-\f1 p}E^\ep[\vae,\vue,\vthe](I)
\end{align}
From Lemma \ref{Be9} and \eqref{est33}, we have that
\begin{align}\label{est59}
 &\|\tau^\ep\|^{l;\f{\beta_0}\ep}_{L^2\Big(I;\dot{B}^{\f d p-(\f1 2-\f1 p)}_{p,1}\Big)}\leq C_4\ep^{\f1 2-\f 1 p}\|(\vae,\vue,\vthe)(t_0)\|_{X^\ep_2}\\
 &+C_4 \ep^{\f1 2-\f1 p}E^\ep[\vae,\vue,\vthe](I)\cdot \left(M^{\ep,\alpha}_{p,q}[\vae,\vue,\vthe](I)+\bigg(M^{\ep,\alpha}_{p,q}[\vae,\vue,\vartheta^\ep;\tau^\ep,\sigma^\ep](I)\bigg)^2\right)\nonumber
\end{align}
Adding \eqref{est55} to \eqref{est56}, and then using \eqref{est60}, \eqref{est59}, Proposition \ref{pro1} yields to \eqref{est58}.\ \  $\Box$
\section{Proof of Theorem \ref{main}}\label{s6}
From \cite{CD1,D1}, if the initial data $(a_0,\vu_0,\vartheta_0)$ satisfy \eqref{cond1} and $\inf_{x\in\mathbb{R}^d} (1+\ep a_0(x))>0$, then there exist a positive time $0<T<\infty$ such that systems \eqref{equ2} exist a unique regular solution $(\vae,\vue,\vthe)$ on $[0,T]$, which obeys
\begin{align}\label{cond5}
 &\vae\in C\Big([0,T];\dot{B}_{2,1}^{\f{d} 2}\cap\dot{B}_{2,1}^{\f d 2-1}\Big),\ \ \vue\in C\Big([0,T];\dot{B}_{2,1}^{\f{d} 2-1}\Big)\cap L^1\Big(0,T;\dot{B}_{2,1}^{\f{d} 2+1}\Big),\nonumber\\
 &\vthe\in C\Big([0,T];\dot{B}_{2,1}^{\f{d} 2-2}\Big)\cap L^1\Big(0,T;\dot{B}_{2,1}^{\f{d} 2}\Big),\ \ \inf_{(t,x)\in[0,T]\times\mathbb{R}^d} (1+\ep a(t,x))>0.
\end{align}

Let $0<T_{1}<\infty$ be a positive time  and $(\vae,\vue,\vthe)$ be a solution to systems \eqref{equ2} on $[0,T_1)$ which satisfy \eqref{cond5} for any $0<T<T_1$, then the solution can be continued beyond $T_1$ provided that
\begin{align}\label{cond6}
\|\vae\|_{L^\infty\left(0,T_1;\dot{B}_{2,1}^{\f d 2}\right)}&+\|\vue\|_{L^1\left(0,T_1;\dot{B}_{2,1}^{\f d 2+1}\right)}+\|\vthe\|_{L^1\left(0,T_1;\dot{B}_{2,1}^{\f d 2}\right)}<\infty,\nonumber\\
 &\inf_{(t,x)\in[0,T_1)\times\mathbb{R}^d} (1+\ep a(t,x))>0.
\end{align}

We define
\[
 T^{\ast} :=\sup\{0<T<\infty:(\vae,\vue,\vthe) \text{ satisfies \eqref{cond5} for } T \}.
\]

Using the above continuation criterion\footnote{The proof of continuation criterion is similar to Proposition 10.10 in \cite{BCD}.}, we will show $T^{\ast}=\infty$.

\vspace{5mm}

Let $(\Theta,\vv)$ be a regular solution to systems \eqref{equ5} supplement with initial data $\Theta_0=\f{\vartheta_0-a_0}{\sqrt{2}}$ and $\vv_0=\mathbb{P}\vu_0$, and be in the class \eqref{assum1}. Note that $\|f\|_{L^{\infty}}\leq C_5(d,p) \|f\|_{\dot{B}_{p,1}^{\f d p}}$, then we set $C_0=\max\{1,C_1,C_2,C_3,C_4,C_5\}$ where $C_1,C_2,C_3,C_4$ appears in the previous lemma and
 Let
\[
\de_0:=\min\left\{\f1 {96C^2_0},\ \ \f1 2\|(\vv,\Theta)\|_{L^q\left(0,\infty;\dot{B}^{\f d p-1+\f2 q}_{p,1}\right)\cap L^1\left(0,\infty;\dot{B}^{\f d p+1}_{p,1}\right)}\right\}.
\]
From Lemma 2.10 in \cite{F}, there exist a positive number $N_0=N_0(\delta_0,p,q,\vv,\Theta)$ and a time sequence $\{T_n=T_n(\delta_0,p,q,\vv,\Theta)\}^{N_0}_{n=0}$ such that
$0=T_0<T_1<\cdots<T_{N_0-1}<T_{N_0}=\infty$ and
\begin{align}\label{finally1}
 \|(\vv,\Theta)\|_{L^q\left(T_{n-1},T_{n};\dot{B}^{\f d p-1+\f2 q}_{p,1}\right)\cap L^1\left(T_{n-1},T_{n};\dot{B}^{\f d p+1}_{p,1}\right)}\leq\de_0,\ \ n=1,2,\cdots,N_0.
\end{align}
Note that if $\ep\leq1$,
\[
\|(a_0,\mathbb{Q}\vu_0,\vartheta_0)\|^{h;\alpha}_{X^\ep_p}\leq C\left(\|a_0\|^{h;\f{\beta_0}\ep}_{\dot{B}^{\f d 2}_{2,1}}+\|(\vartheta_0,\mathbb{Q}\vu_0)\|^{h;\f{\beta_0}\ep}_{\dot{B}^{\f d 2-1}_{2,1}}+\|(a_0,\vartheta_0,\mathbb{Q}\vu_0)\|^{m;\alpha,\f{\beta_0}\ep}_{\dot{B}^{\f d 2-1}_{2,1}}\right),
\]
where the constant $C$ depend on $d,p$.

There exist a positive constant $\alpha_0=\alpha_0(d,p,q,C_0,a_0,\vu_0,\vartheta_0,N_0,\delta_0,\vv,\Theta)$ such that
\begin{align}\label{cond8}
\|(\Theta,\vv)\|^{h;\f{\alpha_0} {64}}_{L^2(I;\dot{B}_{p,1}^{\f d p})\cap L^{q'}(I;\dot{B}_{p,1}^{\f d p-1+\f2 {q'}})}\leq\f{\de_0}{(4C_0)^{N_0}}.
\end{align}

and for $\ep\leq1$,
\begin{align}\label{cond9}
 \|(a_0,\mathbb{Q}\vu_0,\vartheta_0)\|^{h;\alpha_0}_{X^\ep_p}\leq\f{\de_0}{(4C_0)^{N_0}}.
\end{align}

Then for such $\alpha_0$, we select a positive constant $\ep_0=\ep_0(\alpha_0,\delta_0,N_0,q,C_0,a_0,\vu_0,\vartheta_0)\leq\min\{1,\ \ \f{1}{\alpha_0},\ \ \f{\beta_0}{16\alpha_0}\}$ such that
\begin{align}\label{cond10}
 (\alpha_0\ep)^{\f{1}{q(q-1)}}\Big(\|(a_0,\vu_0,\vartheta_0)\|_{X^\ep_2}+\delta_0\Big)\leq \f1 {36}\f{\de_0}{(4C_0)^{N_0}}
\end{align}

for any $0<\ep\leq\ep_0$.

Now we fix a Much number $0<\ep\leq\ep_0$. We will use reductio and absurdum to show $T^\ast\geq T_1$. We define a time
\[
T^\ast_1:=\sup\{0<T<T^\ast, \eqref{finally2}-\eqref{finally8} \text{ hold at } T\}
\]
where
\begin{align}\label{finally2}
E^\ep[\vae,\vue,\vthe](0,T)\leq 4C_0\|(a_0,\vu_0,\vartheta_0)\|_{X^\ep_2},
\end{align}
\begin{align}\label{finally3}
M^{\ep,\alpha}_{p,q}[\vae,\vue-\vv,\vartheta^\ep;\tau^\ep,\sigma^\ep-\Theta](0,T)\leq \f{\de_0}{(4C_0)^{N_0-1}},
\end{align}
\begin{align}\label{finally4}
\|\mathbb{P}\vue-\vv\|_{L^\infty\left(0,T;\dot{B}^{\f d p-1}_{p,1}\right)\cap L^1\left(0,T;\dot{B}^{\f d p+1}_{p,1}\right)}\leq \f1 4 \f{\de_0}{(4C_0)^{N_0-1}}.
\end{align}
\begin{align}\label{finally8}
  &\|\mathbb{P}\vue-\vv\|_{L^\infty\Big(I;\dot{B}^{\f d p-1-(\f1 2-\f1 p)}_{p,1}\Big)\cap L^1\Big(I;\dot{B}^{\f d p+1-(\f1 2-\f1 p)}_{p,1}\Big)}+\|(\tau^\ep,\mathbb{Q}\vue)\|^{l;\f{\beta_0}\ep}_{L^2\Big(I;\dot{B}^{\f d p-(\f1 2-\f1 p)}_{p,1}\Big)}\nonumber\\
 &+\|\sigma^\ep-\Theta\|^{l;\f{\beta_0}\ep}_{L^\infty\Big(I;\dot{B}^{\f d p-1-(\f1 2-\f1 p)}_{p,1}\Big)\cap L^2\Big(I;\dot{B}^{\f d p-(\f1 2-\f1 p)}_{p,1}\Big)}\leq 2\ep^{\f1 2-\f1 p}4C_0\Big(\|(a_0,\vu_0,\vartheta_0)\|_{X^\ep_2}+1\Big)
\end{align}
It is clear that $T^\ast_1>0$. We firstly assume $T^\ast< T_1$. Let $0<T<T^\ast_1$. From Proposition \ref{pro1}, \eqref{finally1}, \eqref{finally1}, we have that
\begin{align*}
 &M^{\ep,\alpha}_{p,q}[\vae,\vue,\vartheta^\ep;\tau^\ep,\sigma^\ep](0,T)\\\
 &\leq M^{\ep,\alpha}_{p,q}[\vae,\vue-\vv,\vartheta^\ep;\tau^\ep,\sigma^\ep-\Theta](0,T)
 +2\|(\vv,\Theta)\|_{L^q\left(0,T;\dot{B}^{\f d p-1+\f2 q}_{p,1}\right)\cap L^1\left(0,T;\dot{B}^{\f d p+1}_{p,1}\right)}\leq3\delta_0,
\end{align*}
\begin{align*}
 \ep\|\vae\|_{L^\infty(0,T;L^\infty)}&\leq C_0\ep\|\vae\|_{L^\infty(0,T;\dot{B}^{\f d p}_{p,1})}\\
 &\leq C^2_0 M^{\ep,\alpha}_{p,q}[\vae,\vue,\vartheta^\ep;\tau^\ep,\sigma^\ep](0,T)\leq 3C^2_0\delta_0\leq \f1 2,
\end{align*}
which implies
\begin{align}\label{est40}
  \f1 2 \leq1+\ep a(t,x)\leq \f3 2, \ \ \forall (t,x)\in [0,T]\times \mathbb{R}^d.
\end{align}

From Lemma \ref{Le5.1},
\begin{align}\label{est32}
 E^\ep[\vae,\vue,\vthe](0,T)&\leq C_0\|(a_0,\vu_0,\vartheta_0)\|_{X^\ep_2}+4C_0\|(a_0,\vu_0,\vartheta_0)\|_{X^\ep_2}(3\delta_0+9\delta^2_0)\nonumber\\
 &\leq 2C_0\|(a_0,\vu_0,\vartheta_0)\|_{X^\ep_2}
\end{align}

From Lemma \ref{Be9}, Lemma \ref{Le5.2}, \eqref{cond8}, \eqref{cond9} and \eqref{cond10}
\begin{align}\label{est34}
&\|(\vae,\mathbb{Q}\vue,\vthe)\|^{h;\alpha}_{N^{\ep,\alpha}_{p,q}(0,T)}\leq C_0\f{\de_0}{(4C_0)^{N_0}}+C_0(\alpha\ep)9\delta^2_0+C_03\delta_0\Big(\f{\de_0}{(4C_0)^{N_0}}+\f{\de_0}{(4C_0)^{N_0-1}}\Big),
\end{align}
\begin{align}\label{est35}
&\|(\tau^\ep,\mathbb{Q}\vue)\|^{l;\alpha}_{L^q\left(0,T;\dot{B}_{p,1}^{\f d p-1+\f2 q}\right)}+\|\sigma^\ep-\Theta\|^{l;\alpha}_{L^q(0,T;\dot{B}^{\f d p-1+\f2 q}_{p,1})\cap L^{q'}(0,T;\dot{B}^{\f d p-1+\f2 {q'}}_{p,1})}\nonumber\\
&\leq C_0(\alpha\ep)^{\f1 q}\|(a_0,\vu_0,\vartheta_0)\|_{X^\ep_2}+C_0(\alpha\ep)^{\f1 q}4C_0\|(a_0,\vu_0,\vartheta_0)\|_{X^\ep_2}(3\delta_0+9\delta^2_0)+C_0(\alpha\ep)^{\f1 {q(q-1)}}\|(a_0,\vu_0,\vartheta_0)\|_{X^\ep_2}\nonumber\\
&\quad+C_0(\alpha\ep)^{\f1 {q(q-1)}}4C_0\|(a_0,\vu_0,\vartheta_0)\|_{X^\ep_2}(3\delta_0+9\delta^2_0)+C_0(\alpha\ep)9\delta^2_0+C_0\delta_0\f{\de_0}{(4C_0)^{N_0}}\nonumber\\
&\quad+C_0\f{\de_0}{(4C_0)^{N_0-1}}(4\delta_0+9\delta^2_0)
\end{align}
\begin{align}\label{est36}
 &(\alpha\ep) E^\ep[\vae,\mathbb{Q}\vue,\vthe](0,T)
  + \left(\|(\tau^\ep,\mathbb{Q}\vue)\|^{l;\alpha}_{{L^q(0,T;\dot{B}^{\f d p-1+\f2 q}_{p,1})}}\right)^\f1 {q-1}\Big(  E^\ep[\vae,\mathbb{Q}\vue,\vthe](0,T)\Big)^\f{q-2}{q-1}\nonumber\\
  &\leq (\alpha\ep)4C_0\|(a_0,\vu_0,\vartheta_0)\|_{X^\ep_2}+C_0(\alpha\ep)^{\f1 {q(q-1)}}\|(a_0,\vu_0,\vartheta_0)\|_{X^\ep_2}\\
  &\quad+C_0(\alpha\ep)^{\f1 {q(q-1)}}4C_0\|(a_0,\vu_0,\vartheta_0)\|_{X^\ep_2}(3\delta_0+9\delta^2_0)\nonumber,
\end{align}
\begin{align}\label{est38}
\|\mathbb{P}\vue-\vv\|_{L^\infty\left(0,T;\dot{B}^{\f d p-1}_{p,1}\right)\cap L^1\left(0,T;\dot{B}^{\f d p+1}_{p,1}\right)}&\leq C_0(\alpha\ep)9\delta^2_0+C_0\f{\de_0}{(4C_0)^{N_0-1}}(4\delta_0+9\delta^2_0)\nonumber\\
&\leq \f5 {96}\f{\de_0}{(4C_0)^{N_0-1}}
\end{align}

Collecting \eqref{est34}-\eqref{est36} and using \eqref{est38}, we get that
\begin{align}\label{est37}
 &M^{\ep,\alpha}_{p,q}[\vae,\vue-\vv,\vartheta^\ep;\tau^\ep,\sigma^\ep-\Theta](0,T)\nonumber\\
 &\leq \f{\de_0}{(4C_0)^{N_0-1}}\left(\f1 4+\delta_0+9C_0\delta^2_0+7C_0\delta_0\right)+C_0(\alpha\ep)^{\f1 {q(q-1)}}18\delta^2_0\nonumber\\
 &\quad+(\alpha\ep)^{\f1 {q(q-1)}}4C_0\|(a_0,\vu_0,\vartheta_0)\|_{X^\ep_2}\left(\f7 4+9C_0\delta_0+27C_0\delta^2_0\right)+2\cdot\f5 {96}\f{\de_0}{(4C_0)^{N_0-1}}\nonumber\\
 &\leq \f{\de_0}{(4C_0)^{N_0-1}}\left(\f1 4+\f1 {96}+\f1 {96}+\f7 {96}\right)+\f{\de_0}{(4C_0)^{N_0-1}}\f1 {96}+\f{\de_0}{(4C_0)^{N_0-1}}\f1 {18}+\f{10} {96}\f{\de_0}{(4C_0)^{N_0-1}}\nonumber\\
  &\leq \f3 4\f{\de_0}{(4C_0)^{N_0-1}}.
\end{align}
From Lemma \ref{Le5.3},
\begin{align}\label{est61}
 &\|\mathbb{P}\vue-\vv\|_{L^\infty\Big(I;\dot{B}^{\f d p-1-(\f1 2-\f1 p)}_{p,1}\Big)\cap L^1\Big(I;\dot{B}^{\f d p+1-(\f1 2-\f1 p)}_{p,1}\Big)}+\|(\tau^\ep,\mathbb{Q}\vue)\|^{l;\f{\beta_0}\ep}_{L^2\Big(I;\dot{B}^{\f d p-(\f1 2-\f1 p)}_{p,1}\Big)}\nonumber\\
 &+\|\sigma^\ep-\Theta\|^{l;\f{\beta_0}\ep}_{L^\infty\Big(I;\dot{B}^{\f d p-1-(\f1 2-\f1 p)}_{p,1}\Big)\cap L^2\Big(I;\dot{B}^{\f d p-(\f1 2-\f1 p)}_{p,1}\Big)}\nonumber\\
 &\leq \ep^{\f1 2-\f1 p}C_0\|(a_0,\vu_0,\vartheta_0)\|_{X^\ep_2}+2C_0\ep^{\f1 2-\f1 p}4C_0(\|(a_0,\vu_0,\vartheta_0)\|_{X^\ep_2}+1)(4\delta_0+9\delta^2_0)+C_0\ep^{\f1 2-\f1 p}4\delta_0(4\delta_0+9\delta^2_0)\nonumber\\
 &+C_0\ep^{\f1 2-\f1 p}4C_0\|(a_0,\vu_0,\vartheta_0)\|_{X^\ep_2}(3\delta_0+9\delta^2_0)\nonumber\\
 &\leq \ep^{\f1 2-\f1 p}4C_0\Big(\|(a_0,\vu_0,\vartheta_0)\|_{X^\ep_2}+1\Big)
\end{align}

It is clear that
\begin{align}\label{est39}
 \|\vae\|_{L^\infty\left(0,T;\dot{B}_{2,1}^{\f d 2}\right)}&+\|\vue\|_{L^1\left(0,T;\dot{B}_{2,1}^{\f d 2+1}\right)}+\|\vthe\|_{L^1\left(0,T;\dot{B}_{2,1}^{\f d 2}\right)}\leq C(T,\ep) E^\ep[\vae,\vue,\vthe](0,T).
\end{align}

Then from \eqref{est32}, \eqref{est38}, \eqref{est37}, \eqref{est61}, we have $T^\ast_1=T^\ast$. Also from \eqref{est40}, \eqref{est39} and continuation criterion, the solution can be continued beyond $T^\ast$ which contradicts the definition of $T^\ast$. So we infer that $T^\ast\geq T^\ast_1\geq T_1$ .

Next we  use reductio and absurdum to show $T^\ast\geq T_k,k=2,3,\dots,N_0$. We define the time
\[
T^\ast_k:=\sup\{T_{k-1}<T<T^\ast, \eqref{finally5}-\eqref{finally9} \text{ hold at } T\}
\]
where
\begin{align}\label{finally5}
E^\ep[\vae,\vue,\vthe](T_{k-1},T)\leq (4C_0)^k\|(a_0,\vu_0,\vartheta_0)\|_{X^\ep_2},
\end{align}
\begin{align}\label{finally6}
M^{\ep,\alpha}_{p,q}[\vae,\vue-\vv,\vartheta^\ep;\tau^\ep,\sigma^\ep-\Theta](T_{k-1},T)\leq \f{\de_0}{(4C_0)^{N_0-k}},
\end{align}
\begin{align}\label{finally7}
\|\mathbb{P}\vue-\vv\|_{L^\infty\left(T_{k-1},T;\dot{B}^{\f d p-1}_{p,1}\right)\cap L^1\left(T_{k-1},T;\dot{B}^{\f d p+1}_{p,1}\right)}\leq \f1 4\f{\de_0}{(4C_0)^{N_0-k}}.
\end{align}
\begin{align}\label{finally9}
  &\|\mathbb{P}\vue-\vv\|_{L^\infty\Big(T_{k-1},T;\dot{B}^{\f d p-1-(\f1 2-\f1 p)}_{p,1}\Big)\cap L^1\Big(T_{k-1},T;\dot{B}^{\f d p+1-(\f1 2-\f1 p)}_{p,1}\Big)}+\|(\tau^\ep,\mathbb{Q}\vue)\|^{l;\f{\beta_0}\ep}_{L^2\Big(T_{k-1},T;\dot{B}^{\f d p-(\f1 2-\f1 p)}_{p,1}\Big)}\nonumber\\
 &+\|\sigma^\ep-\Theta\|^{l;\f{\beta_0}\ep}_{L^\infty\Big(T_{k-1},T;\dot{B}^{\f d p-1-(\f1 2-\f1 p)}_{p,1}\Big)\cap L^2\Big(T_{k-1},T;\dot{B}^{\f d p-(\f1 2-\f1 p)}_{p,1}\Big)}\leq 2\ep^{\f1 2-\f1 p}(4C_0)^k\Big(\|(a_0,\vu_0,\vartheta_0)\|_{X^\ep_2}+1\Big)
\end{align}
Then we assume $T_{k-1}<T^\ast< T_k$. Let $T_{k-1}<T<T^\ast$.
We have that
\begin{align*}
 &M^{\ep,\alpha}_{p,q}[\vae,\vue,\vartheta^\ep;\tau^\ep,\sigma^\ep](T_{k-1},T)\\\
 &\leq M^{\ep,\alpha}_{p,q}[\vae,\vue-\vv,\vartheta^\ep;\tau^\ep,\sigma^\ep-\Theta](T_{k-1},T)
 +2\|(\vv,\Theta)\|_{L^q\left(T_{k-1},T;\dot{B}^{\f d p-1+\f2 q}_{p,1}\right)\cap L^1\left(T_{k-1},T;\dot{B}^{\f d p+1}_{p,1}\right)}\leq3\delta_0,
\end{align*}
\begin{align*}
 \ep\|\vae\|_{L^\infty(T_{k-1},T;L^\infty)}&\leq C_0\ep\|\vae\|_{L^\infty(T_{k-1},T;\dot{B}^{\f d p}_{p,1})}\\
& \leq C^2_0 M^{\ep,\alpha}_{p,q}[\vae,\vue,\vartheta^\ep;\tau^\ep,\sigma^\ep](T_{k-1},T)\leq 3C^2_0\delta_0\leq \f1 2,
\end{align*}
which implies
\begin{align}\label{est41}
  \f1 2 \leq1+\ep a(t,x)\leq \f3 2, \ \ \forall (t,x)\in [T_k-1,T]\times \mathbb{R}^d.
\end{align}
It is clear that
\begin{align*}
&\|(\vae,\vue,\vthe)(T_{k-1})\|_{X^\ep_2}\leq 2 E^\ep[\vae,\vue,\vthe](T_{k-2},T_{k-1}),\\ &\|(\vae,\mathbb{Q}\vue,\vthe)(T_{k-1})\|^{h;\alpha}_{X^\ep_p}\leq M^{\ep,\alpha}_{p,q}[\vae,\vue-\vv,\vartheta^\ep;\tau^\ep,\sigma^\ep-\Theta](T_{k-2},T_{k-1}),\\
&\|(\mathbb{P}\vue-\vv)(T_{k-1})\|_{\dot{B}_{p,1}^{\f d p-1}}\leq \|\mathbb{P}\vue-\vv\|_{L^\infty\left(T_{k-2},T_{k-1};\dot{B}^{\f d p-1}_{p,1}\right)},\\
&\|(\mathbb{P}\vue-\vv)(T_{k-1})\|_{\dot{B}_{p,1}^{\f d p-1-(\f1 2-\f1 p)}}+\|(\sigma^\ep-\Theta)(T_{k-1})\|^{l;\f{\beta_0}\ep}_{\dot{B}_{p,1}^{\f d p-1-(\f1 2-\f1 p)}}\\
&\leq\|\mathbb{P}\vue-\vv\|_{L^\infty\Big(T_{k-2},T_{k-1};\dot{B}^{\f d p-1-(\f1 2-\f1 p)}_{p,1}\Big)}
 +\|\sigma^\ep-\Theta\|^{l;\f{\beta_0}\ep}_{L^\infty\Big(T_{k-2},T_{k-1};\dot{B}^{\f d p-1-(\f1 2-\f1 p)}_{p,1}\Big)}.\nonumber\\
\end{align*}

From Lemma \ref{Le5.1},
\begin{align}\label{est42}
 E^\ep[\vae,\vue,\vthe](T_{k-1},T)&\leq C_0\|(\vae,\vue,\vthe)(T_{k-1})\|_{X^\ep_2}+(4C_0)^{k-1}\|(a_0,\vu_0,\vartheta_0)\|_{X^\ep_2}(3\delta_0+9\delta^2_0)\nonumber\\
 &\leq 2C_0(4C_0)^{k-1}\|(a_0,\vu_0,\vartheta_0)\|_{X^\ep_2}+(4C_0)^{k-1}\|(a_0,\vu_0,\vartheta_0)\|_{X^\ep_2}\nonumber\\
 &\leq \f3 4 (4C_0)^{k}\|(a_0,\vu_0,\vartheta_0)\|_{X^\ep_2},
\end{align}

From Lemma \ref{Le5.2},
\begin{align}\label{est43}
&\|(\vae,\mathbb{Q}\vue,\vthe)\|^{h;\alpha}_{N^{\ep,\alpha}_{p,q}(T_{k-1},T)}\leq C_0\f{\de_0}{(4C_0)^{N_0-(k-1)}}+C_0(\alpha\ep)9\delta^2_0+C_03\delta_0\Big(\f{\de_0}{(4C_0)^{N_0-k}}+\f{\de_0}{(4C_0)^{N_0}}\Big),
\end{align}

\begin{align}\label{est44}
&\|(\tau^\ep,\mathbb{Q}\vue)\|^{l;\alpha}_{L^q\left(T_{k-1},T;\dot{B}_{p,1}^{\f d p-1+\f2 q}\right)}+\|\sigma^\ep-\Theta\|^{l;\alpha}_{L^q(0,T;\dot{B}^{\f d p-1+\f2 q}_{p,1})\cap L^{q'}(t_{k-1},T;\dot{B}^{\f d p-1+\f2 {q'}}_{p,1})}\nonumber\\
&\leq C_0(\alpha\ep)^{\f1 q}2(4C_0)^{k-1}\|(a_0,\vu_0,\vartheta_0)\|_{X^\ep_2}+C_0(\alpha\ep)^{\f1 q}(4C_0)^k\|(a_0,\vu_0,\vartheta_0)\|_{X^\ep_2}(3\delta_0+9\delta^2_0)\nonumber\\
&\quad+C_0\f{\de_0}{(4C_0)^{N_0-(k-1)}}+C_0(\alpha\ep)^{\f1 {q(q-1)}}2(4C_0)^{k-1}\|(a_0,\vu_0,\vartheta_0)\|_{X^\ep_2}\nonumber\\
&\quad+C_0(\alpha\ep)^{\f1 {q(q-1)}}(4C_0)^k\|(a_0,\vu_0,\vartheta_0)\|_{X^\ep_2}(3\delta_0+9\delta^2_0)\nonumber\\
&\quad+C_0(\alpha\ep)9\delta^2_0+C_0\delta_0\f{\de_0}{(4C_0)^{N_0}}
+C_0\f{\de_0}{(4C_0)^{N_0-k}}(4\delta_0+9\delta^2_0)
\end{align}
\begin{align}\label{est45}
 &(\alpha\ep) E^\ep[\vae,\mathbb{Q}\vue,\vthe](T_{k-1},T)
  + \left(\|(\tau^\ep,\mathbb{Q}\vue)\|^{l;\alpha}_{{L^q(T_{k-1},T;\dot{B}^{\f d p-1+\f2 q}_{p,1})}}\right)^\f1 {q-1}\Big(  E^\ep[\vae,\mathbb{Q}\vue,\vthe](T_{k-1},T)\Big)^\f{q-2}{q-1}\nonumber\\
  &\leq (\alpha\ep)(4C_0)^{k}\|(a_0,\vu_0,\vartheta_0)\|_{X^\ep_2}+C_0(\alpha\ep)^{\f1 {q(q-1)}}2(4C_0)^{k-1}\|(a_0,\vu_0,\vartheta_0)\|_{X^\ep_2}\\
  &\quad+C_0(\alpha\ep)^{\f1 {q(q-1)}}(4C_0)^k\|(a_0,\vu_0,\vartheta_0)\|_{X^\ep_2}(3\delta_0+9\delta^2_0)\nonumber,
\end{align}

\begin{align}\label{est46}
&\|\mathbb{P}\vue-\vv\|_{L^\infty\left(T_{k-1},T;\dot{B}^{\f d p-1}_{p,1}\right)\cap L^1\left(0,T;\dot{B}^{\f d p+1}_{p,1}\right)}\nonumber\\
&\leq C_0\f1 4\f{\de_0}{(4C_0)^{N_0-(k-1)}}+C_0(\alpha\ep)9\delta^2_0+C_0\f{\de_0}{(4C_0)^{N_0-k}}(4\delta_0+9\delta^2_0)\nonumber\\
&\leq \f{11} {96}\f{\de_0}{(4C_0)^{N_0-k}}
\end{align}

Collecting \eqref{est43}-\eqref{est45} and using \eqref{est46}, we obtain that
\begin{align}\label{est47}
&M^{\ep,\alpha}_{p,q}[\vae,\vue-\vv,\vartheta^\ep;\tau^\ep,\sigma^\ep-\Theta](T_{k-1},T)\nonumber\\
&\leq \f{\de_0}{(4C_0)^{N_0-k}}\left(\f1 2+\delta_0+7C_0\delta_0+9C_0\delta^2_0\right)+C_0(\alpha\ep)^{\f1 {q(q-1)}}18\delta^2_0\nonumber\\
&\quad+(\alpha\ep)^{\f1 {q(q-1)}}(4C_0)^k\|(a_0,\vu_0,\vartheta_0)\|_{X^\ep_2}\left(\f5 2+9C_0\delta_0+27C_0\delta^2_0\right)+2\cdot\f{11} {96}\f{\de_0}{(4C_0)^{N_0-k}}\nonumber\\
&\leq \f{\de_0}{(4C_0)^{N_0-k}}\left(\f1 2+\f9{96}\right)+\f1{96}\f{\de_0}{(4C_0)^{N_0-k}}+\f1 {12}\f{\de_0}{(4C_0)^{N_0-k}}+\f{22} {96}\f{\de_0}{(4C_0)^{N_0-k}}\nonumber\\
&\leq\f{88} {96}\f{\de_0}{(4C_0)^{N_0-k}}.
\end{align}
\begin{align}\label{est62}
 &\|\mathbb{P}\vue-\vv\|_{L^\infty\Big(T_{k-1},T;\dot{B}^{\f d p-1-(\f1 2-\f1 p)}_{p,1}\Big)\cap L^1\Big(T_{k-1},T;\dot{B}^{\f d p+1-(\f1 2-\f1 p)}_{p,1}\Big)}+\|(\tau^\ep,\mathbb{Q}\vue)\|^{l;\f{\beta_0}\ep}_{L^2\Big(T_{k-1},T;\dot{B}^{\f d p-(\f1 2-\f1 p)}_{p,1}\Big)}\nonumber\\
 &+\|\sigma^\ep-\Theta\|^{l;\f{\beta_0}\ep}_{L^\infty\Big(T_{k-1},T;\dot{B}^{\f d p-1-(\f1 2-\f1 p)}_{p,1}\Big)\cap L^2\Big(T_{k-1},T;\dot{B}^{\f d p-(\f1 2-\f1 p)}_{p,1}\Big)}\nonumber\\
 &\leq C_0 2\ep^{\f1 2-\f1 p}(4C_0)^{k-1}\Big(\|(a_0,\vu_0,\vartheta_0)\|_{X^\ep_2}+1\Big)+C_0 \ep^{\f1 2-\f1 p}2(4C_0)^{k-1}\|(a_0,\vu_0,\vartheta_0)\|_{X^\ep_2}\nonumber\\
 &+C_02\ep^{\f1 2-\f1 p}(4C_0)^{k}\Big(\|(a_0,\vu_0,\vartheta_0)\|_{X^\ep_2}+1\Big)(4\delta_0+9\delta^2_0)+C_0\ep^{\f1 2-\f1 p}4\delta_0(4\delta_0+9\delta^2_0)\nonumber\\
 &+C_0\ep^{\f1 2-\f1 p}(4C_0)^{k}\|(a_0,\vu_0,\vartheta_0)\|_{X^\ep_2}(3\delta_0+9\delta^2_0)\nonumber\\
 &\leq \f3 4\cdot 2\ep^{\f1 2-\f1 p}(4C_0)^{k}\Big(\|(a_0,\vu_0,\vartheta_0)\|_{X^\ep_2}+1\Big)
 \end{align}
It is clear that
\begin{align}\label{est48}
 &\|\vae\|_{L^\infty\left(T_{k-1},T;\dot{B}_{2,1}^{\f d 2}\right)}+\|\vue\|_{L^1\left(T_{k-1},T;\dot{B}_{2,1}^{\f d 2+1}\right)}+\|\vthe\|_{L^1\left(T_{k-1},T;\dot{B}_{2,1}^{\f d 2}\right)}\nonumber\\
 &\leq C(T,\ep) E^\ep[\vae,\vue,\vthe](T_{k-1},T).
\end{align}

Then from \eqref{est42},\eqref{est46},\eqref{est47}, \eqref{est62}, we have that $T^\ast_k=T^\ast$.  From \eqref{est41} and \eqref{est48} and continuation criterion, the solution can be continued beyond $T^\ast$ which contradicts the definition of $T^\ast$. So we infer that $T^\ast\geq T^\ast_k\geq T_k$ and get that
\begin{align*}
 \f1 2 \leq1+\ep a(t,x)\leq \f3 2, \ \ \forall (t,x)\in [0,\infty)\times \mathbb{R}^d,
\end{align*}
\begin{align}\label{est49}
 E^\ep[\vae,\vue,\vthe](0,\infty)&\leq \sum\limits_{1\leq k\leq T_{N_0}} E^\ep[\vae,\vue,\vthe](T_{k-1},T_{k})\nonumber\\
 &\leq \sum\limits_{1\leq k\leq T_{N_0}}(4C_0)^k\|(a_0,\vu_0,\vartheta_0)\|_{X^\ep_2}\leq C\|(a_0,\vu_0,\vartheta_0)\|_{X^\ep_2},
\end{align}

\begin{align}\label{est50}
  M^{\ep,\alpha}_{p,q}[\vae,\vue-\vv,\vartheta^\ep;\tau^\ep,\sigma^\ep-\Theta](0,\infty)&\leq \sum\limits_{1\leq k\leq T_{N_0}}M^{\ep,\alpha}_{p,q}[\vae,\vue-\vv,\vartheta^\ep;\tau^\ep,\sigma^\ep-\Theta](T_{k-1},T)\nonumber\\
  &\leq \sum\limits_{1\leq k\leq T_{N_0}} \f{\de_0}{(4C_0)^{N_0-k}}\leq 2\delta_0,
\end{align}
\begin{align}\label{est51}
  &M^{\ep,\alpha}_{p,q}[\vae,\vue,\vartheta^\ep;\tau^\ep,\sigma^\ep](0,\infty)\\
  &\leq  M^{\ep,\alpha}_{p,q}[\vae,\vue-\vv,\vartheta^\ep;\tau^\ep,\sigma^\ep-\Theta](0,\infty)+2\|(\vv,\Theta)\|_{L^q\left(0,\infty;\dot{B}^{\f d p-1+\f2 q}_{p,1}\right)\cap L^1\left(0,\infty;\dot{B}^{\f d p+1}_{p,1}\right)}\leq C,\nonumber
\end{align}
\begin{align}\label{est52}
 \|\mathbb{P}\vue-\vv\|_{L^\infty\left(0,\infty;\dot{B}^{\f d p-1}_{p,1}\right)\cap L^1\left(0,\infty;\dot{B}^{\f d p+1}_{p,1}\right)}&\leq \sum\limits_{1\leq k\leq T_{N_0}}
   \|\mathbb{P}\vue-\vv\|_{L^\infty\left(T_{k-1},T_{k};\dot{B}^{\f d p-1}_{p,1}\right)\cap L^1\left(T_{k-1},T_{k};\dot{B}^{\f d p+1}_{p,1}\right)}\nonumber\\
  &\leq \f1 4\sum\limits_{1\leq k\leq T_{N_0}}\f{\de_0}{(4C_0)^{N_0-k}}\leq \f1 2 \delta_0.
\end{align}
\begin{align}\label{est63}
  &\|\mathbb{P}\vue-\vv\|_{L^\infty\Big(0,\infty;\dot{B}^{\f d p-1-(\f1 2-\f1 p)}_{p,1}\Big)\cap L^1\Big(0,\infty;\dot{B}^{\f d p+1-(\f1 2-\f1 p)}_{p,1}\Big)}+\|(\tau^\ep,\mathbb{Q}\vue)\|^{l;\f{\beta_0}\ep}_{L^2\Big(0,\infty;\dot{B}^{\f d p-(\f1 2-\f1 p)}_{p,1}\Big)}\nonumber\\
 &+\|\sigma^\ep-\Theta\|^{l;\f{\beta_0}\ep}_{L^\infty\Big(0,\infty;\dot{B}^{\f d p-1-(\f1 2-\f1 p)}_{p,1}\Big)\cap L^2\Big(0,\infty;\dot{B}^{\f d p-(\f1 2-\f1 p)}_{p,1}\Big)}\leq C\ep^{\f1 2-\f1 p}
\end{align}
From \eqref{est49}, we infer that $(\vae,\vue,\vthe)$ is in the class \eqref{cond7} and satisfies \eqref{est1}. Using \eqref{est51}, we obtain that
\begin{align}\label{new3}
 \|\mathbb{Q}\vue\|^{h;\f{\beta_0}\ep}_{L^2\Big(0,\infty;\dot{B}^{\f d p-(\f1 2-\f1 p)}_{p,1}\Big)}&\leq C\ep^{\f1 2-\f1 p}\|\mathbb{Q}\vue\|^{h;\f{\beta_0}\ep}_{L^2\Big(0,\infty;\dot{B}^{\f d p}_{p,1}\Big)}\nonumber\\
 &\leq C\ep^{\f1 2-\f1 p}M^{\ep,\alpha}_{p,q}[\vae,\vue,\vartheta^\ep;\tau^\ep,\sigma^\ep](0,\infty)\leq C\ep^{\f1 2-\f1 p},
\end{align}
\begin{align}\label{est53}
&\|\sigma^\ep-\Theta\|^{h;\f{\beta_0}\ep}_{L^\infty\Big(0,\infty;\dot{B}^{\f d p-2}_{p,1}\Big)\cap L^1\Big(0,\infty;\dot{B}^{\f d p}_{p,1}\Big)}\nonumber\\
&\leq C\ep \left(\ep\|\vae\|^{h;\f{\beta_0}\ep}_{L^\infty\Big(0,\infty;\dot{B}^{\f d p}_{p,1}\Big)}+\f1 \ep\|\vae\|^{h;\f{\beta_0}\ep}_{L^1\Big(0,\infty;\dot{B}^{\f d p}_{p,1}\Big)}\right.\nonumber\\
&\left.+\f1 \ep\|\vthe\|^{h;\f{\beta_0}\ep}_{L^\infty\Big(0,\infty;\dot{B}^{\f d p-2}_{p,1}\Big)\cap L^1\Big(0,\infty;\dot{B}^{\f d p}_{p,1}\Big)}+\|\Theta\|^{h;\f{\beta_0}\ep}_{L^\infty\Big(0,\infty;\dot{B}^{\f d p-1}_{p,1}\Big)\cap L^1\Big(0,\infty;\dot{B}^{\f d p+1}_{p,1}\Big)}\right)\nonumber\\
&\leq C\ep\left(M^{\ep,\alpha}_{p,q}[\vae,\vue,\vartheta^\ep;\tau^\ep,\sigma^\ep](0,\infty)+\|\Theta\|^{h;\f{\beta_0}\ep}_{L^\infty\Big(0,\infty;\dot{B}^{\f d p-1}_{p,1}\Big)\cap L^1\Big(0,\infty;\dot{B}^{\f d p+1}_{p,1}\Big)}\right)\nonumber\\
&\leq C\ep.
\end{align}
Similarly, we have
\begin{align}\label{new2}
\|\tau^\ep\|^{h;\f{\beta_0}\ep}_{L^\infty\Big(0,\infty;\dot{B}^{\f d p-2}_{p,1}\Big)\cap L^1\Big(0,\infty;\dot{B}^{\f d p}_{p,1}\Big)}\leq C\ep.
\end{align}
From \eqref{est63}-\eqref{new2}, we obtain \eqref{est3}.

Finally we point out that \eqref{est50} and \eqref{est52} hold for any $\delta\leq\delta_0$ and sufficient small $\ep$. More specifically, for any $\delta\leq\delta_0$, from Lemma 2.10 in \cite{F}, there exist a positive number $N_\delta$ and a time sequence $\{T_{\delta;n}\}^{N_\delta}_{n=0}$ such that
$0=T_{\delta;0}<T_{\delta;1}<\cdots<T_{\delta;N_\delta-1}<T_{\delta;N_\delta}=\infty$ and
\begin{align*}
 \|(\vv,\Theta)\|_{L^q\left(T_{\delta;n-1},T_{\delta;n};\dot{B}^{\f d p-1+\f2 q}_{p,1}\right)\cap L^1\left(T_{\delta;n-1},T_{\delta;n};\dot{B}^{\f d p+1}_{p,1}\right)}\leq\de,\ \ n=1,2,\cdots,N_\delta.
\end{align*}
Then we can select $\alpha_\delta$ such that for $\ep\leq1$
\begin{align*}
&\|(\Theta,\vv)\|^{h;\f{\alpha_\delta} {64}}_{L^2(I;\dot{B}_{p,1}^{\f d p})\cap L^{q'}(I;\dot{B}_{p,1}^{\f d p-1+\f2 {q'}})}\leq\f{\de}{(4C_0)^{N_\delta}},
 &\|(a_0,\mathbb{Q}\vu_0,\vartheta_0)\|^{h;\alpha_\delta}_{X^\ep_p}\leq\f{\de}{(4C_0)^{N_\delta}}.
\end{align*}
For such $\alpha_\delta$, we also select  $\ep_\delta\leq\min\{\ep_0,\ \ \f{1}{\alpha_\delta},\ \ \f{\beta_0}{16\alpha_\delta}\}$ such that
\[
(\alpha_\delta\ep)^{\f{1}{q(q-1)}}\Big(\|(a_0,\vu_0,\vartheta_0)\|_{X^\ep_2}+\delta\Big)\leq \f1 {36}\f{\de}{(4C_0)^{N_\delta}}
\]
for any $0<\ep\leq\ep_\delta$.

In a completely similar way to the previous argument, we can get a similar result to \eqref{est50} and \eqref{est52}, i.e.,
\begin{align}
  &M^{\ep,\alpha_\delta}_{p,q}[\vae,\vue-\vv,\vartheta^\ep;\tau^\ep,\sigma^\ep-\Theta](0,\infty)\leq 2\delta,\nonumber\\
 &\|\mathbb{P}\vue-\vv\|_{L^\infty\left(0,\infty;\dot{B}^{\f d p-1}_{p,1}\right)\cap L^1\left(0,\infty;\dot{B}^{\f d p+1}_{p,1}\right)}\leq \f1 2 \delta,\label{est65}
\end{align}
for any $0<\ep\leq\ep_\delta$.
Note that
\begin{align}\label{est64}
  &\ep\|a\|_{L^\infty(0,\infty;\dot{B}^{\f d p}_{p,1})}+\|\mathbb{Q}\vue\|_{L^q(0,\infty;\dot{B}^{\f d p-1+\f2 q}_{p,1})\cap L^{q'}(0,\infty;\dot{B}^{\f d p-1+\f2 {q'}}_{p,1})}\nonumber\\
  &\leq CM^{\ep,\alpha_\delta}_{p,q}[\vae,\vue-\vv,\vartheta^\ep;\tau^\ep,\sigma^\ep-\Theta](0,\infty)\nonumber\\
  &\leq C\delta
\end{align}
where the constant $C$ depends on $d,p,q$. From \eqref{est65} and \eqref{est64}, we obtain \eqref{est66}.\ \ $\Box$
\vspace{5mm}

\centerline{\bf Acknowledgements}
\vspace{2mm}
The author would like to thank Professor Yongzhong Sun and Professor Yang Li for many helpful discussions.

\vspace{5mm}

\centerline{\bf Conflict of interest}
\vspace{2mm}
the  author states that there is no conflict of interest.

\vspace{5mm}

\centerline{\bf Data Availability Statement}
\vspace{2mm}
Data sharing not applicable to this article as no datasets were generated or
analyzed during the current study.


\end{document}